\documentclass[12pt]{article}
\usepackage{amsfonts}
\usepackage{epsfig}
\title{The Octagonal PET II: The Topology of the Limit Sets}
\author{Richard Evan Schwartz \thanks{\hskip 5 pt Supported by 
N.S.F. Research Grant DMS-0072607}}

\newtheorem{theorem}{Theorem}[section]

\newtheorem{lemma}[theorem]{Lemma}

\newtheorem{corollary}[theorem]{Corollary}

\def\startproof{{\bf {\medskip}{\noindent}Proof: }}

\def\endproof{$\spadesuit$  \newline}

\def\R{\mbox{\boldmath{$R$}}}% 

\begin{document}
\maketitle

\begin{abstract}
This paper is a sequel to [{\bf S0\/}].  In
[{\bf S0\/}], we studied 
a $1$-parameter family of polygon
exchange transformations.  We showed
that this family is completely renormalizable.
In this paper, we use the renormalization
scheme in [{\bf S0\/}] to give a complete
classification of topological types
of the limit sets which arise in
these PETs.
\end{abstract}

\section{Introduction}

\subsection{Background}

A {\it polytope exchange transformation\/} (or PET) is
defined by a polytope $X$ which has been partitioned
in two ways into smaller polytopes:
$$X=\bigcup_{i=1}^m A_i=\bigcup_{i=1}^m B_i.$$
For each $i$ there is some vector $V_i$ such that
$B_i=A_i+V_i$.  That is, some translation carries
$A_i$ to $B_i$.  One then defines a map $f: X \to X$
by the formula $f(x)=x+V_i$ for all $x \in {\rm int\/}(A_i)$.
It is understood that $f$ is not defined for points
in the boundaries of the small polytopes.  The
inverse map is defined by $f^{-1}(y)=y-V_i$ for all
$y \in {\rm int\/}(B_i)$.

The simplest examples of PETs are $1$-dimensional
systems, known as {\it interval exchange transformations\/}
(or IETs).
These systems have been extensively studied
in the past $30$ years, and there are close connections
between IETs and other areas of mathematics such
as Teichmuller theory.
See, for instance, [{\bf Y\/}] and [{\bf Z\/}] and the many references
mentioned therein.

Some examples of polygon exchange maps have been
studied in [{\bf AG\/}], [{\bf AKT\/}], [{\bf H\/}],
[{\bf Hoo\/}],
[{\bf LKV\/}],
[{\bf Low\/}], [{\bf S2\/}],
[{\bf S3\/}], and [{\bf T\/}].
Some definitive theoretical work concerning the
(zero) entropy of such systems is done in
[{\bf GH1\/}], [{\bf GH2\/}], and [{\bf B\/}].

The {\it Rauzy renormalization\/} [{\bf R\/}]
gives a satisfying renormalization theory for
the family of IETs all having the same number of
intervals in the partition.
Often, renormalization phenomena are found
in individual PETs, and in the papers
[{\bf S0\/}], [{\bf S1\/}] and [{\bf Hoo\/}],
there were constructed families of PETs
which had a renormalization theory at least
vaguely similar to the Rauzy renormalization.

In [{\bf S0\/}], we introduced a family of PETs 
in every even dimension.
In dimension $2n$, the objects
are indexed by $GL_n(\R)$.  In the
two-dimensional setting, there is a $1$-parameter
family. We studied this family in detail,
and in particular worked out a renormalization
scheme for the family. 
In this paper, we use the renormalization scheme
to work out the topology of the associated limit sets
of the $2$-dimensional examples.

In [{\bf S0\/}] we give a somewhat fuller discussion
of the background and context for our examples.

\subsection{Construction of the PET}

Here we recall the basic construction given in [{\bf S0\/}],
at least in the $2$-dimensional case.
Our construction depends on a parameter $s \in (0,\infty)$.
The $8$ parallalograms in Figure 1.1 
are the orbit of a single parallelogram $P$ under a dihedral group
of order $8$.  Two of the sides of $P$ are
determined by the vectors
$(2,0)$ and $(2s,2s)$. For $j=1,2$, let $F_j$ denote the parallelogram
centered at the origin and translation equivalent to the ones
in the picture labeled $F_j$.  Let $L_j$ denote the lattice
generated by the sides of the parallelograms labeled $L_j$.
(Either one generates the same lattice.)

\begin{center}
\resizebox{!}{3in}{\includegraphics{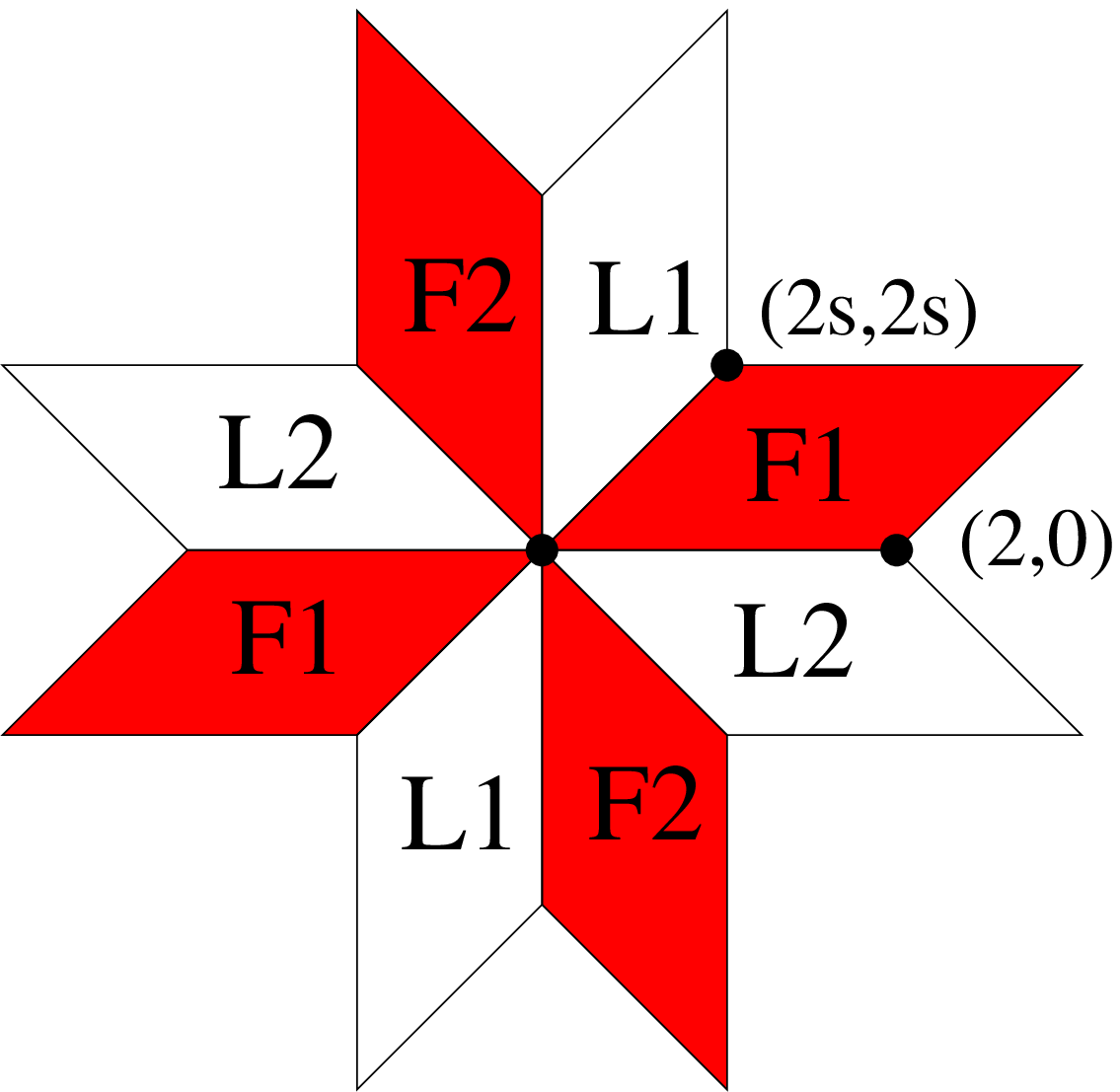}}
\newline
{\bf Figure 1.1:\/} The scheme for the PET.
\end{center}

It turns out that
$F_i$ is a fundamental domain for
$L_j$, for all $i,j \in \{1,2\}$.
Given $p \in F_j$ we let
\begin{equation}
f'(p)=p+V_p \in F_{3-j}, \hskip 30 pt
V_p \in L_{3-j}.
\end{equation}
The choice of $V_p$ is almost always
unique, on account of $F_{3-j}$ being
a fundamental domain for $L_{3-j}$.
When the choice is not unique,
we leave $f'$ undefined.
When $p \in F_1 \cap F_2$ we have
$V_p=0 \in L_1 \cap L_2$.
We define $f=(f')^2$, which preserves both
$F_1$ and $F_2$, and we set $X=F_1$.
Our system is
$f: X \to X$, which we denote by $(X,f)$.

\subsection{Prior Results}

The results in
[{\bf S0\/}] (and in this paper)
ultimately follow from a renormalization
scheme associated to our family.  Define the 
{\it renormalization map\/} $R: (0,1) \to [0,1)$ as
follows.
\begin{itemize}
\item $R(x)=1/(2x)-{\rm floor\/}(1/(2x))$ if $x<1/2$.
\item $R(x)=1-x$ if $x>1/2$.
\end{itemize}
This map will appear in many of our results.
We explain its direct connection to the PETs
in \S 2.

A {\it periodic tile\/} for $(X,f)$
as a maximal
convex polygon on which $f$ and its iterates
are completely defined and periodic.
We call the union $\Delta$ of the periodic
tiles the {\it tiling\/}. 
In [{\bf S0\/}] we proved, among other things,
the following results about the tiling.

\begin{theorem}
\label{tiling}
When $s$ is irrational, $\Delta_s$ is a
full measure set consisting entirely of
squares and semi-regular octagons.
\begin{itemize}
\item $\Delta_s$ consists entirely of
squares iff $R^n(s)<1/2$ for all $n$.
\item $\Delta_s$ has finitely many octagons
iff $R^n(s)>1/2$ for finitely many $n$.
\item $\Delta_s$ has infinitely many octagons
iff $R^n(s)>1/2$ for infinitely many $n$.
\end{itemize}
\end{theorem}

The fact that the union of tiles in $\Delta_s$
has full measure allows us to define the limit set
\footnote{For a general polygon exchange map, the
limit set is probably best defined as the set of
weakly aperiodic points. We call a point
$p \in X$ {\it weakly aperiodic\/} if
there is a sequence
$\{q_n\}$ converging to $p$ with the
following property. The first iterates
of $f$ are defined on $q_n$ and
the points $f^k(q_n)$ for $k=1,...,n$
are distinct.} $\widehat \Lambda_s$ to be the
set of points $p$
such that
every neighborhood of $p$ intersects
infinitely many tiles of $\Delta$.  

In [{\bf S0\/}] we proved, among other things,
the following results about the limit set.
\begin{theorem}
Let $s \in \R$ be irrational.
The projection of $\widehat \Lambda_s$
onto any line parallel to an $8$th root
of unity contains a line segment.  Hence
$\widehat \Lambda_s$ has Hausdorff
dimension at least $1$.  Moreover,
$\widehat \Lambda_s$ is 
not contained in a finite union of lines.
\end{theorem}

We also proved that our $1$-parameter
family had what one might call hidden
hyperbolic symmetry.

\begin{theorem}
The Hausdorff dimension of
$\widehat \Lambda_s$, as a
function of the parameter, is invariant
under the action of the $(2,4,\infty)$
hyperbolic reflection triangle group
acting on the parameter space by linear
fractional transformations.
\end{theorem}

The results in [{\bf S0\/}] are somewhat
more detailed than the ones listed above,
and the reader should refer to [{\bf S0\/}]
for the more detailed versions of the results above.

\subsection{Topology of the Limit Sets}

Here is the main result of this paper.
\begin{theorem}[Main]
Let $s \in (0,1)$ be irrational.
\begin{enumerate}
\item $\widehat \Lambda_s$ is a disjoint union of two
arcs if and only if $R^n(s)<1/2$ for all $n$.
\item $\widehat \Lambda_s$ is a finite forest \footnote{A 
{\it finite forest\/} is a finite disjoint union
of finite trees}
if and only if $R^n(s)>1/2$ for finitely many $n$.
\item $\widehat \Lambda_s$ is a Cantor set if and only
if $R^n(s)>1/2$ infinitely often.
\end{enumerate}
\end{theorem}

Combining the Main Theorem with
Theorem \ref{tiling}, we get the
following corollary.

\begin{corollary}
Let $s \in (0,1)$ be irrational.
\begin{enumerate}
\item $\widehat \Lambda_s$ is a disjoint union of two
arcs if and only if $\Delta_s$ contains only squares.
\item $\widehat \Lambda_s$ is a finite forest if and only
if $\Delta_s$ contains finitely many octagons.
\item $\widehat \Lambda_s$ is a Cantor set if and only if
$\Delta_s$ contains infinitely many octagons.
\end{enumerate}
\end{corollary}

Figures 1.2 shows Statement 1 of the Main Theorem in
action.

\begin{center}
\resizebox{!}{2.4in}{\includegraphics{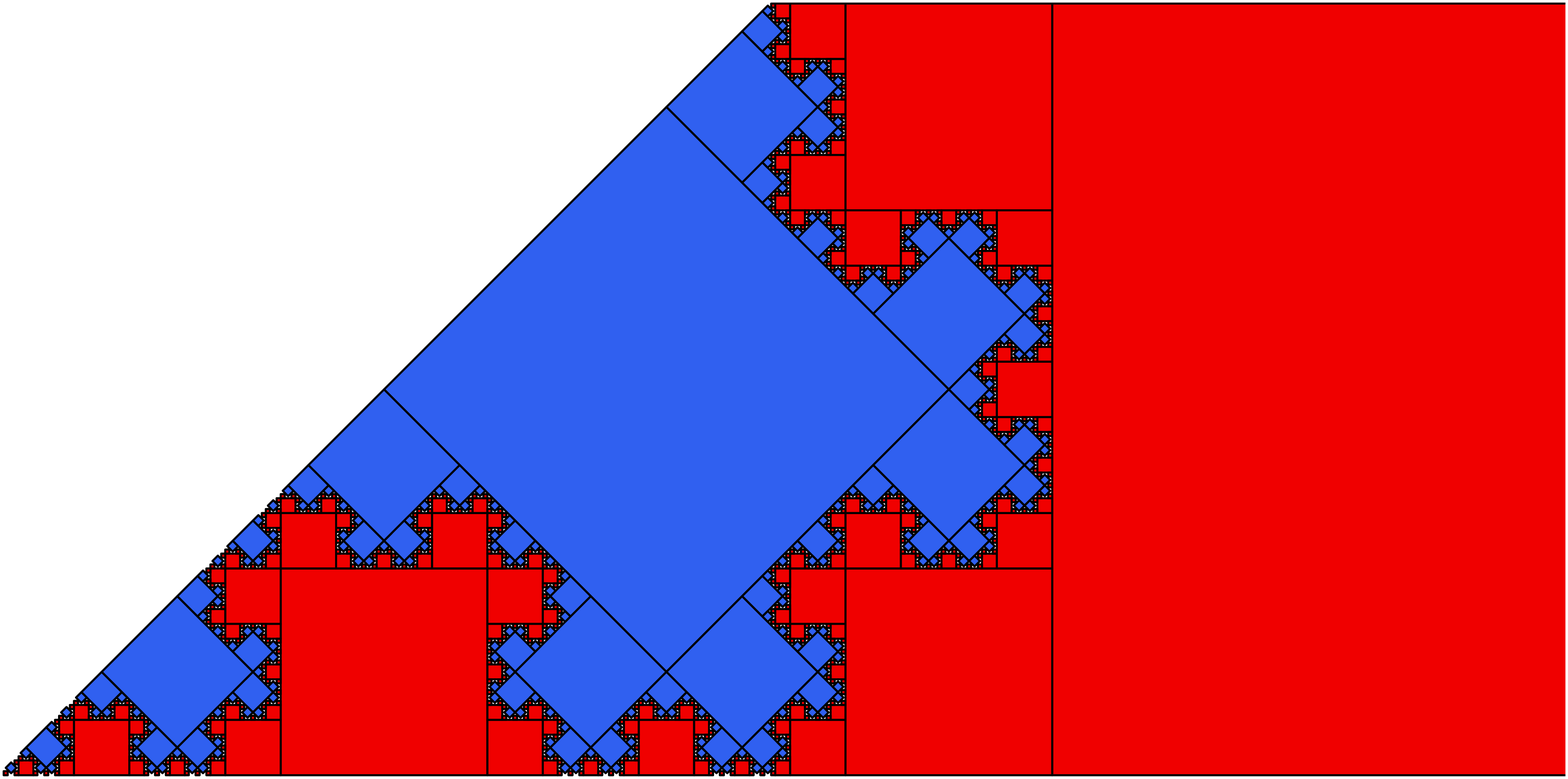}}
\newline
{\bf Figure 1.2:\/} The left half of $\Delta_s$ for
$s=\sqrt 3/2-1/2$.
\end{center}

Figure 1.3 shows Statement 2 of the Main Theorem in
action.  The parameter is such that its even expansion
is $3,1,3,1,2,2,...$.  A $k$ in the $n$th position
indicates that $R^n(s) \in (1/(k+1),1/k)$.

\begin{center}
\resizebox{!}{2.4in}{\includegraphics{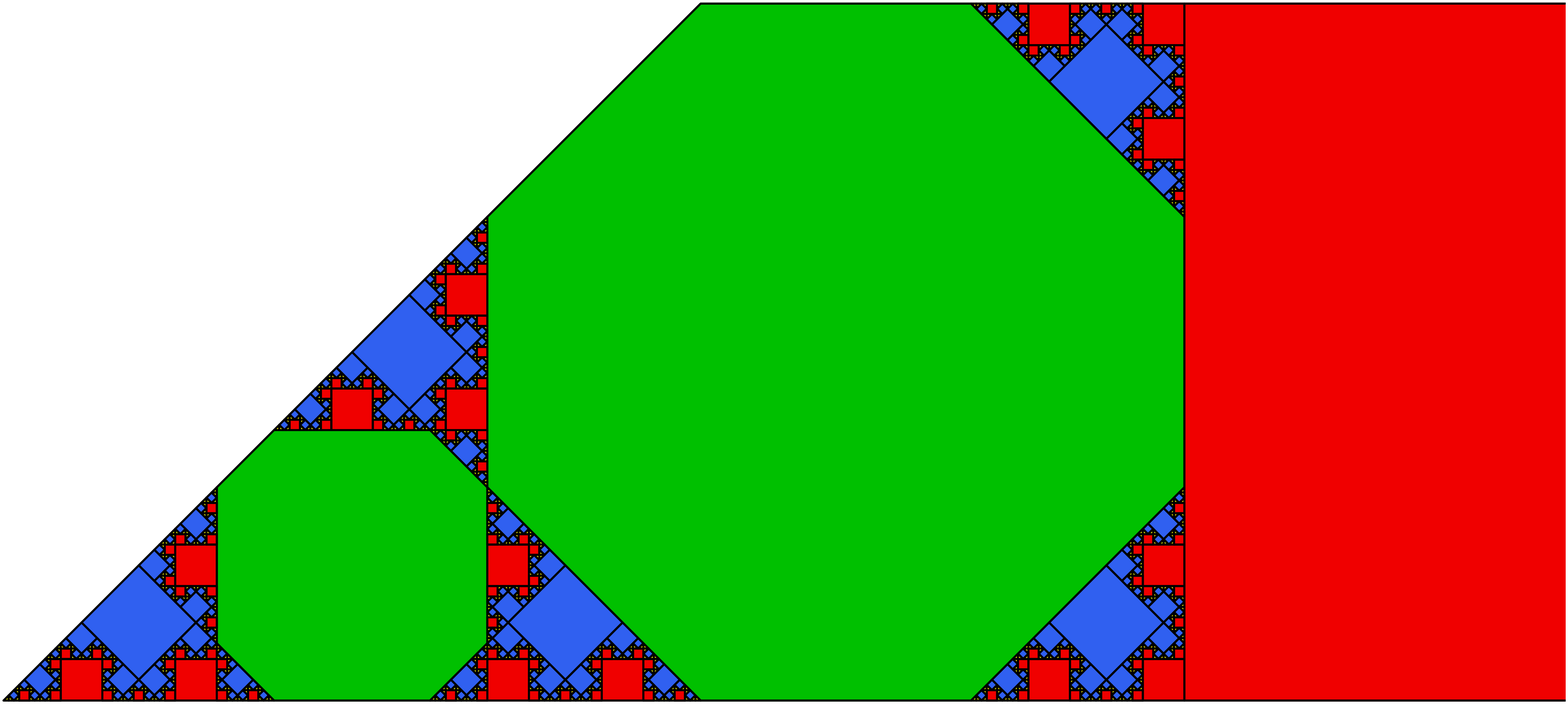}}
\newline
{\bf Figure 1.3:\/} The left half of $\Delta_s$ for
$s=(3,1,3,1,2,2,2...)$.
\end{center}

Figure 1.4 shows Statement 3 of the Main Theorem in

\begin{center}
\resizebox{!}{2.4in}{\includegraphics{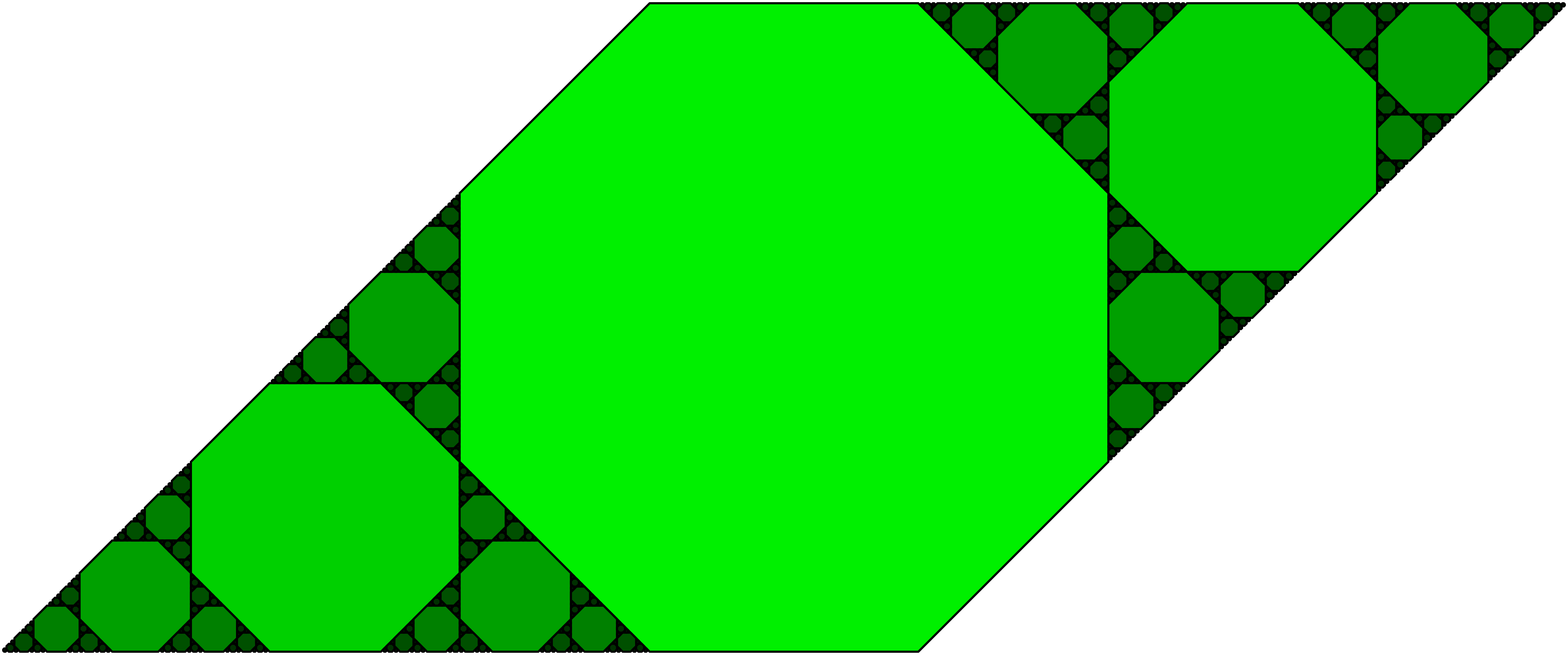}}
\newline
{\bf Figure 1.3:\/} The left half of $\Delta_s$ for
$s=\sqrt 2/2$.
\end{center}

In all these figures,
the right half is a rotated image
of the left half.
Intuitively, the squares tend to line up in simply connected
chunks whose complementary regions are curves and the
octagons both split and splice the curves.

This paper is organized as follows. 
In \S 2 we will
recall the renormalization
scheme from [{\bf S0\/}], as well as
some symmetry results.   Following
the summary given in \S 2, the
rest of this paper is self-contained.
In \S 3-4 we prove some geometric
and combinatorial results about how the
tiles of $\Delta_s$ fill $X_s$ and
how $\widehat \Lambda_s$ intersects
various lines of symmetry.
In \S 5-7 we prove the Main Theorem, one
statement per chapter.

We would like to say a word about the nature of our
proofs.  For the sake of giving a readable exposition,
there are a number of routine geometric calculations,
concerning polygonal regions in the plane, which we omit.
Such calculations are all exercises in plane geometry.
Rather than write out these calculations,
we illustrate them with extensive pictures
from our java program OctaPET.   It is certainly
best to read this paper while using OctaPET; this is
how we wrote the paper. 

A common mistake made by beginning students is to
try to prove a general statement by just considering
one example.  We do not mean to make this mistake
here, even though superficially some of our proofs
look like this. In a written paper dealing
with a $1$-parameter family of systems, we cannot
illustrate the picture for every parameter.  The pictures
we do show are typical for the given parameter interval,
and the written arguments we give only make statements
which hold for all the relevant parameters.
\newline

I would like to thank Nicolas Bedaride,
Pat Hooper, Injee Jeong,
John Smillie, and
Sergei Tabachnikov for interesting
conversations about topics related to
this work.  Some of this work
was carried out at ICERM in Summer 2012,
and some was
carried out during my sabbatical in
2012-13.
This sabbatical was funded from many sources.
I would like to thank the National
Science Foundation, All Souls College, Oxford,
the Oxford Maths Institute,
the Simons Foundation, the Leverhulme Trust, and
Brown University for their support during this time
period.

\newpage
\section{Preliminaries}
\subsection{Closedness}

We start with an essentially obvious result, 
which we state for the record.

\begin{lemma}
$\widehat \Lambda_s$ is closed.
\end{lemma}

\startproof
Let $\{p_n\}$ be a sequence of points in
$\widehat \Lambda_s$ converging to
some $q$.  We want to show that 
$q \in \widehat \Lambda_s$.  We need to
show that every open neighborhood $U$ of
$q$ contains infinitely many tiles of
$\Delta_s$.  For some $n$ we have
$p_n \in U$.  But then some neighborhood
$V$ of $p_n$ lies in $U$. But $V$ contains
infinitely many tiles of $\Delta_s$.
Hence, so does $U$.
\endproof

\subsection{The Even Expansion}

We call an irrational number
$s \in (0,1)$ is {\it oddly even\/}
if it has continued fraction expansion
$[0,a_1,a_2,a_3,...]$ with $a_k$ even
for all odd $k$.  In [{\bf S0\/}] we
showed that $s$ is oddly even if and
only if $R^n(s)<1/2$ for all $n$.
We call the rational $s$ {\it oddly even\/} if
the (finite) orbit $\{R^n(s)\}$ avoids
$(1/2,1)$.   In [{\bf S0\/}], we proved
the following easy-to-believe lemma.

\begin{lemma}
\label{appx}
Let $s$ be an oddly even irrational.
There is a sequence $\{r_n\}$ of oddly even rational
numbers which converges to $s$.
\end{lemma}

When $s$ is irrational, we define
the {\it even expansion\/} of $s$ to
be $n_0,n_1,n_2,...$ where
\begin{equation}
\label{interval}
\frac{1}{n_k+1}<R^k(s)<\frac{1}{n_k}.
\end{equation}
Since $R$ maps each interval $(1/2n,1/(2n-1)]$ onto
$[1/2,1)$, an irrational $s \in (0,1)$
is oddly even if and only if its even expansion
consists entirely of even numbers.

We can define the even expansion of $s \in (0,1)$
when $s$ is rational.  The last nonzero term in
the $R$-orbit of $s$ is $1/2q$ for
some $q=1,2,3...$ and by convention we take
the corresponding term in the even expansion to be $1/2q$.
The preceding terms lie in the interiors of the
intervals in Equation \ref{interval}.
Wnen $r_n \to s$ and $s$
is irrational,
even expansions of $r_n$ converges to the even expansion of $s$.
\subsection{Rational Limits}

When we speak of the convergence of compact sets, we
mean to use the Hausdorff topology.  A sequence
of compact sets $\{P_n\}$ converges to a compact
set $P_{\infty}$ if, for every $\epsilon>0$ there
is some $N$ such that $n>N$ implies that
the Hausdorff distance from $P_n$ to
$P_{\infty}$ is less than $\epsilon$.
The Hausdorff distance between compact sets
$S_1$ and $S_2$ is the infimal $d$ such that
$S_j$ is contained in the $d$-neighborhood
of $S_{3-j}$ for $j=1,2$.

Let $\{s_n\}$ be a sequence of rationals
converging to an irrational parameter $s \in (0,1)$.
Let $\Delta_n$ be the perodic tiling associated
to $s_n$ and let $\Delta_{\infty}$ be the periodic tiling
associated to $s$. 

\begin{lemma}[Approximation]
\label{converge}
Let $P_{\infty}$ be a tile of $\Delta_{\infty}$.
Then for all large $n$, there is a tile
$P_n$ of $\Delta_n$ such that $\{P_n\}$ converges
to $P_{\infty}$.
\end{lemma}

\startproof
We proved this in [{\bf S0\/}, \S 3,3]. 
Here is a sketch.
When $s$ is irrational, we show that the
periodic points are stable under perturbation.
Hence, there is some sequence $\{P_n\}$ such
that $\lim P_n$ contains an open subset of
$P_{\infty}$.
We also show that the set of periodic points
in $\Delta_t$ with the same combinatorics 
(essentially the list
of vectors we use when defining the map $f_n$ and
its iterates, written in a parameter-independent basis)
is a single tile of $\Delta_t$, as long as
$t \not = 1/n$ for $n=1,2,3...$.
That is, different
tiles have different associated dynamical combinatorics.
This uniqueness forces $\lim P_n=P_{\infty}$.
\endproof

\begin{lemma}
Let $P_n$ be a tile of $\Delta_n$
and $\{P_n\}$ converges to a polygon $P_{\infty}$.
Then $P_{\infty}$ is a tile of $\Delta_{\infty}$.
\end{lemma}

\startproof
We give the argument in [{\bf S0\/}, \S 6.2].
Here is a sketch.
There is a uniform lower bound
on the size of $P_n$ and a uniform upper bound on the
area of $X_n$. This puts a uniform upper bound on the
period of $P_n$. Passing to a subsequence, we can assume
that the combinatorics of the orbit
is independent of $n$.  By continuity, any point
in $P_{\infty}$ is a periodic point for $f_n$
when $n$ is sufficiently large.  This shows that
there is a periodic tile $P'_{\infty}$ of $\Delta_{\infty}$
which contains $P_{\infty}$. 
Applying the Convergence Lemma, there is some sequence
$\{P'_n\}$ such that $P'_n \to P'_{\infty}$.
Eventually $P'_n$ and $P_n$ overlap, and so we must
have $P_n=P'_n$. But then $P'_{\infty}=P_{\infty}$.
\endproof

We abbreviate the results in this section by saying that
the tilings $\{\Delta_n\}$ converge to the tiling
$\Delta_{\infty}$.

\subsection{Renormalization}

In this section we explain the significance
of the renormalization map $R$ which appears
in our results.

We use the notation from \S 1. So, $F_1=X$ is the
domain of our system and $F_2$ is the image of
$F_1$ rotated by $\pi/2$.  The map $f$ is
trivial on the intersection $F_1 \cap F_2$.
We call this intersection the {\it trivial tile\/}.

Suppressing the parameter $s$, we define
\begin{equation}
Y=F_1-F_2 = X-F_2 \subset X.
\end{equation}
$Y$ is the portion of $X$ outside the trivial tile.

We call a subset $S \subset X$ {\it clean\/} if
$\partial S$ does not intersect the interior
of any tile of $\Delta$.
Here is the main result from [{\bf S0\/}].

\begin{theorem}
\label{renorm}
Suppose $s \in (0,1)$ and $t=R(s) \in (0,1)$.
There is a clean set $Z_s \subset X_s$
and a similarity $\phi_s: Y_t \to Z_s$
such that $\phi_s(\Delta_t \cap Y_t)=\Delta_s \cap Z_s$.
\begin{enumerate}
\item $\phi_s$ commutes with reflection in the origin and
maps the acute vertices of $X_t$ to the
acute vertices of $X_s$.
\item  When $s<1/2$, the restriction
of $\phi_s$ to each component of $Y_t$
is an orientation reversing similarity,
with scale factor $s \sqrt 2$.
\item When $s<1/2$, either half of
$\phi_s$ extends to the trivial tile of
$\Delta_t$ and maps it to a tile in $\Delta_s$. 
\item When $s<1/2$, the only nontrivial
orbits which miss $Z_s$ are
contained in the 
$\phi_s$-images of the trivial tile of
$\Delta_t$.
These orbits have period $2$.
\item When $s>1/2$ the restriction
of $\phi_s$ to each component of $Y_s$
is a translation.
\item When $s>1/2$, all
nontrivial orbits intersect $Z_s$.
\end{enumerate}
\end{theorem}

Figures 2.1 and 2.2 illustrate the result for
$s<1/2$.  Figures 2.3 and 2.4 show the Main Theorem in
action for $s>1/2$.

\begin{center}
\resizebox{!}{1.3in}{\includegraphics{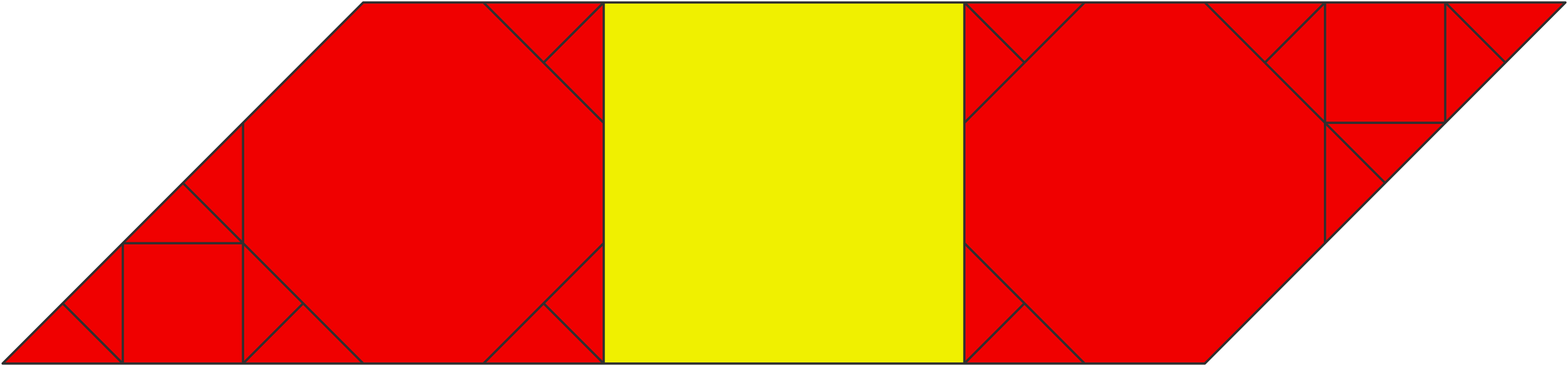}}
\newline
{\bf Figure 2.1:\/} $Y_t$ in red for $t=3/10=R(5/13)$.
\end{center}

\begin{center}
\resizebox{!}{1.4in}{\includegraphics{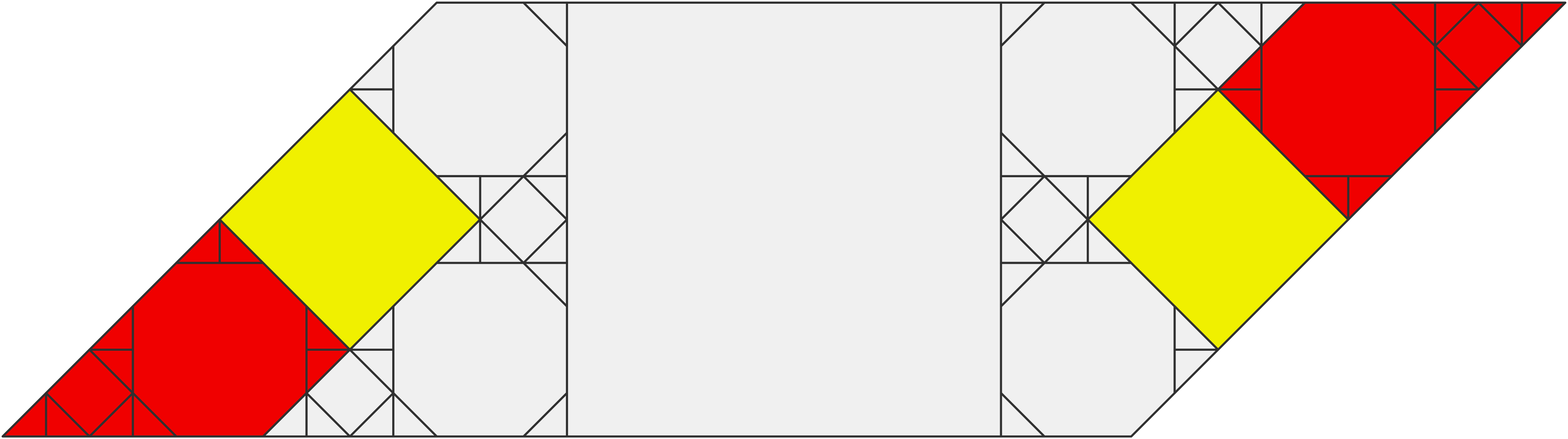}}
\newline
{\bf Figure 2.2:\/} $Z_s$ in red $s=5/13$.
\end{center}

\begin{center}
\resizebox{!}{1.35in}{\includegraphics{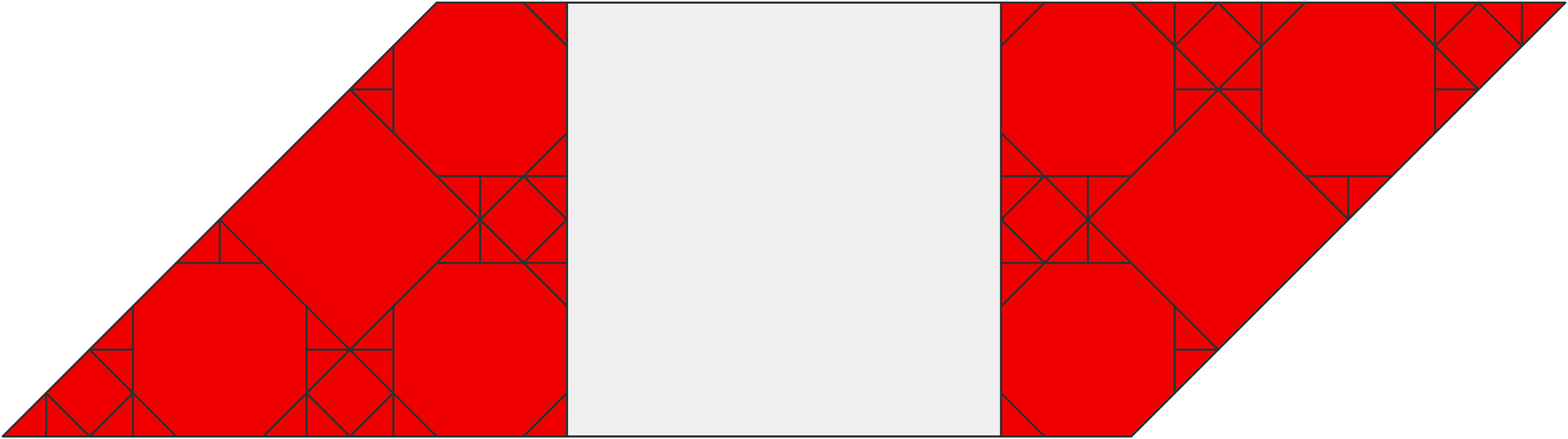}}
\newline
{\bf Figure 2.3:\/} $Y_t$ in red for $t=R(8/13)=5/13$.
\end{center}

\begin{center}
\resizebox{!}{2in}{\includegraphics{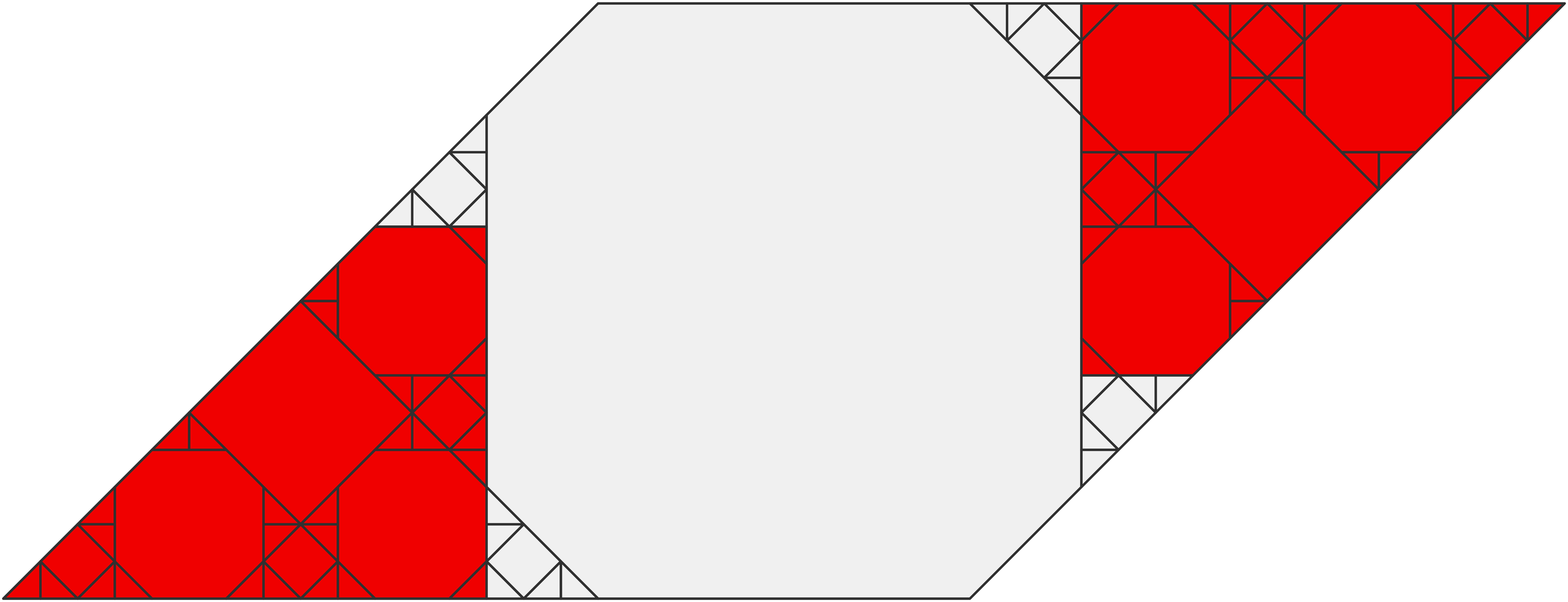}}
\newline
{\bf Figure 2.4:\/} $Z_s$ in red for $s=8/13$.
\end{center}
\subsection{Elementary Symmetry}

In [{\bf S0\/}] we established a number of symmetries
of the system.  The ones listed in this section have
easy proofs, though we do not give them here.
\newline
\newline
{\bf Rotation Symmetry:\/}
$\Delta_s$ is invariant under reflection in the origin:
\begin{equation}
\iota(x,y)=(-x,-y),
\end{equation}
In view of this symmetry, we will often
draw only the left half of the picture.
\newline
\newline
{\bf Inversion Symmetry:\/}
Let $t=1/(2s)$.  There is a similarity carrying
$\Delta_s$ to $\Delta_t$.  The similarity is
orientation reversing, fixes the origin, and
interchanges lines of slope $0$ with lines of
slope $1$.  We called this fact the
Inversion Lemma in [{\bf S0\/}]. The Inversion Lemma
gives shape to the renormalization map $R$.
\newline
\newline
{\bf Central Tiles:\/}
When $s \in (1/2,1)$, the intersection
$F_1 \cap F_2$ is an octagon, which we
call the {\it central tile\/}.
When $s \leq 1/2$ or $s \geq 1$, 
the intersection 
$F_1 \cap F_2$ is a square. This square
generates a grid in the plane, and
finitely many squares in this grid lie in
$X=F_1$.  We call these squares the
{\it central tiles\/}.  See Figure 2.5.

\begin{center}
\resizebox{!}{2.4in}{\includegraphics{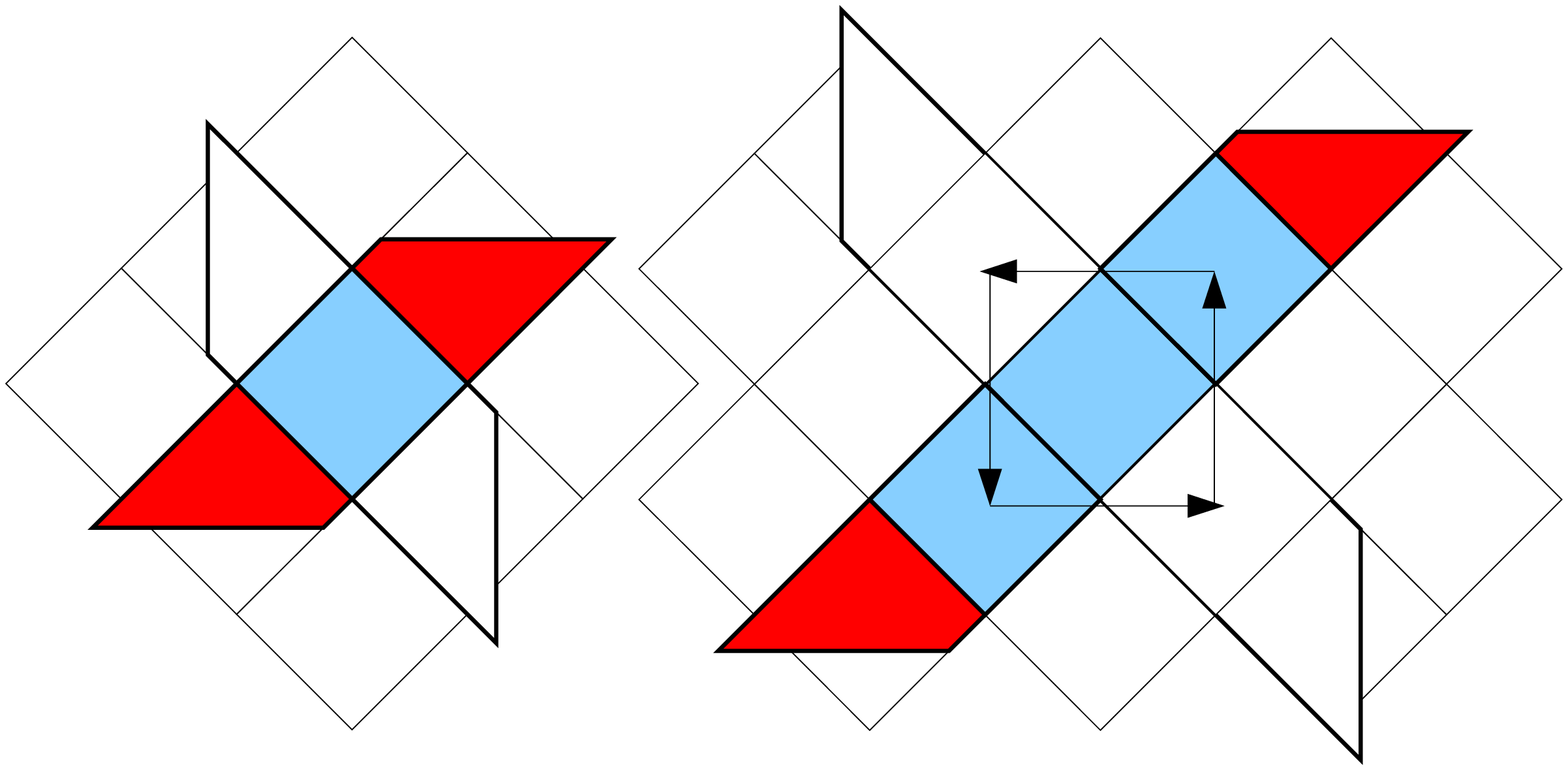}}
\newline
{\bf Figure 2.5:\/} 
The central tiles (blue) for $s=5/4$ and $t=9/4$.
\end{center}

For any relevant object $S \subset X_s$, the set
$S^0$ denotes the portion of $S$ lying to
the left of the central tiles of $\Delta_s$.
We will use this notational convention repeatedly.

\subsection{Insertion Symmetry}
\label{insertion}

When $s \geq 1$ and $t=s+1$, the
intersections 
$\Delta_s^0=\Delta_s \cap X_s^0$ and
$\Delta_t^0\Delta_t \cap X_s^0$ are isometric.
(We give the argument in [{\bf S0\/}].)
Combining this fact with the inversion symmetry,
we get the following result, also
stated in [{\bf S0\/}]. 

\begin{lemma}[Insertion]
If $s<1/2$
and $t=s/(2s+1)$ then 
$\Delta_s \cap X_s^0$ and
$\Delta_t \cap X_t^0$ are similar.
\end{lemma}

\begin{center}
\resizebox{!}{1.6in}{\includegraphics{Pix/renorm4.ps}}
\newline
{\bf Figure 2.6:\/} $\Delta_s$ for $s=5/13$.
\end{center}

\begin{center}
\resizebox{!}{1.4in}{\includegraphics{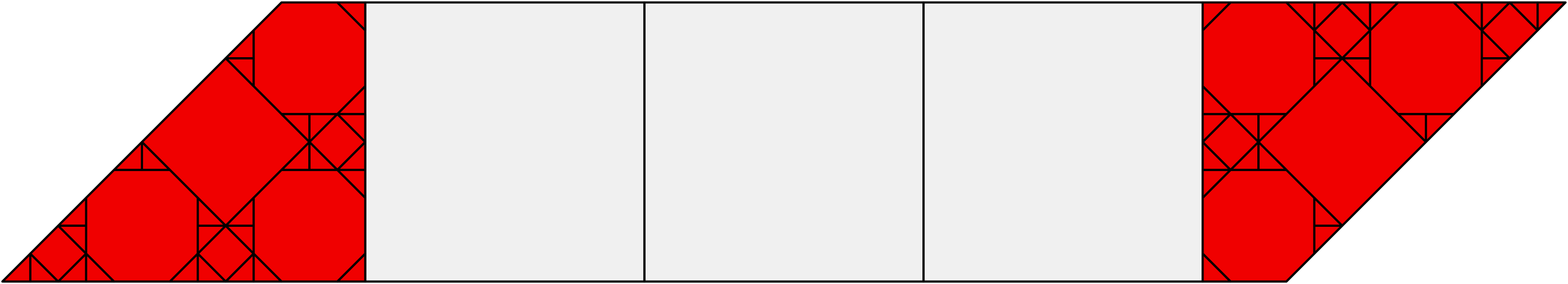}}
\newline
{\bf Figure 2.7:\/} 
$\Delta_t$ for $t=5/23$.
\end{center}

The reason for the name of this result
is that the picture for $t$ is obtained
from the picture for $s$ just by inserting
$2$ new central squares.  The reader can see
this in action by comparing Figures 2.6 and 2.7.

In view of the Insertion Lemma, we will often
describe our results for $s \in [1/4,1)$.
The one other advantage of the Insertion
Lemma, which we exploited in [{\bf S0\/}] is
that often a statement for parameters
in $[1/4,1)$ involves a finite computational
proof whereas the statement for all
parameters in $(0,1)$ would require
(if computed directly) an infinite computation.

\subsection{Bilateral Symmetry}
\label{symm}

The kind of symmetry described in this section
looks obvious from the pictures, but it required
a nontrivial computational proof in
[{\bf S0\/}]. We will describe the symmetry for
$s \in [1/4,1)$.  

We say that a line $L$ is a {\it line of symmetry\/}
for $\Delta_s$ if 
\begin{equation}
\Delta_s \cap \bigg(X_s \cap R_L(X_s)\bigg)
\end{equation}
is invariant under the reflection $R_L$ in $L$.
Note that $X_s$ itself need not be invariant
under $R_L$.  Consider the following lines.
\begin{itemize}
\item Let $H$ be the line $y=0$. 
\item Let $V$ be the line $x=-1$.
\item For $D_s$ be the line of slope $-1$ through $(-s,-s)$.
\item For $s \in [1/4,1/2]$, $E_s$ be the line of slope $-1$ $(-3s,-3s)$.
\item For $s \in [1/2,1]$, $E_s$ is the line of slope $1$ through $(-s,-s)$.
\end{itemize}
We call these the {\it fundamental lines of symmetry\/}.
For each line $L$ above, the line $\iota(L)$ is
also a line of symmetry.  However, we will not
usually refer to these other lines.

$H$, $V$, and $D_s$ and $E_s$ respectively are the lines of symmetry for
$$
A_s=(X_s \cap R_H(X_s))^0, \hskip 15 pt
B_s=X_s \cap R_V(X_s),$$
\begin{equation}
\label{hex}
P_s=X_s \cap R_D(X_s), \hskip 5 pt
Q_s=X_s \cap R_D(X_s).
\end{equation}

\begin{center}
\resizebox{!}{1.7in}{\includegraphics{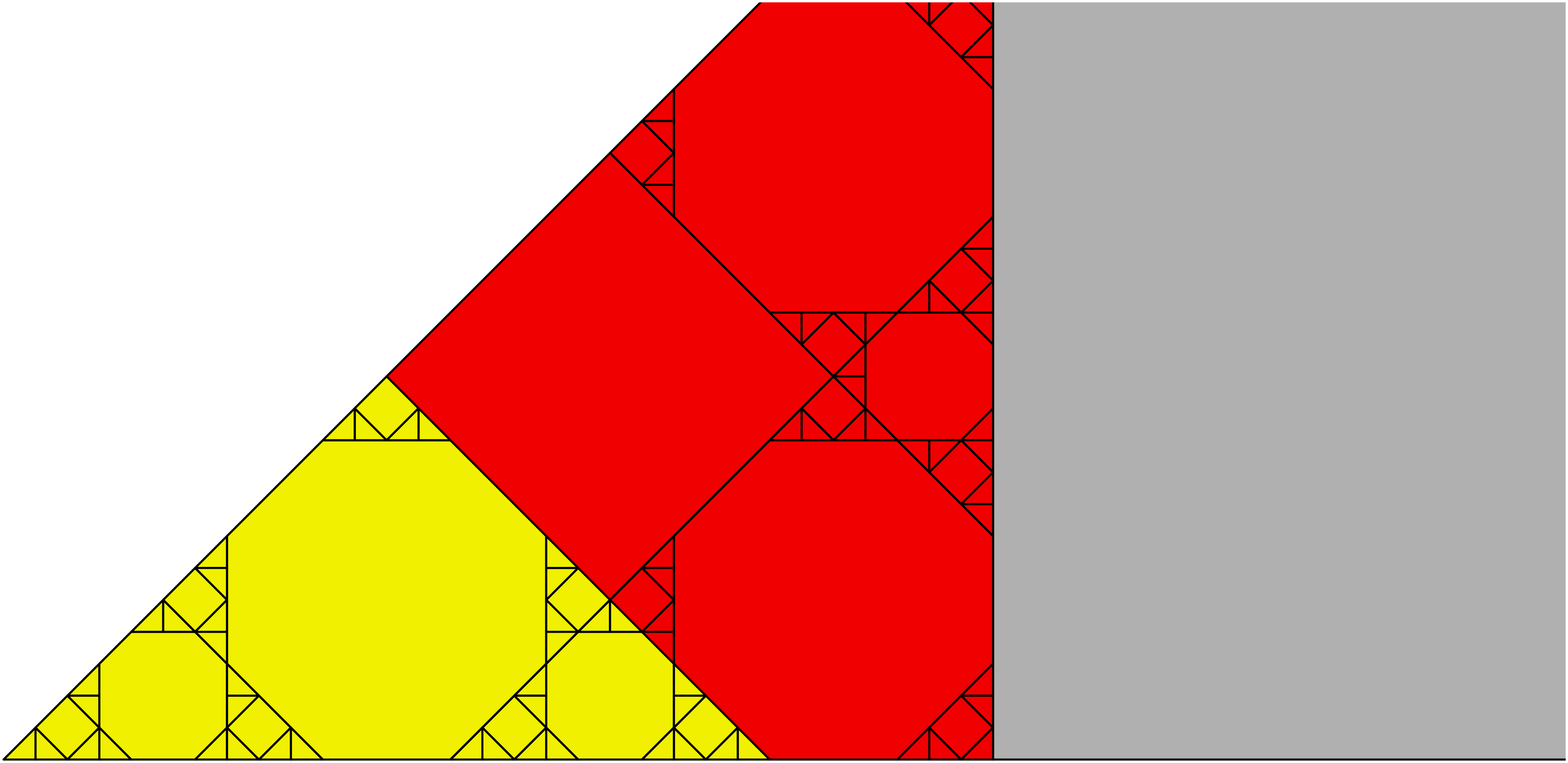}}
\newline
{\bf Figure 2.8:\/} 
$A_s$ (red) and $B_s$ (yellow) for $s=12/31$.
\end{center}

\begin{center}
\resizebox{!}{1.7in}{\includegraphics{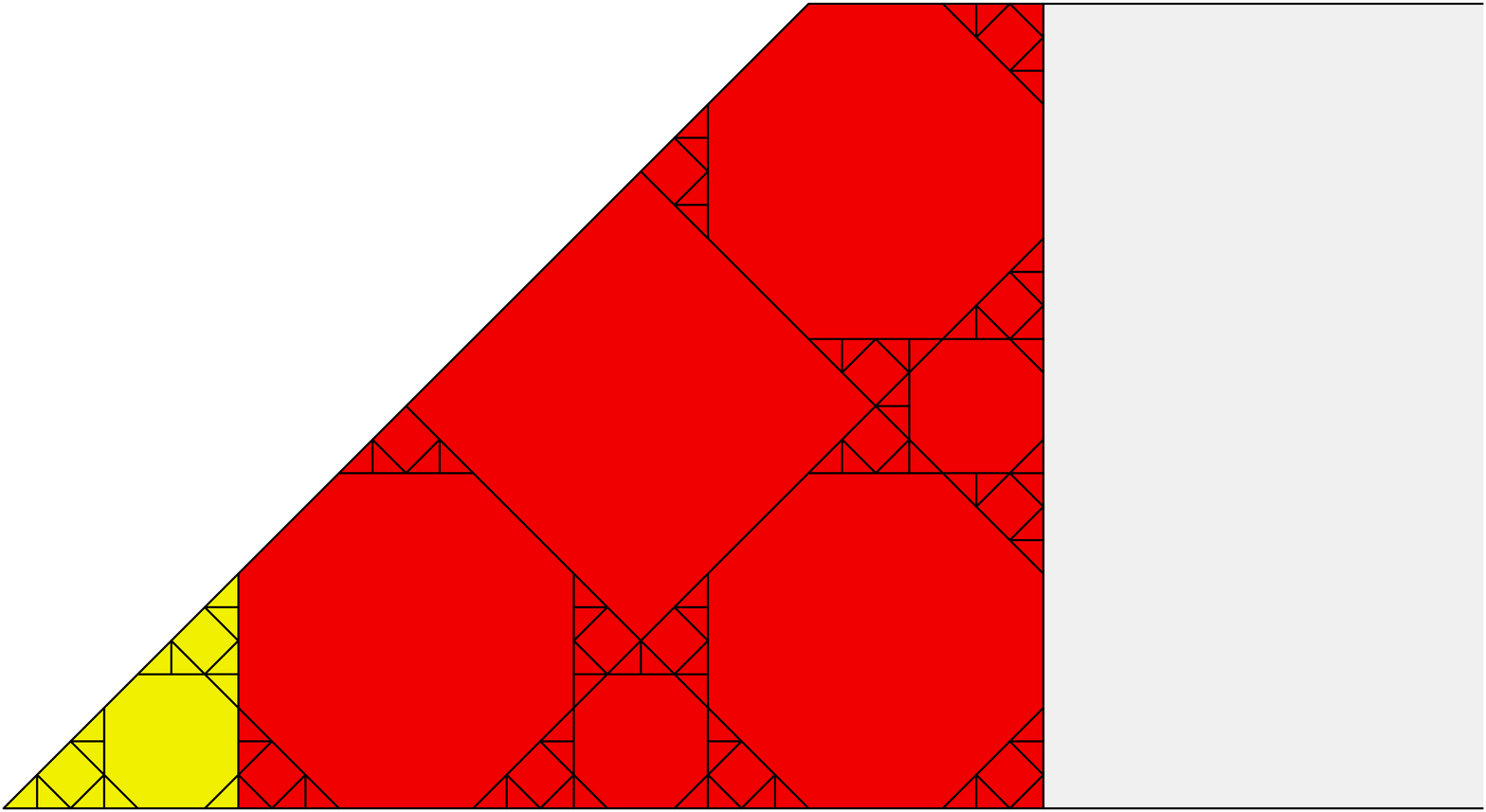}}
\newline
{\bf Figure 2.9:\/} 
$P_s$ (red) and $Q_s$ (yellow) for $s=12/31$.
\end{center}

\begin{center}
\resizebox{!}{1.7in}{\includegraphics{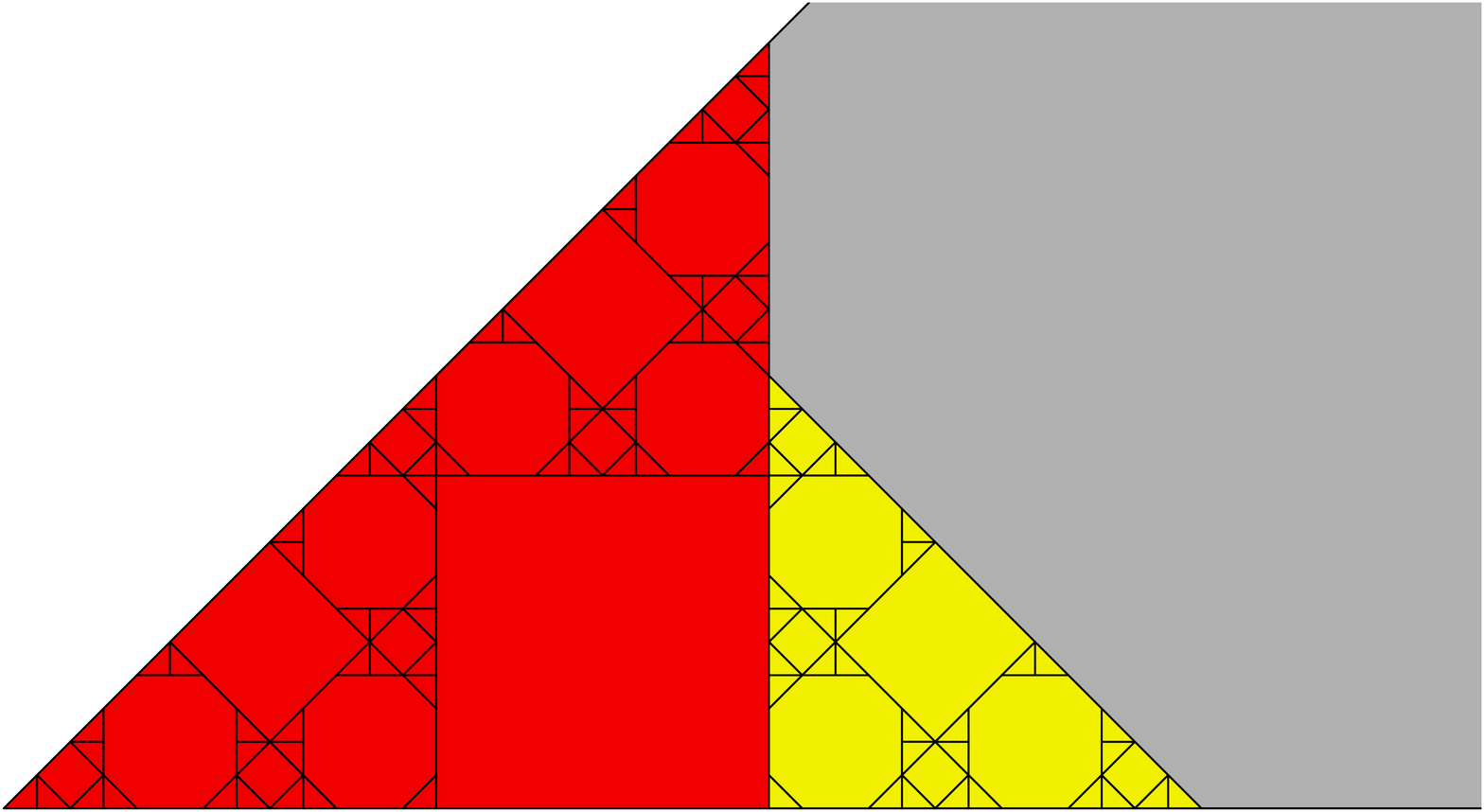}}
\newline
{\bf Figure 2.10:\/} 
$P_s$ (red) and $Q_s$ (yellow) for $s=18/23$.
\end{center}

\noindent
{\bf Remarks:\/} \newline
(i)
In case $s \in (1/2,1)$, the bilateral symmetry behaves 
nicely with respect to the inversion symmetry.
The value $r=1/(2s)$ also lies in $(1/2,1)$.  There
is a similarity carrying $\Delta_r$ to $\Delta_s$,
and this similarity carries the regions
$(A_r,B_r,P_r,Q_r)$ to the regions
$(Q_s,P_s,B_s,A_s)$.  So, for $s \in (1/2,1)$, the
symmetries associated to $(P_s,Q_s)$ are equivalent
to the ones associated to $(A_s,B_s)$. This is
not true for $s <1/2$.
\newline
(ii)
Looking at the figures, the reader will probably be able
to find other regions of bilateral symmetry.  I think
that all the bilateral symmetry one sees in the picture
is a consequence of the ones already mentioned, together
with renormalization.
\newline
(iii)
In [{\bf S0\/}] we defined
$A_s$ to be the hexagon $X_s \cap R_H(X_s)$,
but here it seems better just to take the
portion of this hexagon which lies to the
left of the central tiles of $\Delta_s$.
This notation change causes no troubles.
\newline
(iv) The sets $A_s,B_s,P_s,Q_s$, which we call
{\it symmetric sets\/}, are all clean. This follows
from the fact that each side of one of these sets
either lies in $\partial X_s$ or can be moved to
$\partial X_s$ by the bilateral symmetry.

\subsection{Squares in the Tiling}

Here we sketch the proof of a result from
[{\bf S0\/}] which illustrates 
some of the power of the results mentioned
above.

\begin{lemma}
\label{box}
When $s$ is rational and oddly even,
$\Delta_s$ contains only squares and
right-angled isosceles triangles.
When $s$ is irrational and oddly even,
$\Delta_s$ contains only squares.
\end{lemma}

\startproof
When $s=1/2$ we check that
$\Delta_s$ consists entirely of squares and right-angled
isosceles triangles.  From the Insertion Lemma, the
same result holds
for $s=1/(2n)$ for $n=2,3,4...$ In general, 
the result follows from induction on the
length of the orbit $R^n(x)$ and
Theorem \ref{renorm}.  In [{\bf S0\/}] we
showed that the triangles vanish in the irrational
limit.  So, when $s$ is irrational and
oddly even, $\Delta_s$ consists only of squares.
\endproof

$\Delta_s$ contains two kinds of squares.  We say that a
{\it box\/} is a square whose sides are parallel to the
coordinate axes.  We say that a {\it diamond\/} is a
square whose sides have slope $\pm 1$.  Then
$\Delta_s$ consists entirely of boxes and diamonds.
This follows from induction and the same kin dof proof
given in Lemma \ref{box}.
\newpage

\section{The Pattern of Filling}

\subsection{The First Half}
\label{fill1}

In this chapter we elaborate on Theorem
\ref{renorm}, explaining more precisely how
$\Delta_s$ is composed of parts of
$\phi_s(\Delta_t \cap Y_t)$.  Here $t=R(s)$
and $\phi_s$ are as in Theorem \ref{renorm}.

We first consider the case when $s<1/2$.
Let $t=R(s)$ and $u=R(t)$.
\begin{itemize}
\item If $t>1/2$ let $K=1$.
\item If $t<1/2$ and $u<1/2$, let $K={\rm floor\/}(1/(2t))$.
\item If $t<1/2$ and $u>1/2$, let $K=1+{\rm floor\/}(1/(2t))$.'
\end{itemize}
The number $K=K(s)$ plays an important role in our constructions.
We call it the {\it layering constant\/}

Let $U_s$ be the image, under of
$\phi_s$, of the trivial tile in $\Delta_t$.
This is the dark red tile in Figures 3.1 and 3.2 below.
Define
\begin{equation}
\Psi_s^0=Z_s^0 \cup U_s.
\end{equation}

Let $\tau_s$ denote the subset of
$Z_s^0$ lying beneath the line extending the
top right edge of $U_s$.  Here
$\tau_s$ is the colored region in
Figures 3.1 and 3.2.

\begin{center}
\resizebox{!}{2.5in}{\includegraphics{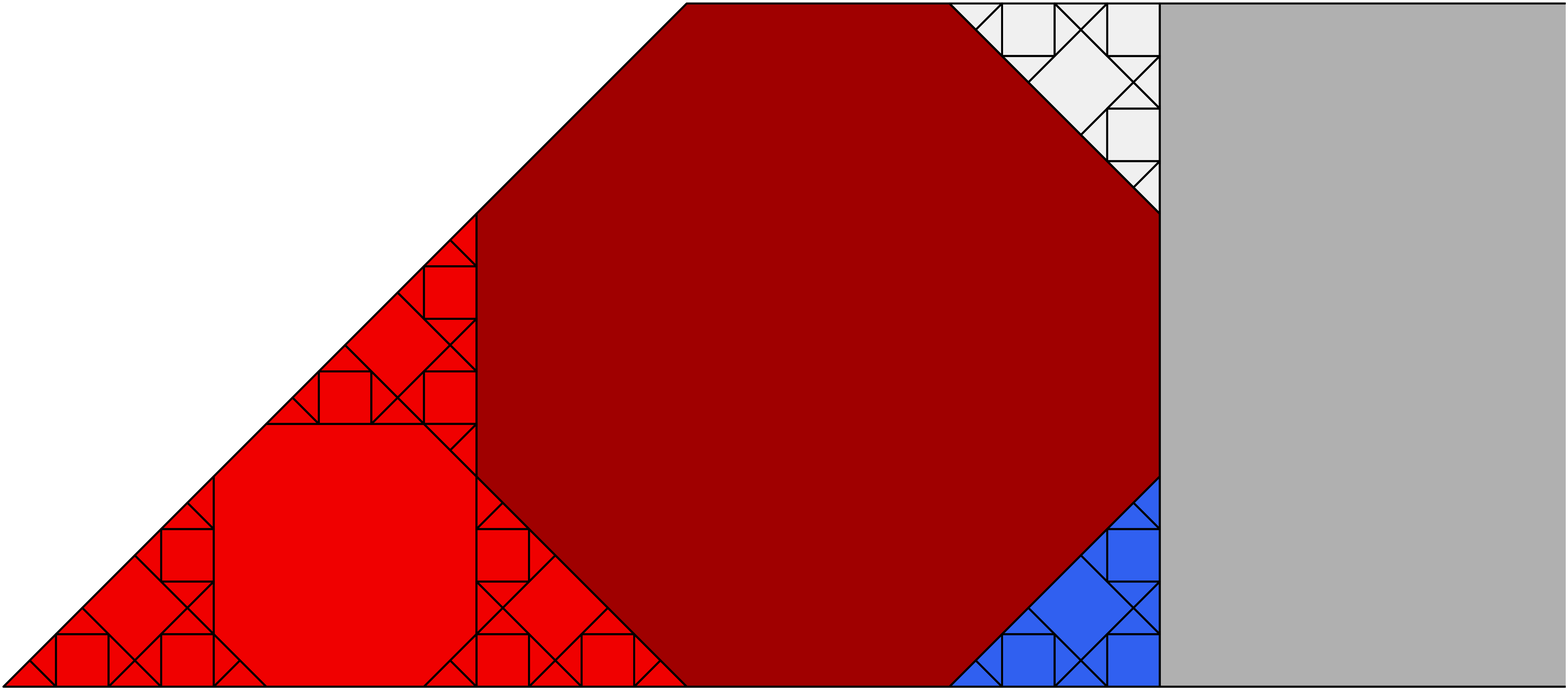}}
\newline
{\bf Figure 3.1:\/} $\Psi_s^0$ (red) and $U_s$ (dark red) and
$\tau_s$ (red, blue) for $s=13/44$.
\end{center}

\begin{center}
\resizebox{!}{2.1in}{\includegraphics{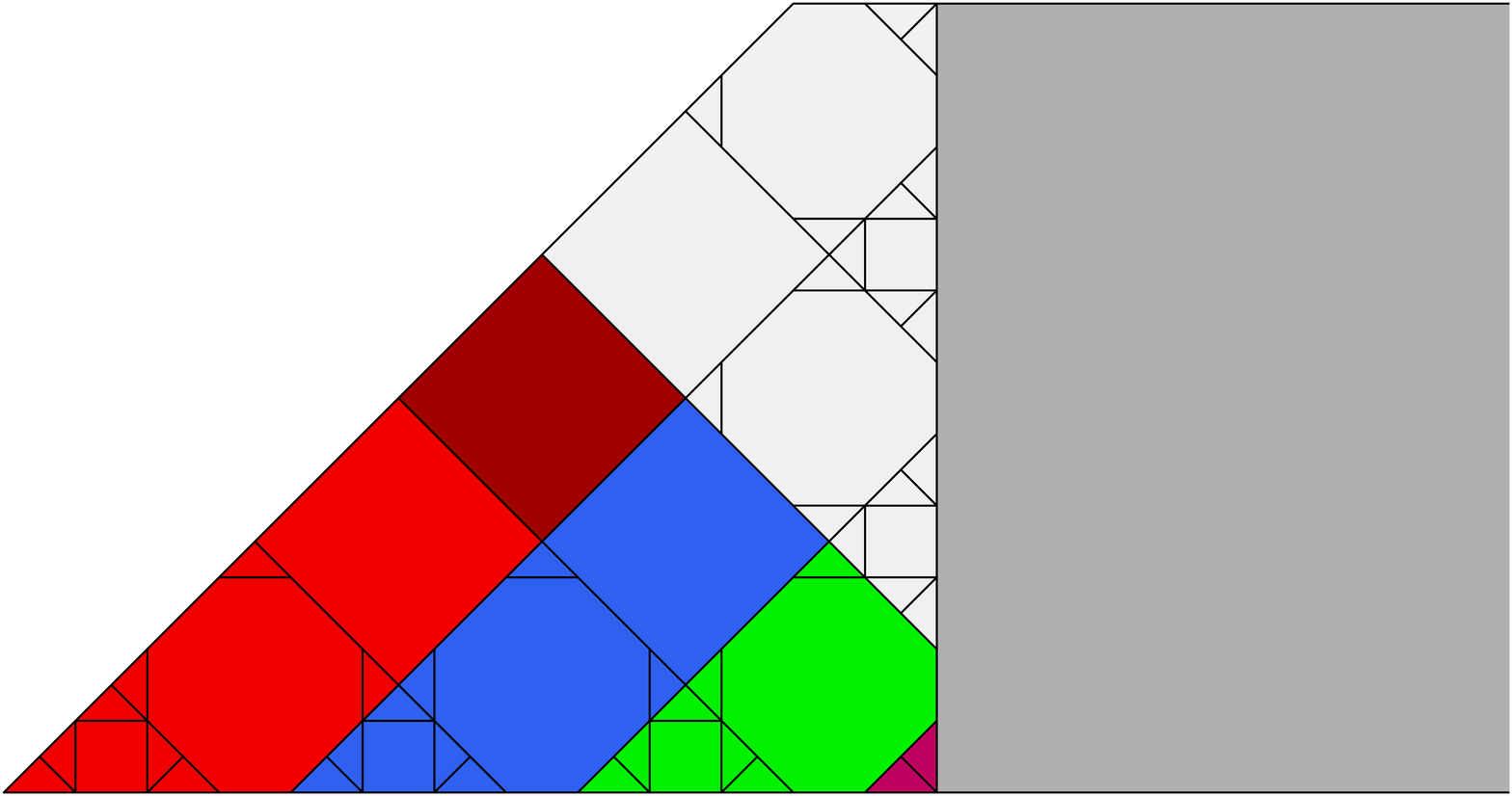}}
\newline
{\bf Figure 3.2:\/} 
$\Psi_s^j$ for $j=0,1,2,3$ (red, blue, green, magenta) for
$s=11/26$.
\end{center}

Let $T_s$ denote
the transformation which translates by
a vector pointing in the positive $x$ direction
and having length equal to the length of
the bottom side of $\Psi_s^0$.
Define
\begin{equation}
\label{tube}
\Psi_s^j=T_j(\Psi_s^0) \cap \tau_s.
\end{equation}

We have
\begin{equation}
\label{trifill}
\label{floor}
\tau_s=\bigcup_{j=0}^K \Psi_s^j.
\end{equation}
For larger $j$, the sets in Equaation \ref{trifill}
are empty.  

The whole tiling is determined by the
tiling inside $\tau_s$ and symmetry.
\begin{equation}
\label{iota2}
X_s = X_s^0 \cup {\rm central\ tiles\/} \cup
\iota(X_s^0), \hskip 30 pt
X_s^0 = \tau_s \cup R_D(\tau_s).
\end{equation}
The second equation is a consequence
of th fact that the top right boundary of
$\tau_s$ is parallel to the fundamental symmetry
line $D_s$ and lies above it.

The following result explains the structure of 
$\Delta_s$ inside $\tau_s$.

\begin{lemma}[Filling]
$T_s^{-j}(\Delta_s \cap \Psi_s^j)=\Delta_s \cap T^{-j}(\Psi_s^j)$
for all $j=1,...,K$.
\end{lemma}

\startproof
What the result means is that the map $T_s^{-j}$
respects the tiling, at least on the relevant domain.

By the Insertion Lemma, we can take
$s \in (1/4,1/2)$.  Referring to the basic
map of our system, the map
$f_s: X_s \to X_s$ is defined in terms of
a partition of $X_s$.  On each piece of this
partition, $f_s$ is a translation.  Let
$\Omega_s$ be the piece of the partition
which shares the lower left vertex of $X_s^0$.
A routine calculation shows that
the restriction of $f_s$ to $\Omega_s$
is $T_s$.  

When $s \in (1/3,1/2)$, the set
$\Omega_s$ is the quadrilateral with the
following properties. 
\begin{itemize}
\item The left edge of $\Omega_s$ is 
contained in the left side of $X_s$.
\item The bottom edge of $\Omega_s$ is
contained in the bottom edge of $X_s$.
\item The top edge of $\Omega_s$ lies
in the same line as the top edge of $Z_s^0$.
\item The right edge $e_s$ of $\Omega_s$ is vertical
and has the property that $T_s(e)$ lies in the
left edge of the leftmost central tile of $\Delta_s$.
\end{itemize}
When $s \in (1/4,1/3)$, the set
$\Omega_s$ is a triangle whose
left, right, and bottom sides are as above.
Our argument works the same in either
case.

Figures 4.3-4.5 illustrate the picture
for $s=19/60,11/30,9/20.$
The diagonal line $L_s$ bisects the yellow
square and the green diamond in each picture.

\begin{center}
\resizebox{!}{2in}{\includegraphics{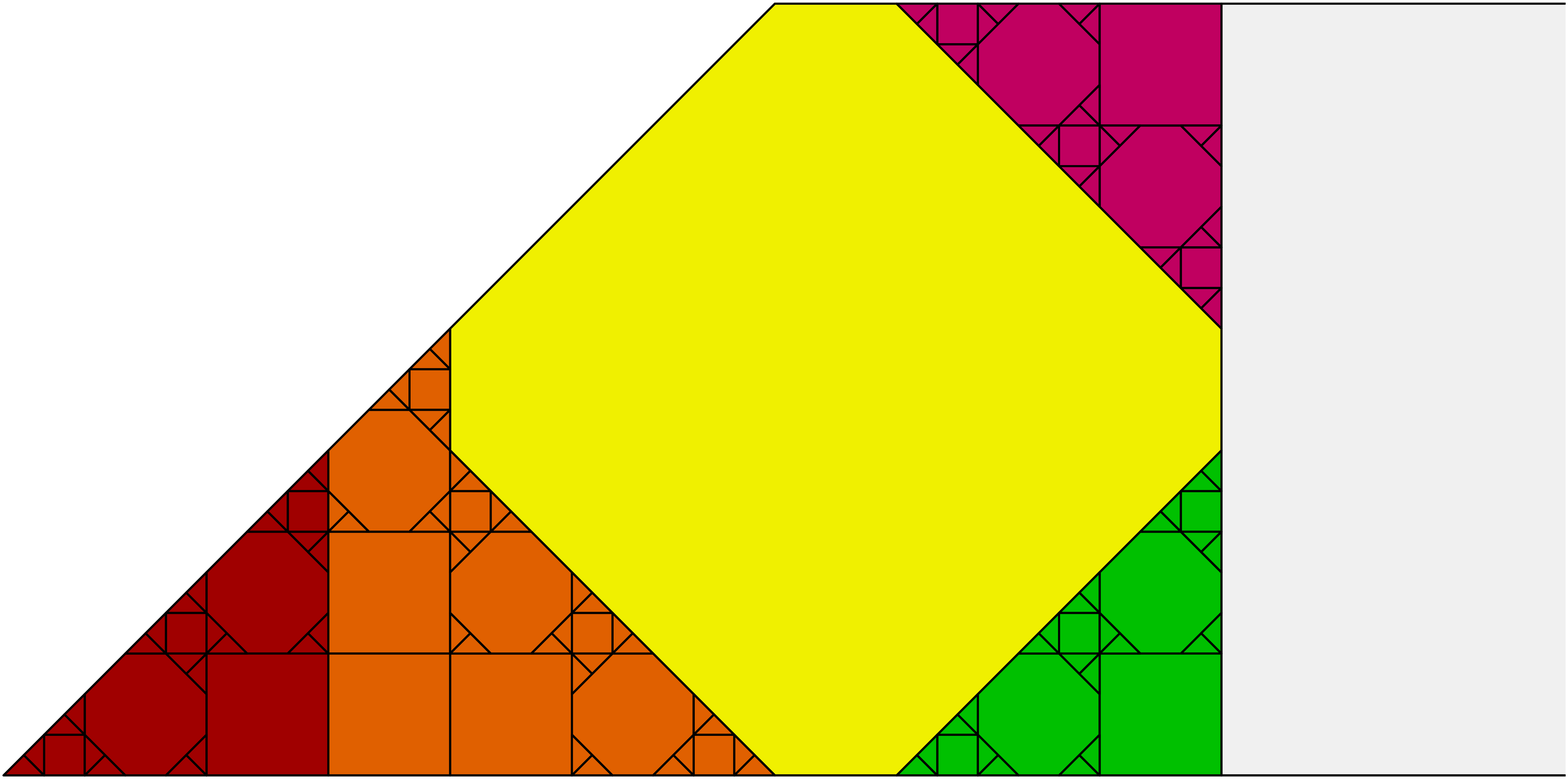}}
\newline
{\bf Figure 3.3:\/} 
$Z_s^0$ (red, orange),
$\Omega_s$ (red) and
$f_s(\Omega_s)$ green
\end{center}

\begin{center}
\resizebox{!}{2in}{\includegraphics{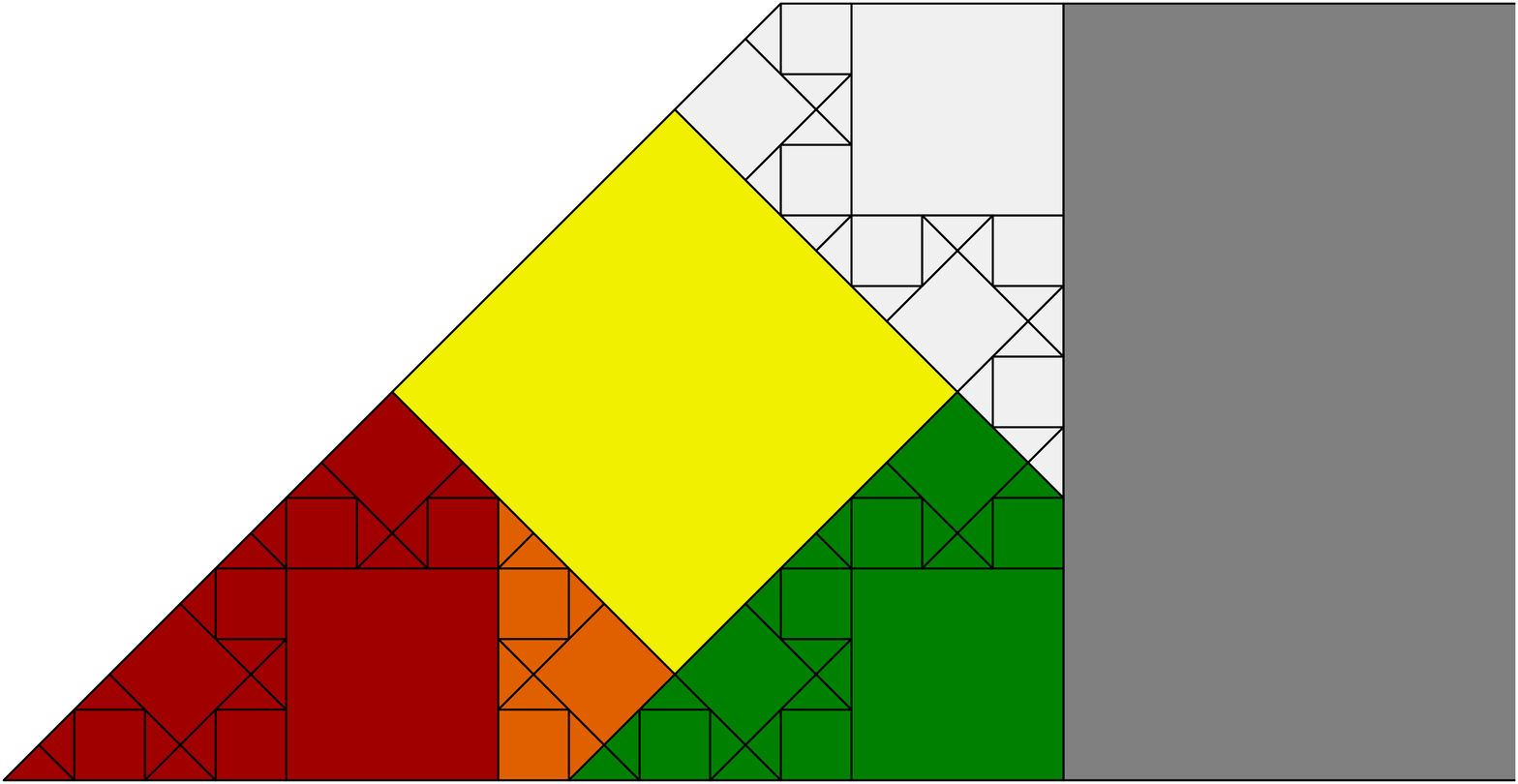}}
\newline
{\bf Figure 3.4:\/} $Z_s^0$ (red, orange) and
$\Omega_s$ (red) and $f_s(\Omega_s)$ (green).
\end{center}

\begin{center}
\resizebox{!}{2.8in}{\includegraphics{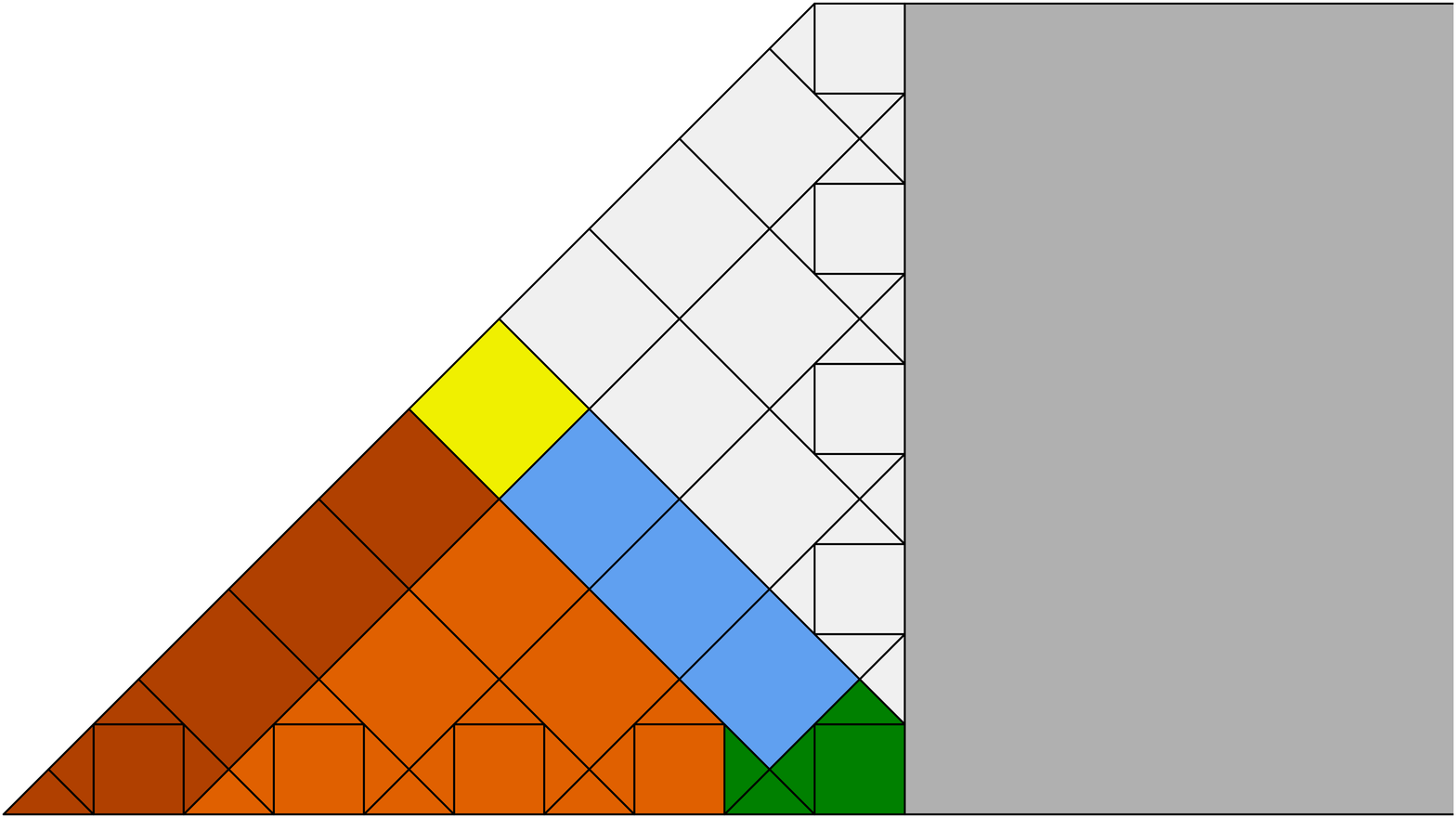}}
\newline
{\bf Figure 3.5:\/} $Z_s^0$ (red),
$\Omega_s$ (red, orange), and
$f_s(\Omega_s)$ (green, orange, blue).
\end{center}

A routine calculation shows that
\begin{equation}
T_s^{-1}(\tau_s-\Psi_s^0) \subset \Omega_s.
\end{equation}
If we start with any point $p \in \tau_s-\Psi_s^0=Z_s^0 \cup U_s$,
we see from the shape of $\Omega_s^0$ that
the iterates $f_s^{-j}(p)$ are defined and lie in
$\Omega_s^0$, for each $j=0,1,2,...$ until we reach some $k \leq K$ such that
$f_s^{-k}(p) \subset Z_s^0$.
Our result follows from this observation and from the
fact that $\Delta_s$ is $f_s$-invariant.
\endproof

Say that a {\it central tile of\/} $\Delta_s \cap \Psi_s^0$ is
the image of a central tile of $\Delta_t$ under the map
$\phi_s$.  For instance, in Figure 3.5, there are $4$
central tiles, one yellow and $3$ dark red.  The central
tiles all have the same size, and lie in
$\Psi_s^0$.

For $j<K$, the tiling
$\Delta_s \cap \Psi_s^j$ has a simple
description.  
\begin{itemize}
\item Start with the tiling $\Delta_s \cap \Psi_s^0$.
\item Chop off the top $j$ central tiles.
\item Translate by $T_s^j$.
\end{itemize}
The tiling $\Delta_s \cap \Psi_s^K$ is more complicated,
but we do not need to know it explicitly.

\subsection{The Second Half}

Now we explain the picture for $s \in (1/2,1)$.
This time, we define
\begin{equation}
K={\rm floor\/}\Big(\frac{1}{2-2s}\Big)
\end{equation}
Again, we call $K$ the layering constant.

This time, the right edge of $Z_s^0$ lies on 
the same line as the left edge of the central
tile of $\Delta_s$.  Let $\delta_s$ denote
this line.  Let $T_s$ denote the translation by
the vector which is positivel proportional
to $(1,1)$ and whose length is the same as the
length of the left side of $Z_s^0$.
Let $\tau_s$ be the
region of $X_s^0$ lying to the left of
$\delta_s$.  In Figure 3.6, the region
$\tau_s$ is colored blue/red/yellow.

\begin{center}
\resizebox{!}{2in}{\includegraphics{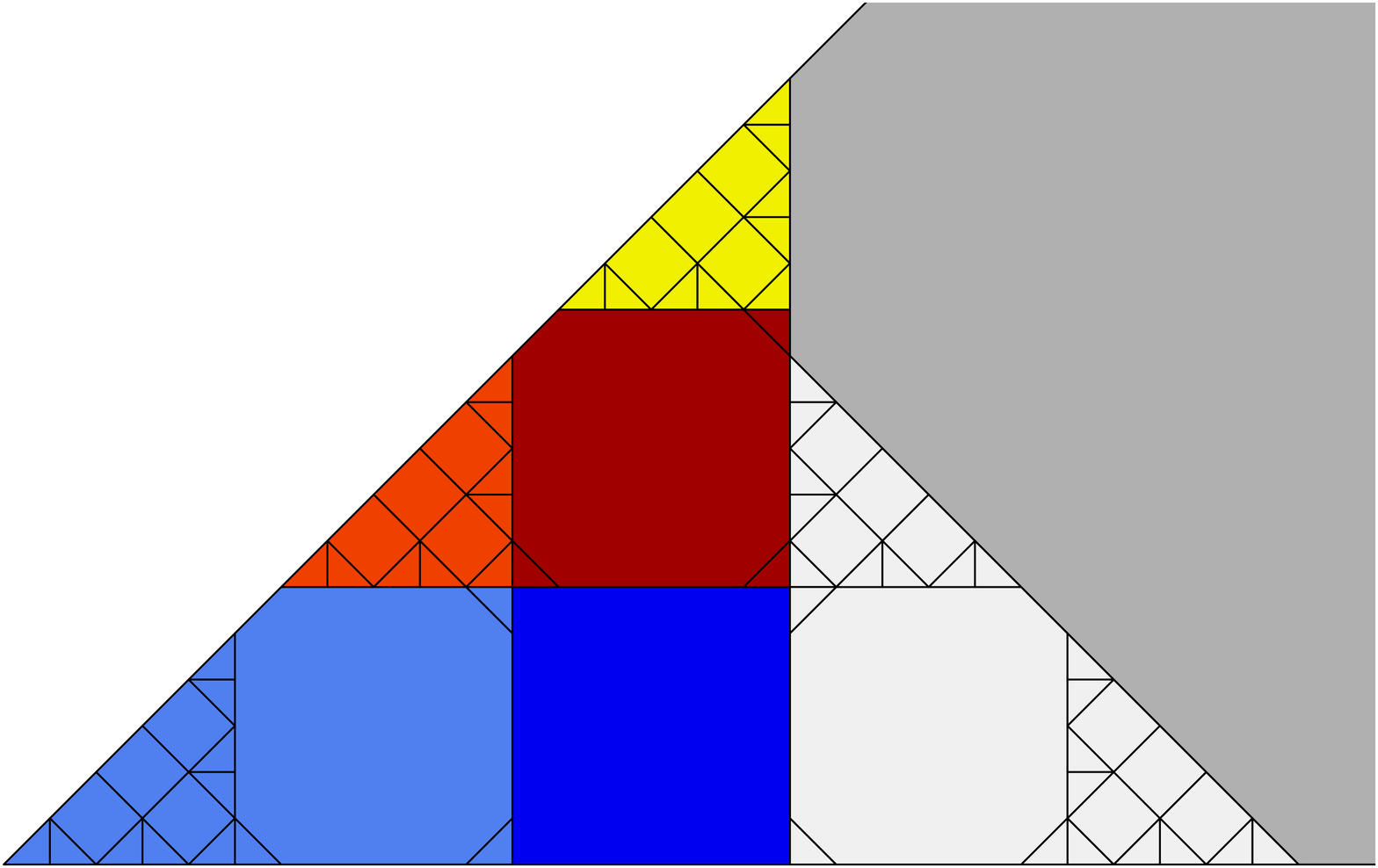}}
\newline
{\bf Figure 3.6\/}: $\Psi_s^j$ for $j=1,2,3$ (red, blue, yellow)
for $s=14/17$.
\end{center}

Define
\begin{equation}
\Psi_s^j=T_s^j(Z_s) \cap \tau_s.
\end{equation}
We have
\begin{equation}
\label{trifill2}
\label{floor2}
\tau_s=\bigcup_{j=0}^K \Psi_s^j.
\end{equation}

This time we have
\begin{equation}
\label{iota3}
X_s = X_s^0 \cup {\rm central\ tiles\/} \cup
\iota(X_s^0), \hskip 30 pt
X^0_s \subset \tau_s \cup R_V(\tau_s),
\end{equation}
where $V$ is the vertical line of symmetry, $x=-1$.

With these definitions, the Filling Lemma holds
{\it verbatim\/}, and the proof is essentially
the same.  This time $\Omega_s$ is a right isosceles
triangle, and the left edge of $f_s(\Omega_s$ lies
in $\delta_s$, and $f_s$ translates diagonally
along the vector that generates the left side of $Z_s$.
In Figure 3.6, $\Omega_s$ is the union of the light
red and light blue tiles.

\subsection{Pyramids}

Now we explain more of the structure.  We will
concentrate on the case $s<1/2$.  

Let $K=K(s)$ be the layering constant.
Say that a
{\it pyramid\/} of size $k$ is a configuration of
diamonds having the structure indicated in Figure
3.7 for $k=1,2,3$. 
We refer to the longest (diagonal) row of
diamonds as the {\it base\/} if the pyramid.

\begin{center}
\resizebox{!}{1.5in}{\includegraphics{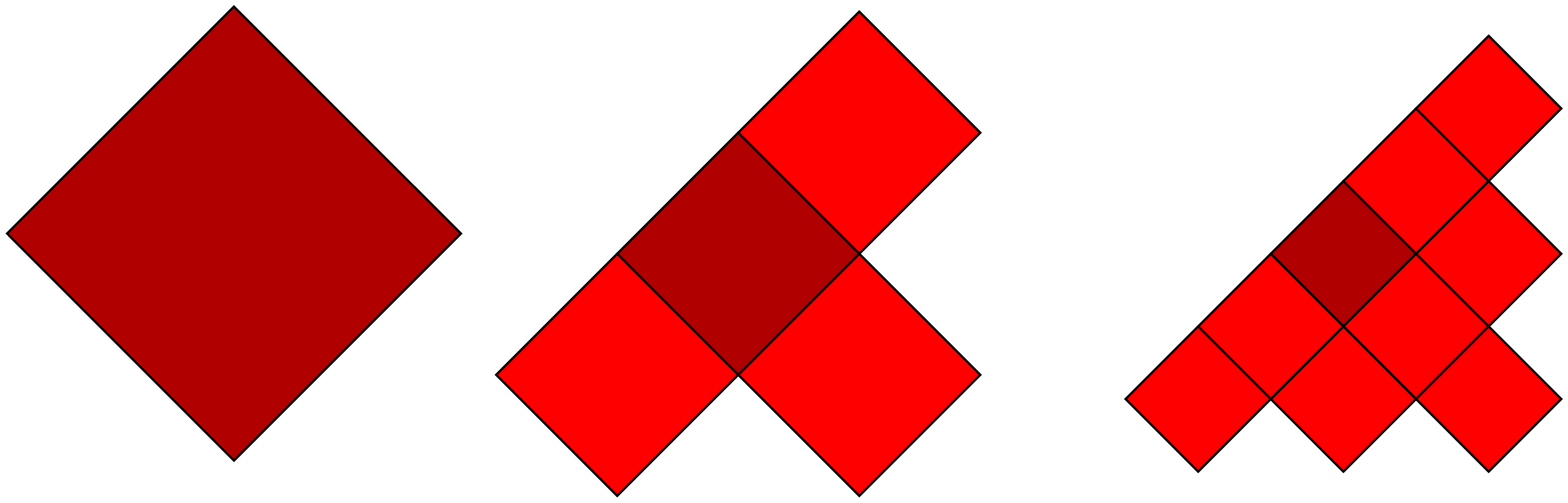}}
\newline
{\bf Figure 3.7:\/} 
Pyramids for $k=1,2,3$.
\end{center}

Figure 3.8 shows a red pyramid of size $2$ contained
in $\Delta_s$ for $s=30/73$.  In the next section, we
will explain the meaning of the blue squares shown
in the figure.

\begin{center}
\resizebox{!}{3in}{\includegraphics{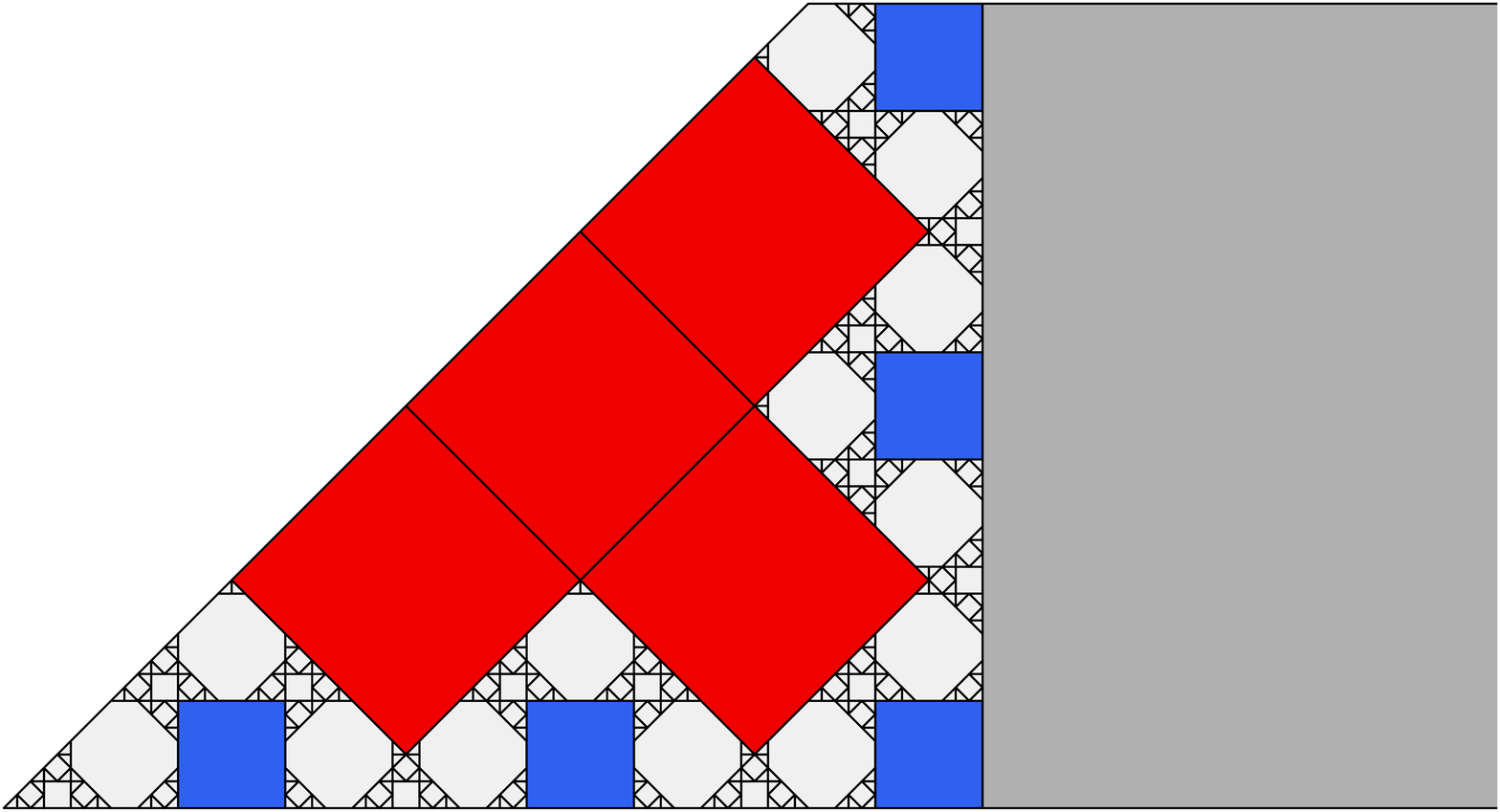}}
\newline
{\bf Figure 3.8:\/} 
The pyramid. $s=30/73$, $R(s)=13/60$ and
$R^2(s)=4/13$.
\end{center}

Say that an {\it extended pyramid\/} of size $K$ is the
union of polygons obtained by taking the outer squares
in each row and chopping off the corners so as to leave
semi-regular octagons.  Figure 3.9 shows an example
of an extended pyramid of size $K$.

\begin{center}
\resizebox{!}{2.1in}{\includegraphics{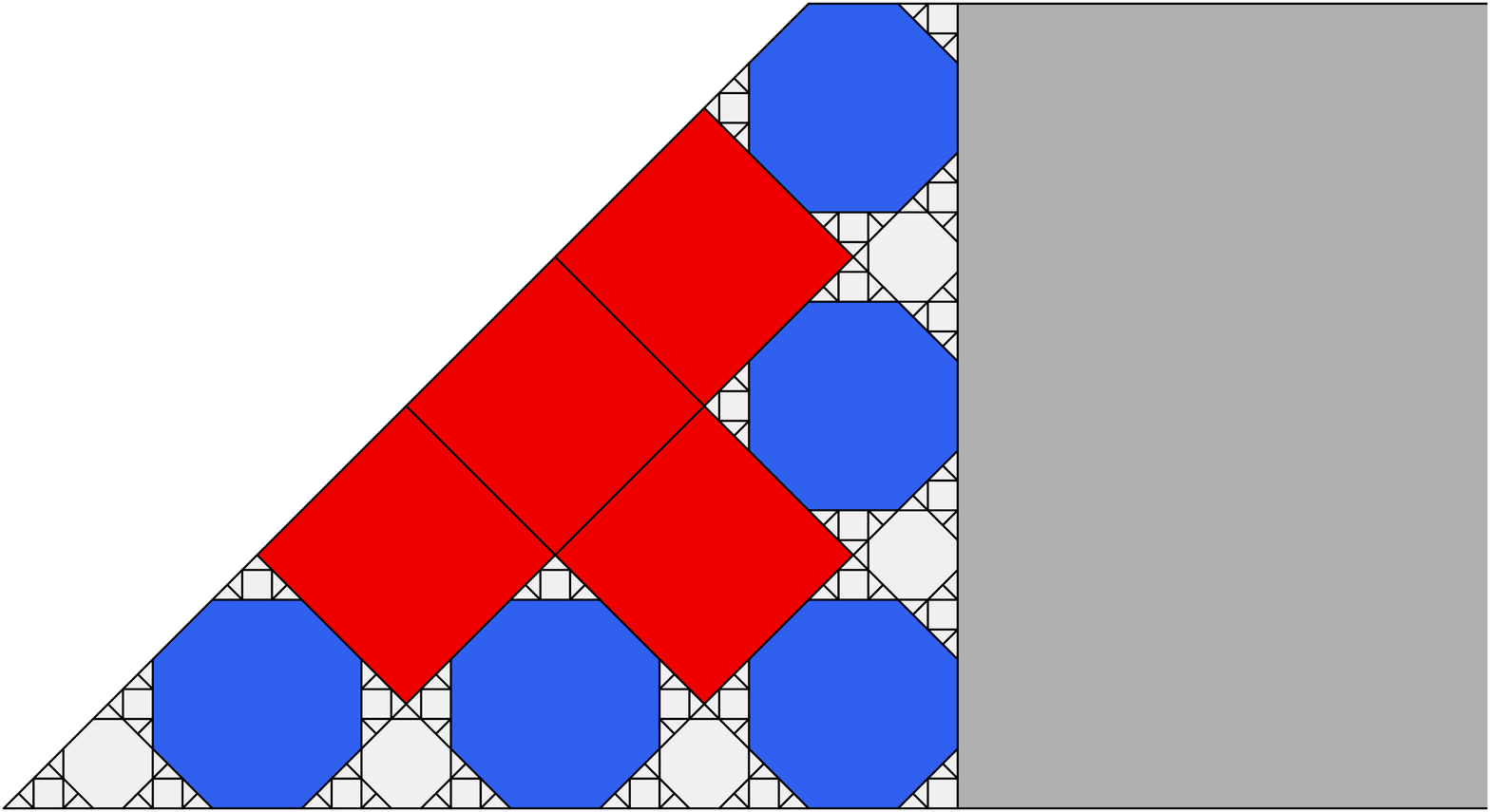}}
\newline
{\bf Figure 3.9:\/} 
{\small Extended Pyramid. $s=27/64$,
$R(s)=5/27$ and $R^2(s)=7/10$.\/}
\end{center}

Say that $s<1/2$ has {\it type 0\/} if $R(s)<1/2$
and $R^2(s)<1/2$.  Otherwise, say that $s$ has
{\it type 1\/}.

\begin{lemma}
\label{pyramid}
Suppose $s$ has type $0$ (respectively type $1$).
Then $\Delta_s$ contains an ordinary (respectively 
extended) pyramid of size $K$.
\end{lemma}

\startproof
Suppose first that $R^2(s)<1/2$.
The bottom $K$ squares in the base of
the desired pyramid are guaranteed by
Theorem \ref{renorm} and the Insertion
Lemma (applied to $\Delta_t$).  The
bottom half of the pyramid is then guaranteed
by the Filling Lemma.  The top half is then
guaranteed by the bilateral symmetry corresponding
to the diagonal line $D_s$.

The proof is essentially the same when
$R^2(s)>1/2$.  What happens here is that
the tiles in $\Delta_t$ just to the left
of the leftmost central tile 
is an octagon which has the same height as the
central tiles.
have the same width as the central tiles.
Here $t=R(s)$.
If we simply keep track of this additional
octagon, we get the same structure as in
the other case, except that these octagons
replace the outer layer of diamonds.
\endproof

\noindent
{\bf Remark:\/}
One has a result similar to Lemma \ref{pyramid} when
$s>1/2$.  When $s>1/2$ and $R^2(s)<1/2$ we get an
ordinary pyramid of size $K-1$.  When 
$s>1/2$ and $R^2(s)>1/2$ we get an extended pyramid
of size $K$.
\subsection{The Octagrid}

\noindent
{\bf Remark:\/}
The reader only interested in Statement 1
of the Main Theorem can skip the rest of this chapter.
\newline

In this section, we consider the case
when $s<1/2$ and $R(s)<1/2$.  One can
do something similar for other parameter
ranges, but we do not need the construction
otherwise.

First we consider the case
$R^2(s)>1/2$.
Thanks to the equalities between
the widths, the distance between the center
of an octagon and the center of an adjacent
square in the extended pyramid is the same
as the distance between two adjacent squares.
Put another way, the union of the centers of
the tiles in the extended pyramid is contained
in a square grid.

When $R^2(s)<1/2$, there is a
similar phenomenon.  This is where the
blue squares in Figure 3.8 (and Figure 3.10) come in.
These squares have the 
following description.  The bottom left
blue square is the image, under $\phi_s$,
of the central tile of $Z_t$.  The remaining
blue tiles exist as a consequence of the
Filling Lemma and bilateral symmetry.
We call the union of the pyramid and
these blue squares the {\it extended pyramid\/}.
  The union of the centers
of the tiles of the extended pyramid, in this
case also, is part of a square grid.
This follows from a routine calculation,
which we omit.

In either case, we
form the {\it octagrid\/} by taking the
union of horizontal, vertical, and diagonal
(meaning slope=$\pm 1$) lines through the
centers of the tiles of the extended pyramid.
The octagrid chops up $\Delta_s$ into what
turn out to be usefully small pieces.

\begin{center}
\resizebox{!}{2.3in}{\includegraphics{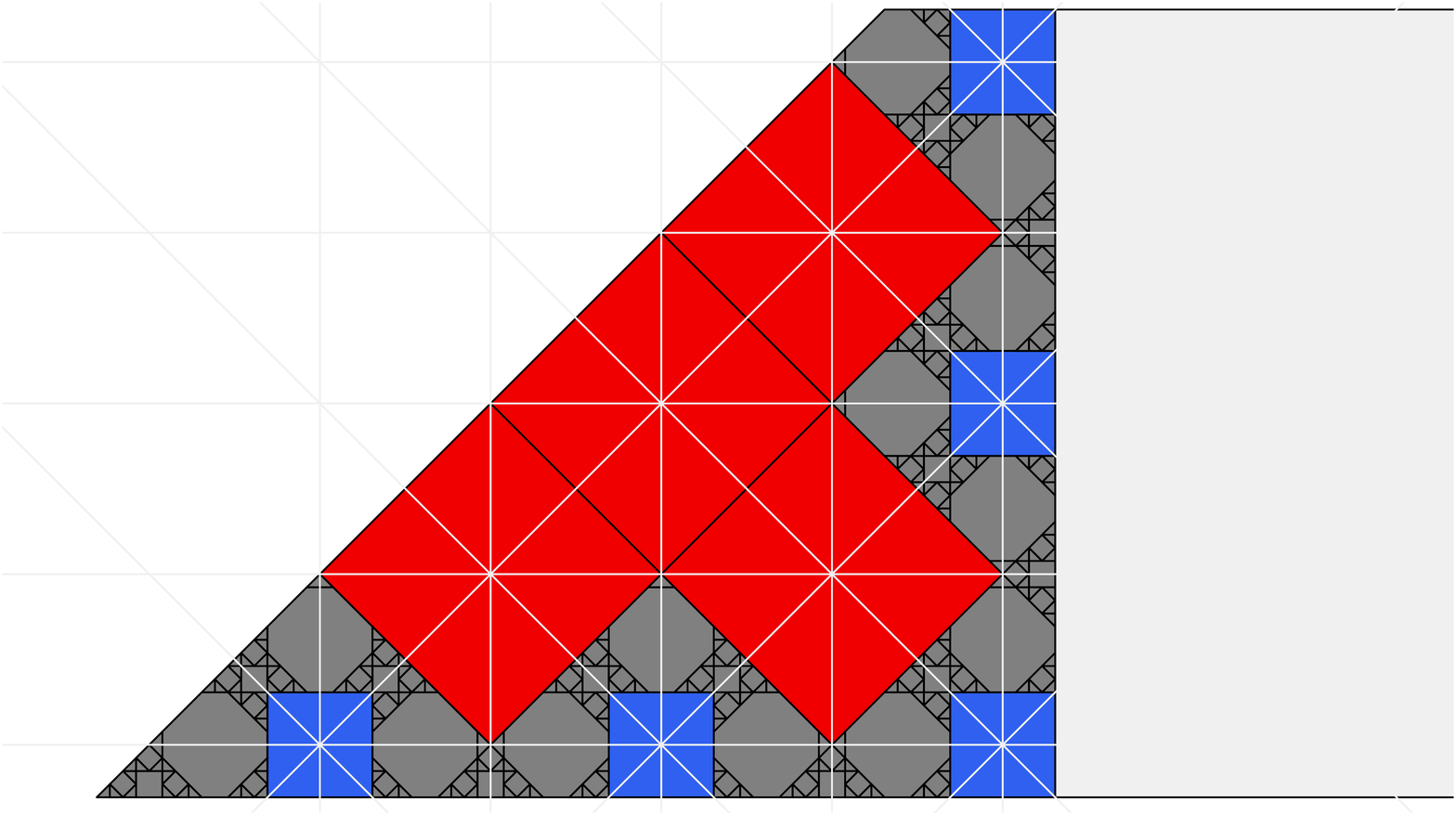}}
\newline
{\bf Figure 3.10\/}: The octagrid for $s=30/73$.
\end{center}

There are finitely many lines in the octagrid.
These lines partition $X_s^0$ into finitely many
bounded convex regions.  We call these regions
{\it octagrid components\/}.  Each
octagrid component is some open convex polygon.

\begin{lemma}
\label{octagrid}
Let $G$ be any octagrid component.  There is
an isometry carrying $\Delta_s \cap G$ to
a subset of $\Delta_s \cap Z_s^0$.
\end{lemma}

\startproof
The fundamental line $D_s$ of symmetry is contained
in the octagrid. The reflection $R_d$ carries each
octagrid component above $D$ into some octagrid
component below $D$.  Hence, by symmetry, it suffices
to prove our result when $G$ lies beneath $D_s$.

$G$ lies in the region $\tau_s$ from the Filling Lemma.
Applying the map $T_s^{-1}$ from Equation \ref{tube}
as many times as we can, we can assume that (the
interior of) $G$ intersects
$Z_s^0$.  If $G \subset Z_s^0$, we are finished.
If $G$ is contained in one of the squares of the pyramid,
we are finished.

The only remaining possibility is that $G$ has a vertex
on the centerline of $Z_s^0$ (the line $\phi_s(H)$)
but extends over the right edge of $Z_s^0$.  The
problem is that the diagonal edge of $G$ is perpendicular
rather than parallel to the centerline of $Z_s^0$.
To fix this problem, we consider $G'=T_s^j \circ R_V(G)$
for the largest value of $j$ we can take.
The interior of $G'$ intesects $Z_s^0$ and
the diagonal edge of $G'$ is parallel to
$Z_s^0$.  Hence, $G' \subset Z_s^0$.
Since we have only applied symmetries of the
tiling, we see that $\Delta_s \cap G$ and
$\Delta_s \cap G'$ are isometric.
\endproof

\begin{center}
\resizebox{!}{2.8in}{\includegraphics{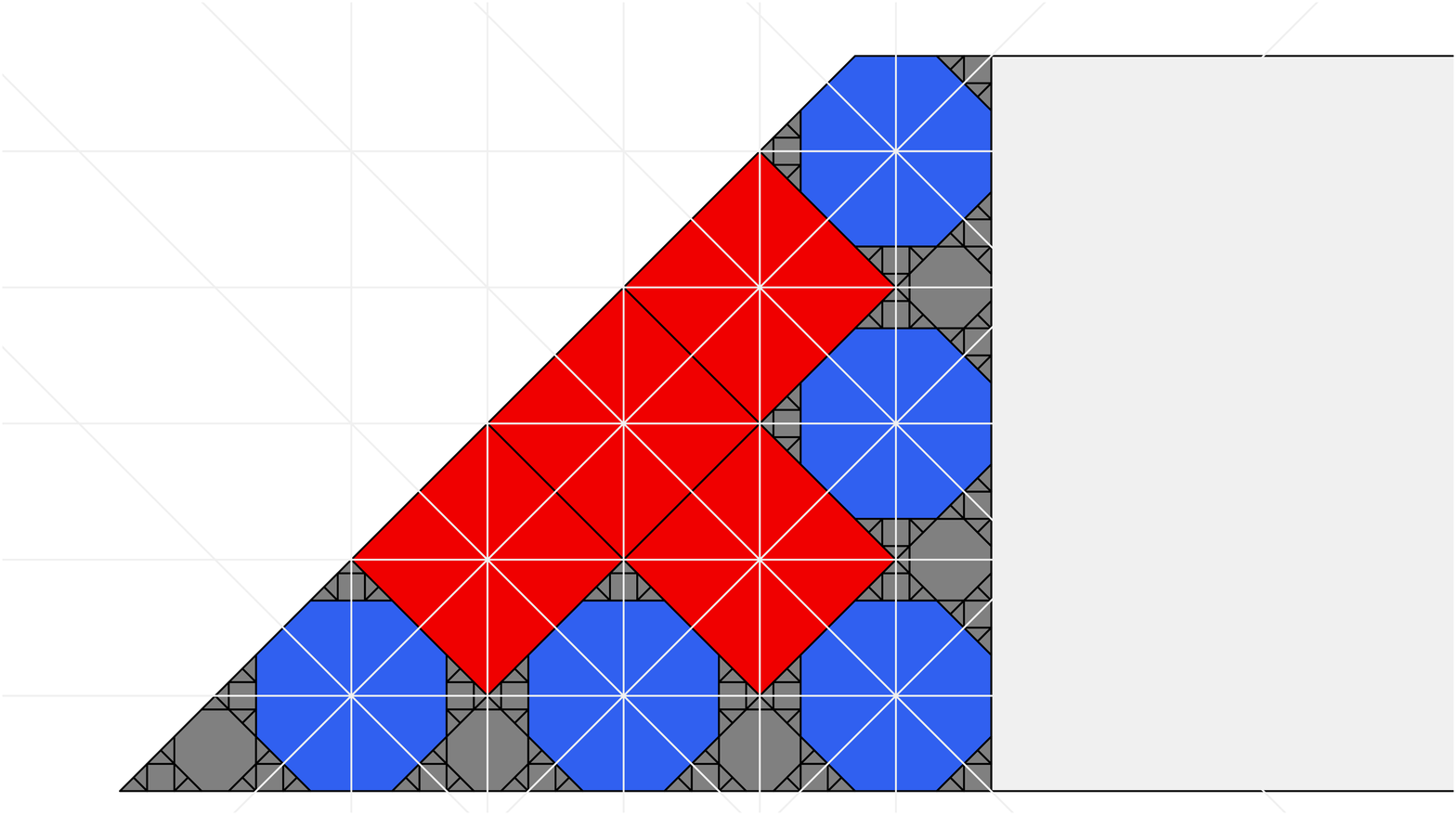}}
\newline
{\bf Figure 3.11\/}: The octagrid for $s=27/64$.
\end{center}

\newpage

\section{Geometric Lemmas}

In this chapter we prove $3$ geometric lemmas,
which we call the Shield Lemma, the Pinching Lemma,
and the Covering Lemma.  Only Statement 3 of the
Main Theorem requires the Covering Lemma.

\subsection{The Shield Lemma}

Let $A_s$ be the symmetric piece from \S \ref{symm}.

\begin{lemma}[Shield]
Let $s$ be irrational and oddly even.
Every point of $\partial A_s$, except
the $2$ vertices having obtuse angles, is contained
in the edge of square.  Those points which
belong to the boundaries of more than one
tile are the vertices of pairs of adjacent squares.
\end{lemma}

Figure 4.1 illustrates the Shield Lemma.
When $s$ is rational, there are $2$ small
triangles touching the obtuse vertices
of $A_s$.  These triangles vanish in the
irrational limit.

\begin{center}
\resizebox{!}{3in}{\includegraphics{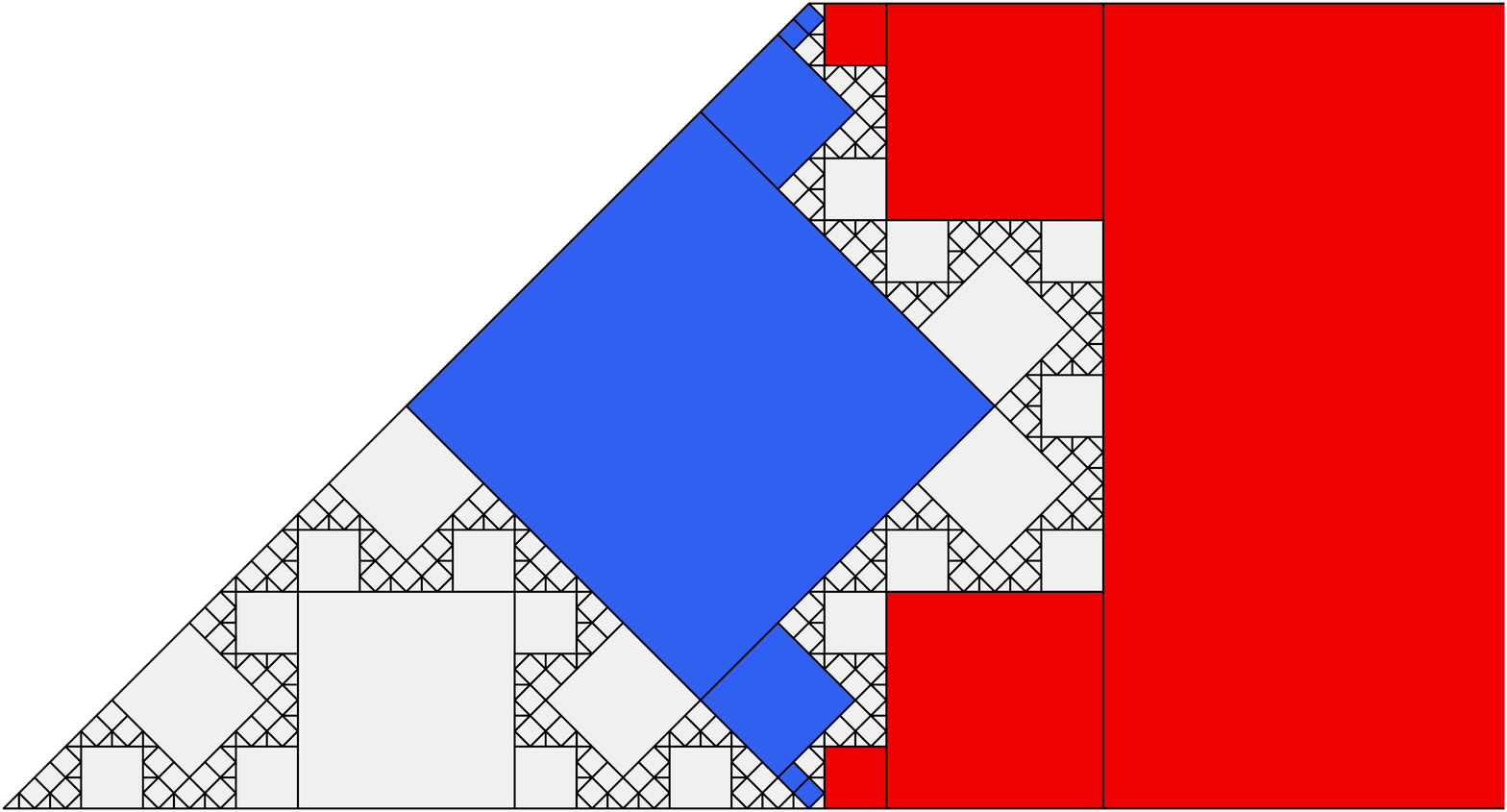}}
\newline
{\bf Figure 4.1:\/}
$A_s$ for $s=26/71$.
\end{center}

We define the {\it shield\/} $\Sigma_s$ 
to be the union of the top left edge of
$A_s$ and the left half of the top edge.
$\Sigma_s$ is the union of two line segments.
The top left vertex $\nu_s$ of $X_s$ (and $A_s$) is the
place where the two line segments join.
By symmetry, it suffices to prove
the Shield Lemma for the points of
$\partial A_s$ contained in the shield.
We analyze the picture in the rational
case and then take limits.   In this section
we work out how $\Delta_s$ sits in 
$A_s$ when $s$ is rational and oddly even.

We say that a square of $\Delta_s$
{\it abuts\/} $\Sigma_s$ if
$\overline T \cap \Sigma_s$ is nonempty.
In this case, the intersection is a line
segment which we call the {\it contact\/}
between the square and the shield.
The {\it radius\/} of $T$ 
is the distance from the center of $T$ to a corner of $T$.
We call a radius $\rho$ {\it realized\/}, if a
square of $\Delta_s$ having radius $\rho$ abuts
the shield. 

\begin{lemma}
\label{abut}
Let $s$ be an oddly even rational.  The following is true.
\begin{enumerate}
\item Exactly one
triangle of $\Delta_s$ abuts $\Sigma_s$
and the segment of contact contains $\nu_s$
as an endpoint.
\item The squares which
abut $\Sigma_s$ occur in monotone decreasing size,
largest to smallest, as one moves from
an endpoint of $\Sigma_s$ to $\nu_s$.
\item The number of squares of each size is determined
by the even expansion of $s$.
\item Let $\rho$ be a realized radius.
Some square of radius $\rho$, which abuts
$\Sigma_s$, has a vertex within $\rho$ of
$\nu_s$.
\end{enumerate}
\end{lemma}

Let us assume Lemma \ref{abut} for now, and
finish the proof of the Shield Lemma.
Let $\{r_n\}$ be a sequence of oddly even rationals
which converges to $s$.  Given the convergence
of tilings described above, we see that the union of square
tiles abutting $\Sigma_s$ is the Hausdorff limit
of the union of square tiles abutting $\Sigma_{r_n}$,
as $n \to \infty$.  The size of the single triangle
in the picture for the rational parameters tends
to $0$.  Hence, every point of $\Sigma_s-\nu_s$
is contained in a segment of contact for some square.
The main point to worry about is that somehow there
is a point $p\in \Sigma_s-\nu_s$, with the
following property:  As $n$ tends to $\infty$,
the square whose segment of contact contains
$p$ tends to $0$ in size.  This unfortunate
situation cannot occur because it would violate
Item 4 of Lemma \ref{abut}.
\newline
\newline
{\bf Proof of Lemma \ref{abut}:\/}
The proof goes by induction on the length of the orbit
$\{R^n(s)\}$.  When $s=1/2$ the result holds by
inspection.  The case $s=1/2n$ follows from the
Insertion Lemma.
For $s \not = 1/2n$, let $t=R(s)$.
By induction, all the properties of the lemma hold
for $\Delta_t$.

\begin{center}
\resizebox{!}{3.2in}{\includegraphics{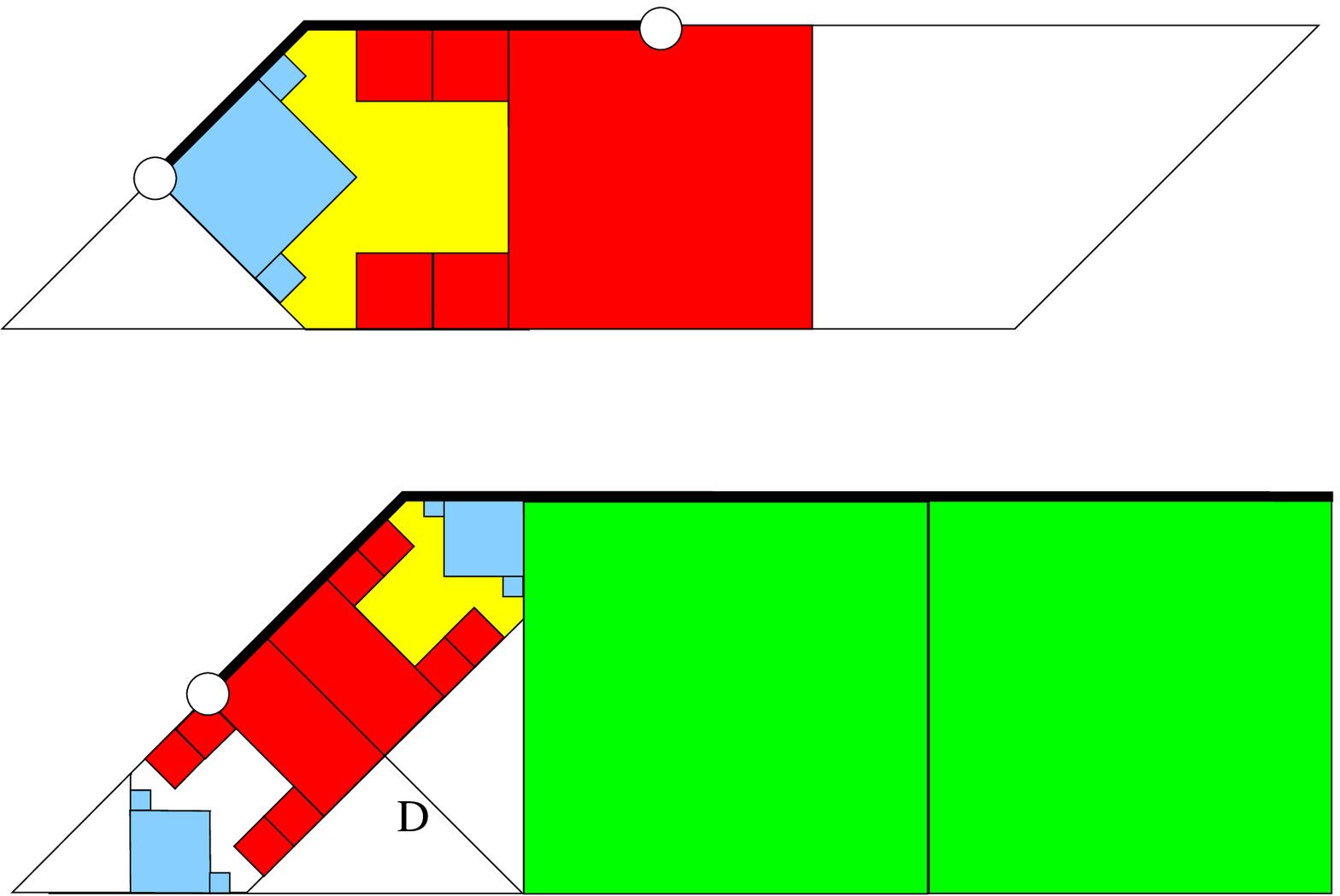}}
\newline
{\bf Figure 4.2:\/}
Inherited structure
\end{center}

Let $D=D_s$ be the diagonal line of bilateral symmetry
and let $R_D$ denote reflection in $D$.
The yellow region at the top left of Figure 4.2 is 
$A_t$.
The big blue square in the top of Figure 4.2,
which is the central square of
$\Psi_t^0$, is mapped by 
$R_D \circ \phi_s$ to a square which
abuts the (green) leftmost central square of $\Delta_s$.

The pattern of tiles abutting
the shield $\Sigma_s$ is the same as the pattern of
squares abutting the shield $\Sigma_t$, except that
some green ones have been appended.  The number of
green ones is determined by the number $n_0$ in the
even expansion of $s$. All the points in our lemma 
follow from this structure.
\endproof

\begin{corollary}
\label{good2}
Let $S_s$ denote the left half
of $\widehat \Lambda_s$.
There exists a convex set $D_2$ such that
$S_s \cap {\rm interior\/}(A_s) \subset D_2$
and $D_2$ intersects $\partial A_s$ only at
the two obtuse vertices.
\end{corollary}

\startproof
We simply slice off from $A_s$ suitably chosen
neighborhoods of the square tiles which abut the
edges of $A_s$.  With a little care (i.e., by making
these neighborhoods shrink very rapidly as we
approach the vertices) we can make the resulting
set convex.
\endproof
\subsection{The Pinching Lemma}

Let $S_s$ denote the left half of
the limit set $\widehat \Lambda_s$.

\begin{lemma}[Pinching]
At most one point of $S_s$ lies on each
fundamental line of bilateral symmetry.
\end{lemma}

We argue by contradiction.
Let $g(s)$ denote the diameter
of $X_s^0$.
Say that a {\it counterexample\/} is a
quadruple $\Omega=(L,p_1,p_2,s)$, where $L=L_s$ is a fundamental
line of symmetry and
$p_1 \not = p_2 \in L \cap S_s$.
We call $s$ the {\it parameter\/} of the counterexample.
We define
\begin{equation}
\label{harmony}
\lambda(\Omega)=\frac{\|p_1-p_2\|}{g(s)}.
\end{equation}
If $\Omega'$ is another counterexample obtained
from $\Omega$ using either the Insertion Lemma
or (when $s>1/2$) the Inversion Lemma, we
have $\lambda(\Omega)=\lambda(\Omega')$.  Moreover
\begin{equation}
\lambda=\sup_{\Omega} \lambda(\Omega) \leq 1
\end{equation}

\begin{lemma}
\label{lazy2}
For any $\epsilon>0$ there exists a
counterexample $\Omega$, whose parameter
is less than $1/2$, such that $\lambda(\Omega)>\lambda-\epsilon$.
\end{lemma}

\startproof
Certainly, there exists a counterexample $\Omega=(L,p_1,p_2,s)$
such that $\lambda(\Omega)>\lambda-\epsilon$.
If $s<1/2$ we are finished.  If $s>\sqrt 2/2$ we can
apply the Inversion Lemma to reduce to the case
$s \in (1/2,\sqrt 2/2]$.  In this case, we have $K=1$,
where $K$ is the layering constant.

Let $t=R(s)=1-s$.
When $K=1$ we have the following facts.
\begin{itemize}
\item $R_D \circ \phi_s(B_t)=A_s$ (red) is a neighborhood of $H$ in $A_s$.
\item $R_D \circ \phi_s(A_t)$ (green) is a neighborhood of $V$ in $B_s$.
\item $\phi_s(P_t)$ (light green) is a neighborhood of $D_s$ in $P_s$.
\item $R_V \circ \phi_s(Q_t)$ (yellow) is a neighborhood of $E_s$ in $Q_s$.
\end{itemize} 
The colors refer to the various regions
in $\Delta_s$ shown in Figure 4.3.
In all cases, we can use the relevant map to
pull back the counterexample isomerically to $\Delta_t$,
and we get a counterexample involving $t$ with the
same diameter.
 But then the $\lambda$-value of this
counterexample does not decrease.  Hence, the new
counterexample $\Omega'$ satisfies $\lambda(\Omega')>\lambda-\epsilon$
and has parameter $t<1/2$.
\endproof

\begin{center}
\resizebox{!}{2in}{\includegraphics{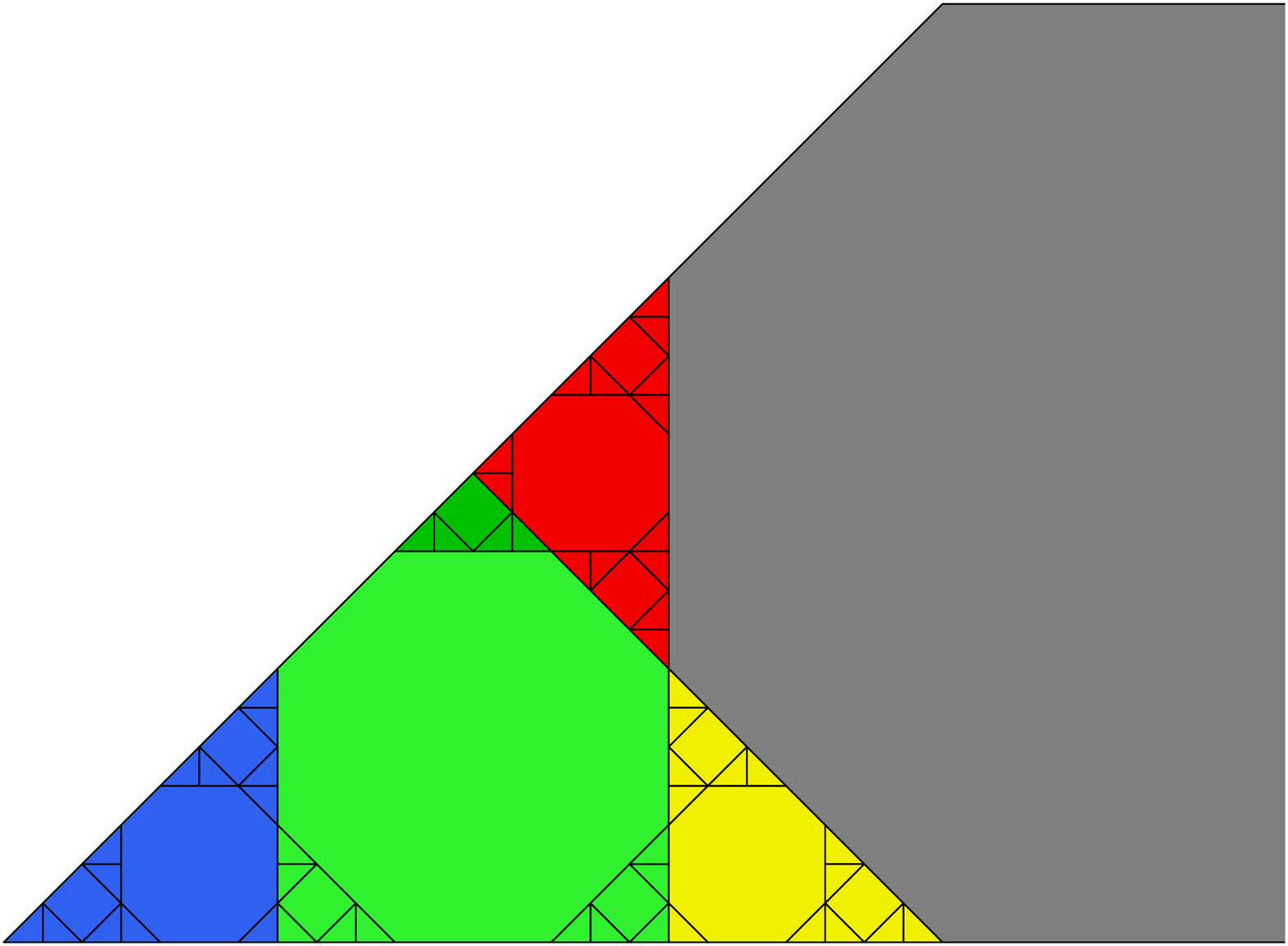}}
\newline
{\bf Figure 4.3:\/} 
Various regions in $\Delta_s$ for $s=12/17$.
\end{center}

\begin{center}
\resizebox{!}{1.5in}{\includegraphics{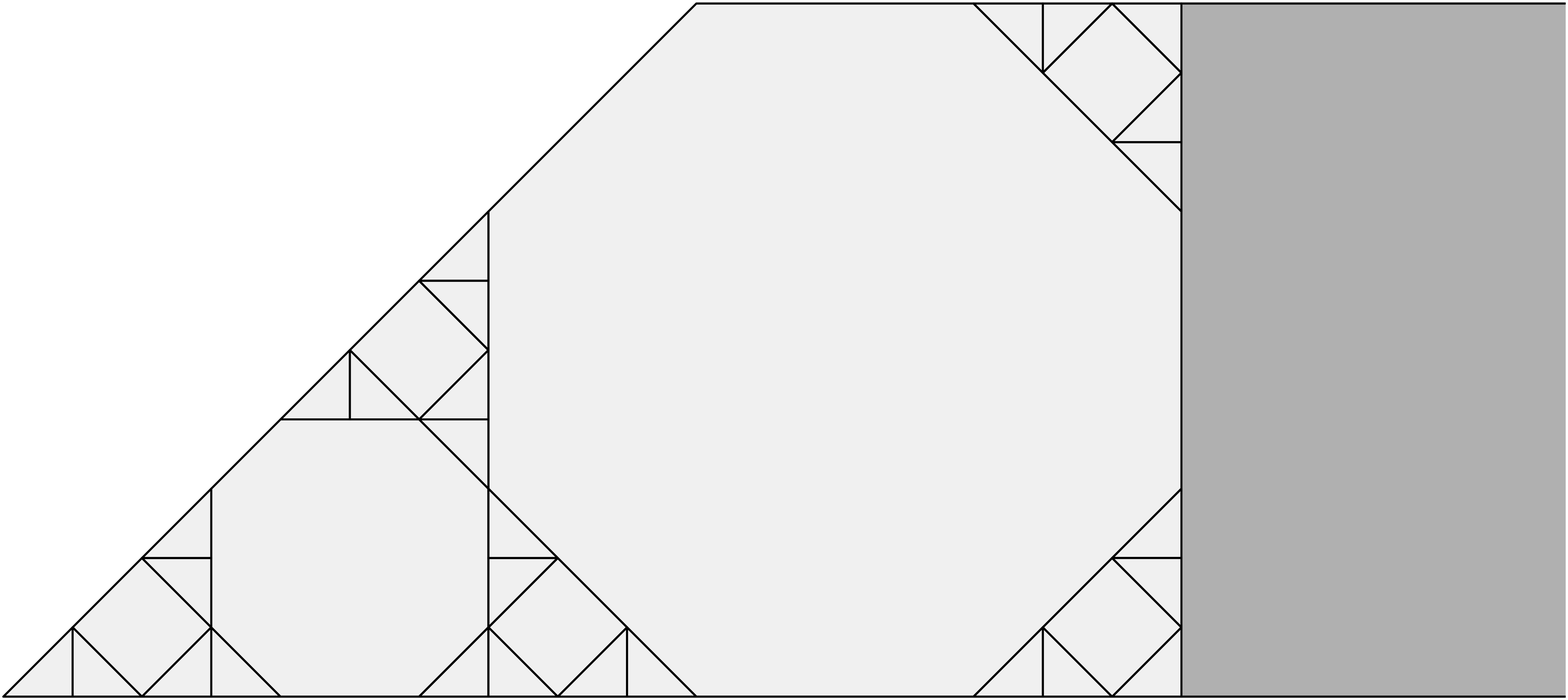}}
\newline
{\bf Figure 4.4:\/} 
The left half of $\Delta_t$ for $t=R(s)=5/17$.
\end{center}

By the previous result, and the Insertion Lemma, we
can find a counterexample $\Omega$ having parameter
$s \in (1/4,1/2)$ such that $\lambda(\Omega)$ is as
close as we like to $\lambda$.  

Consider a counterexample $\Omega$ having
$s \in (1/4,1/3)$.
For $s$ in this range, we have $t=R(s)>1/2$.
The blue octagon in $\Delta_s$
is the image under $\phi_s$, of the
trivial tile in $\Delta_t$.
Call this octagon $O_s$.

\begin{center}
\resizebox{!}{1.6in}{\includegraphics{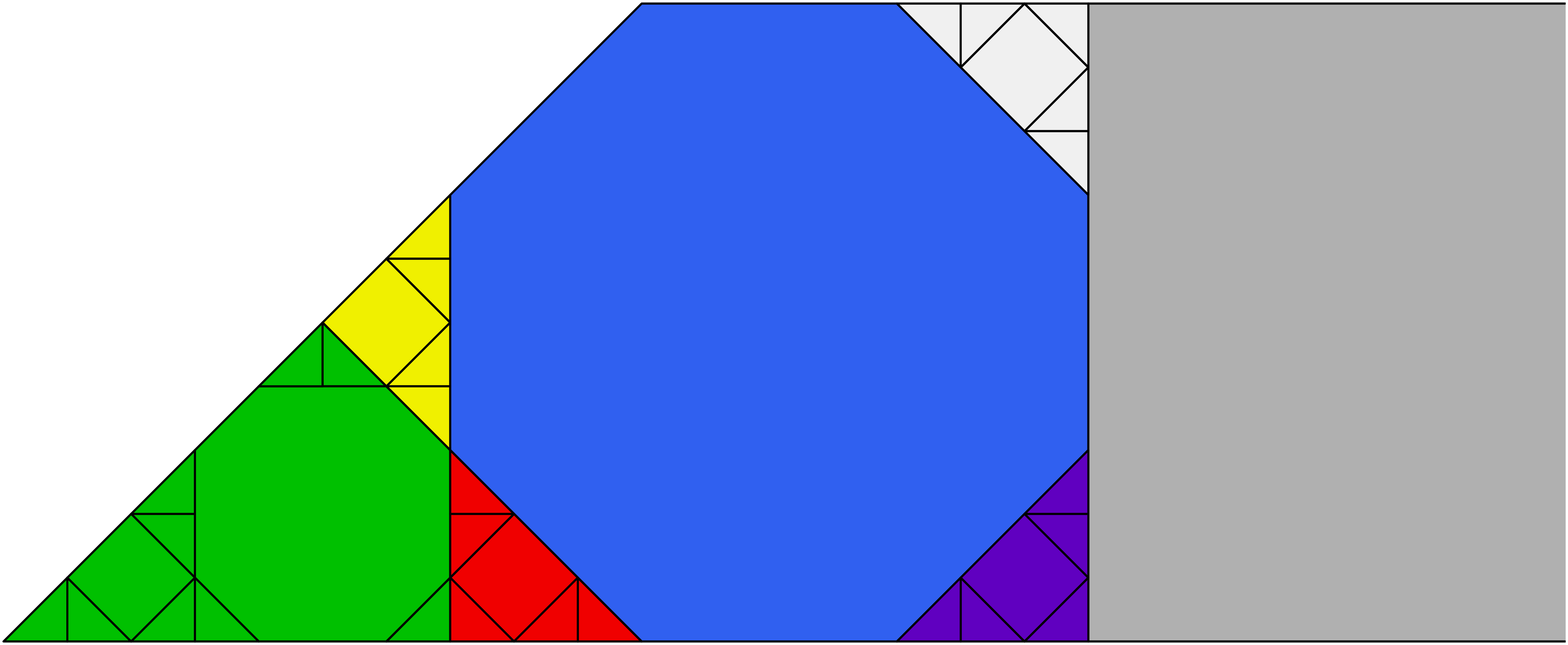}}
\newline
{\bf Figure 4.5:\/} 
$O_s$ (blue) and other regions for for $s=5/17$.
\end{center}

\begin{center}
\resizebox{!}{1.7in}{\includegraphics{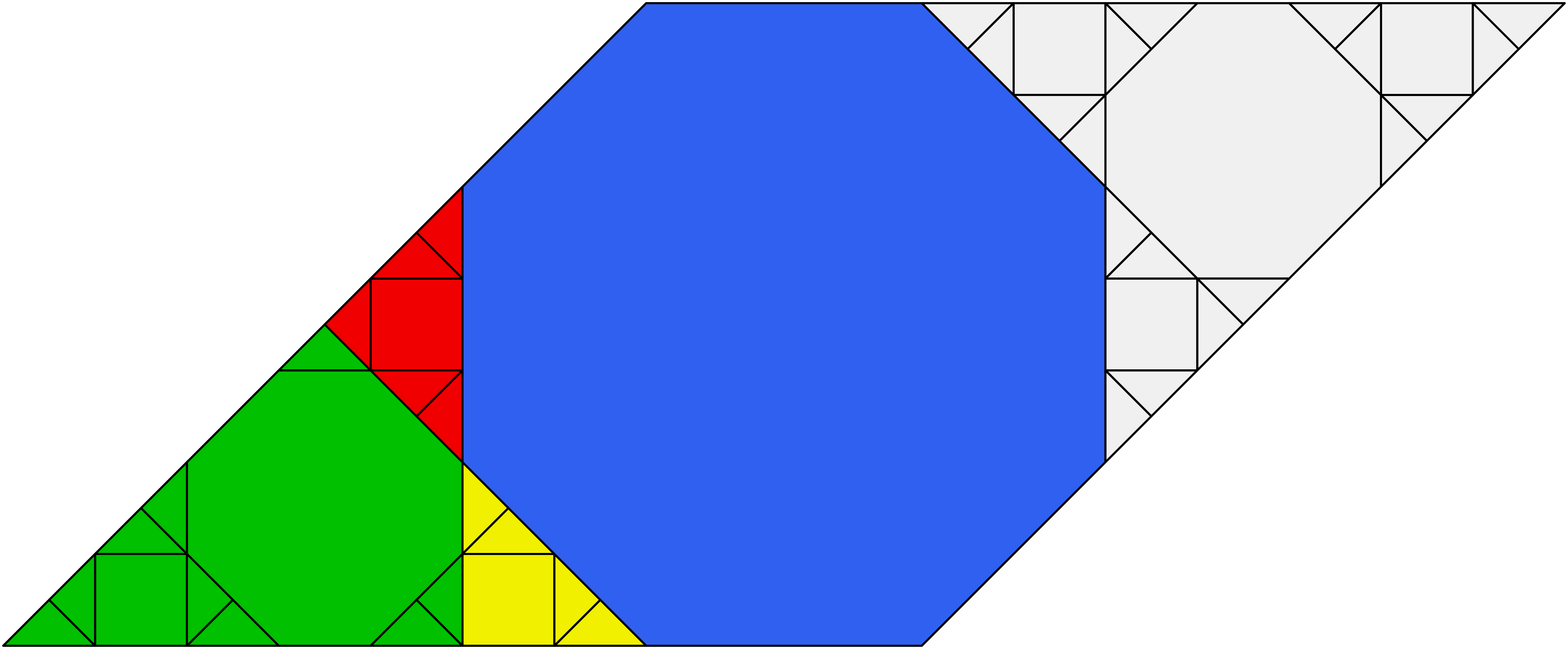}}
\newline
{\bf Figure 4.6:\/} 
$\Delta_t$ for $t=R(5/17)=7/10$.
\end{center}

We have the following facts.
\begin{itemize}
\item $\phi_s(Q_t)$ (yellow) is a
neighborhood of $H^0-O_s$.
\item $\phi_s(P_t)=B_s$ (red, green).
\item $R_H \circ  R_D \circ \phi_s(A_t)$ (purple)
is a neighborhood of $D_s-O_s$.
\item $\phi_s(B_t)=Q_s$ (green, yellow).
\end{itemize}
The colors refer to Figures 5.4 and 5.5.
In all cases, we can pull the counterexample
back by the relevant similarity to get a
counterexample associated to the parameter $t$.
Call the new counterexample $\Omega'$.
Since $Y_t^0=X_t^0$ for $t>1/2$, we have
$$\lambda(\Omega')=\frac{\|\phi_s^{-1}(p)-\phi_s^{-1}(q)\|}{{\rm diam\/}(Y_t^0)}=
\frac{\|p-q\|}{{\rm diam\/}(\phi_s(X_t^0))}=$$
\begin{equation}
\frac{\|p-q\|}{{\rm diam\/}(Z_s^0)}=
\lambda(\Omega)\frac{{\rm diam\/}(X_s^0)}{{\rm diam\/}(Z_s^0)}=
K_s \lambda(\Omega).
\end{equation}
Here the constant $K_s>1$ is uniformly bounded away from
$1$, by compactness.  But this is a contradiction.
We could choose $\lambda(\Omega)>\lambda/K_s$ and
then we would have $\lambda(\Omega')>\lambda$.

Now consider a counterexample with $s \in (1/2,1/3)$.
The argument is very similar to what we just did.
This time, we let $O_s$
be the pyramid from Lemma \ref{pyramid}.
In this case, the same remarks apply about the
points of the counterexample being disjoint from
(the interior of) $O_s$.
We have the following facts.
\begin{itemize}
\item $R_D \circ T_s^{j} \circ \phi_s(P_t)$ (red)
is a neighborhood of $H^0-O_s$ for a suitable 
choice of $j$.  Here $T_s^j$ is as in Equation \ref{tube}.
\item $T_s^{j} \circ \phi_s(Q_t)$ (blue)
is a neighborhood of $V-O_s$.
\item $R_V \circ \phi_s(A_t^0)$ (yellow)
is a neighborhood of $D-O_s$.
\item $E_s=\phi_s(B_t)$ (green).
\end{itemize}

\begin{center}
\resizebox{!}{2.4in}{\includegraphics{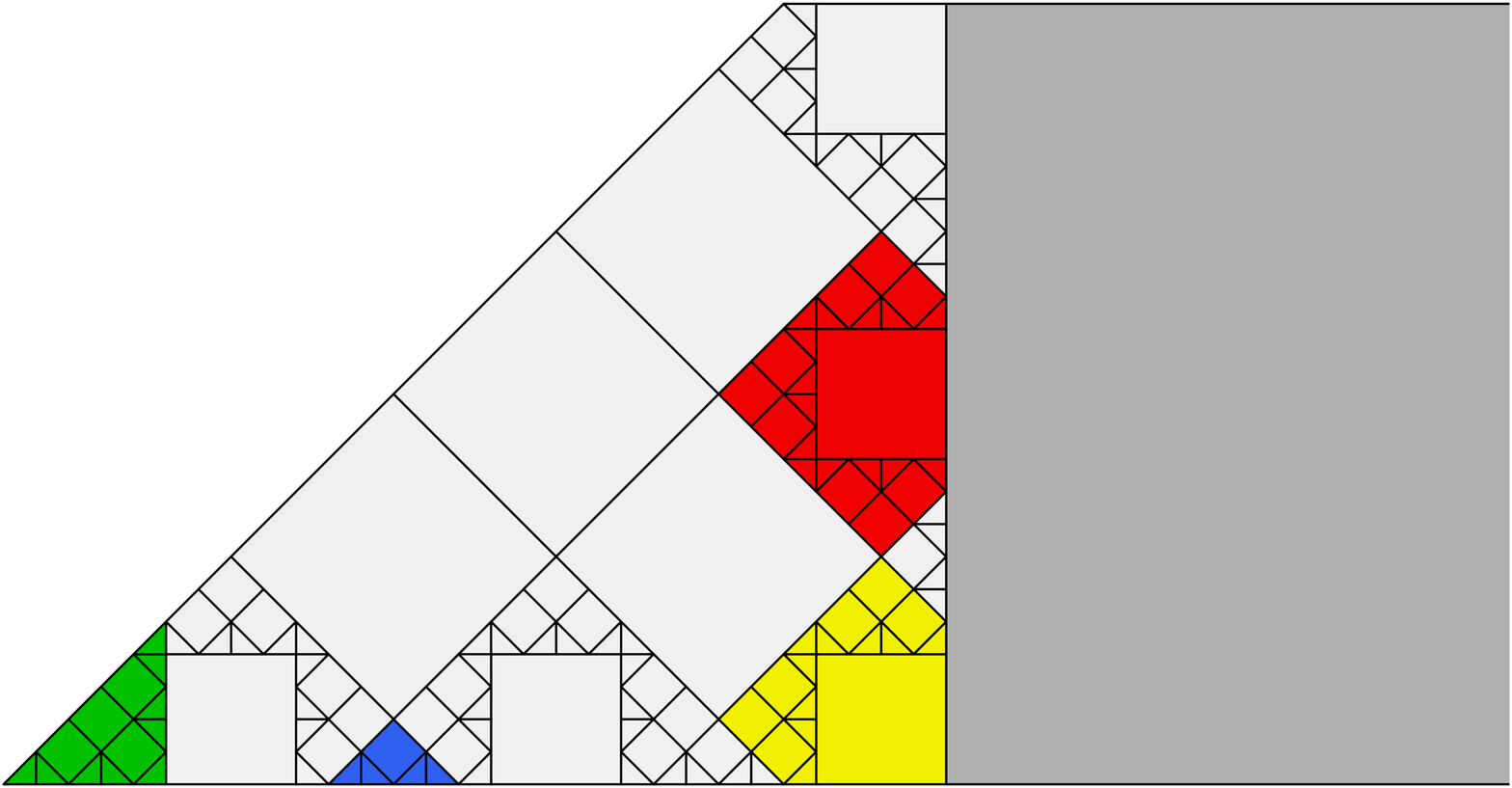}}
\newline
{\bf Figure 4.7:\/} 
Some regions of $\Delta_s$ for $s=12/29$.
\end{center}

\begin{center}
\resizebox{!}{1.6in}{\includegraphics{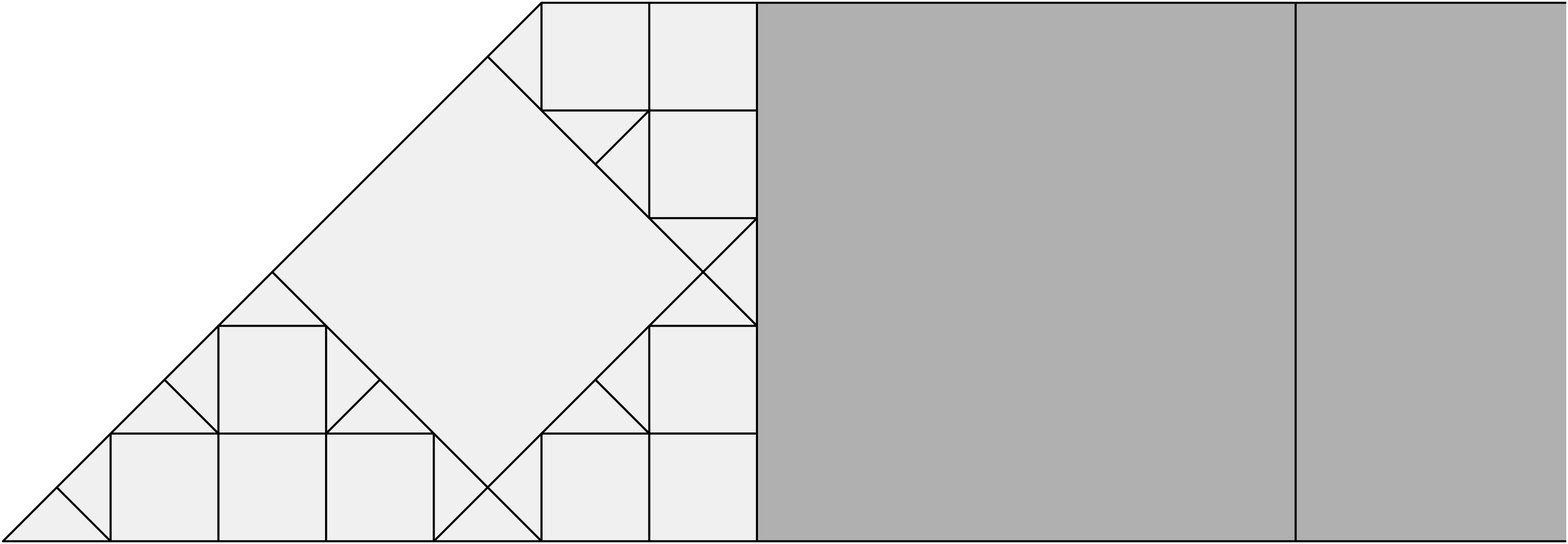}}
\newline
{\bf Figure 4.8:\/} 
The left half of $\Delta_t$ for $t=R(12/29)=5/24$.
\end{center}

This time, we get
\begin{equation}
\lambda(\Omega')=\lambda(\Omega)\frac{{\rm diam\/}(X_s^0)}{{\rm diam\/}(\phi_s(X_t^0))}=
L_s \lambda(\Omega)
\end{equation}
Here $L_s>1$ is uniformly bounded away from $1$.
We get the same contradiction as before.
This completes the proof of the Pinching Lemma.
\subsection{The Covering Lemma}

We say that an $\epsilon$-{\it patch\/} is a triple
$(K,\psi,u)$ where
\begin{itemize}
\item $u \in (0,1)$.
\item $K$ is one of the $4$ symmetric sets
$A_u,B_u,P_u,Q_u$ from \S \ref{symm}.
\item $\psi: K \to X_s$ is a similarity which
contracts by some factor $\lambda \leq \epsilon$.
\item $\psi(\Delta_u \cap K)=\Delta_s \cap \psi(K)$.
\end{itemize}
The last condition means that $\psi$ gives a bijection
between tiles of
$\Delta_u \cap K$ and tiles of $\Delta_s \cap \Psi(K)$.
When we have an $\epsilon$-patch $(K,\psi,u)$, we are
recognizing a small portion of $\Delta_s$ as being a similar
copy of a large portion of $\Delta_u$.  When the choice
of $\epsilon$ is not relevant to the discussion,
we will just say {\it patch\/} in place of
$\epsilon$-patch.  

We present two versions of our result.
In these versions, the constants $m$ and $n$
depend on everything in sight.  They are
meant to (possibly) vary from case to case.
Here is the first version.

\begin{lemma}
\label{cov0}
Suppose $s \in (0,1)$ is irrational, and
let $t=R(s)$.  Then for each 
set $K_s \in \{A_s,B_s,P_s,Q_s\}$, we have
$$K_s = \bigcup_{i=1}^m \alpha_i \cup
\bigcup_{j=1}^n \beta_j,$$
where $\alpha_i$ is a tile of
$\Delta_s$ and $\beta_j$ is the image
of some $\epsilon$-patch
$(K_j,\psi_j,t)$.  Here
$\epsilon$ is the scale factor of $\phi_s$,
the map from Theorem \ref{renorm}.
\end{lemma}

Iterating Lemma \ref{cov0} for the sequence
$\{R^n(x)\}$ and using the fact that
$X_s^0=A_s \cup B_s=P_s \cup Q_s$, we
get the following corollary.

\begin{corollary}[Covering]
\label{cov1}
For any $\epsilon>0$ we have
\begin{equation}
\label{patcov}
\Delta_s = \bigcup_{i=1}^m \alpha_i \cup
\bigcup_{j=1}^n \beta_j,
\end{equation}
where $\alpha_i$ is a tile of
$\Delta_s$ and $\beta_j$ is the image
of some $\epsilon$-patch
$(K_j,\psi_j,u)$.  Moreover, we
can take $u=R^k(s)$ for any sufficiently large $s$.
Equation \ref{patcov} likewise
holds for $\Delta_s \cap K_s$ in place of $\Delta_s$ for
each symmetric piece $K_s$.
\end{corollary}

Lemma \ref{cov0} has the same kind of proof as the Pinching Lemma.
Let $K=K(s)$ be the layering constant for $s$ which
appears in the Filling Lemma.

\begin{lemma}
Lemma \ref{cov0} is true when $s>1/2$ and
$K(s)>1$.
\end{lemma}

\startproof
Figure 4.9 illustrates a typical picture for $s$
in this range.
\begin{center}
\resizebox{!}{1.8in}{\includegraphics{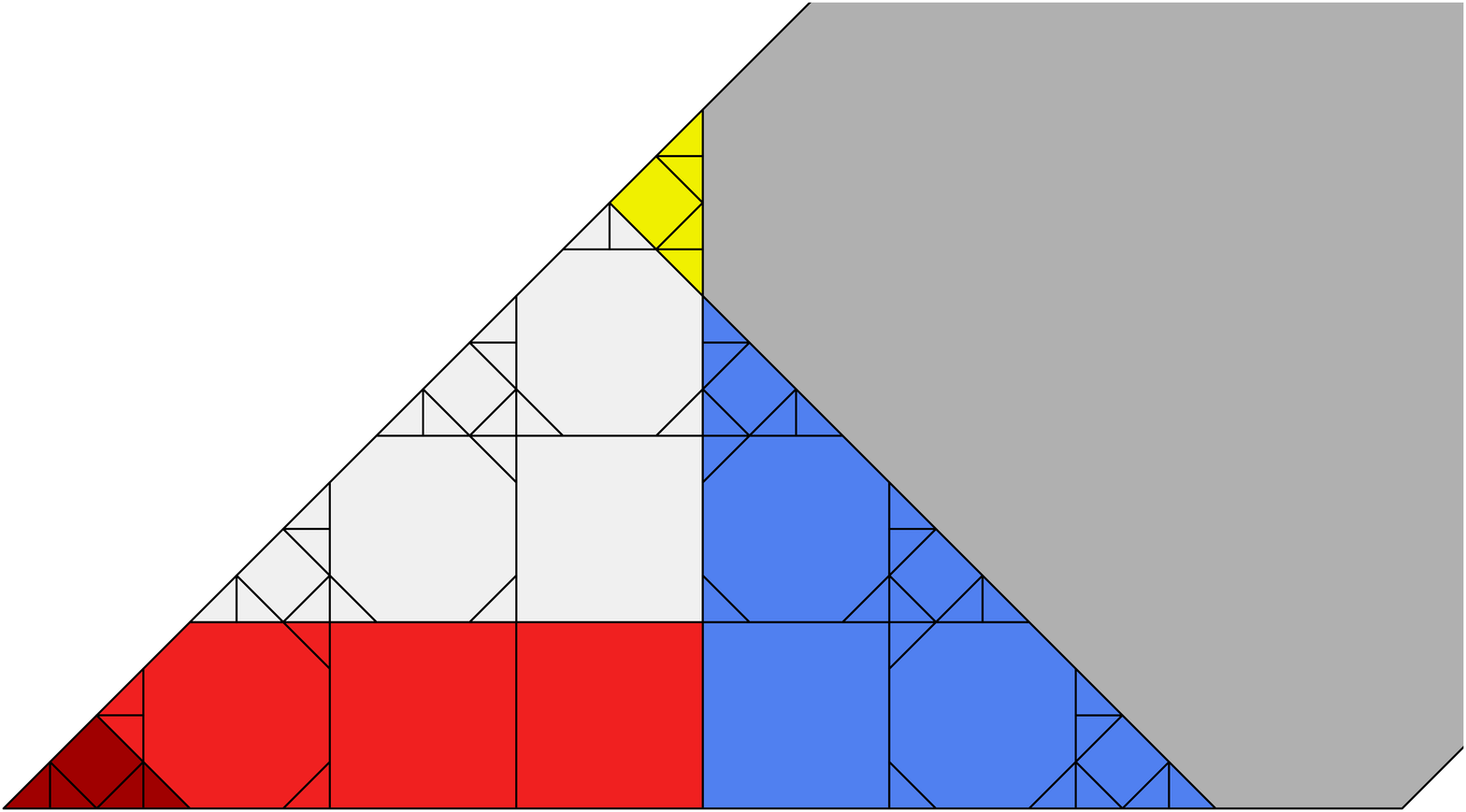}}
\newline
{\bf Figure 4.9:\/} The tiling $\Delta_s^0$ for $s=13/15$.
Here $K(s)=3$.
\end{center}

We have the following equations.
\begin{itemize}
\item $A_s=R_D \circ \phi_s(B_t)$.
$A_s$ is colored yellow in Figure 4.9, and
$\phi_s(B_t)$ is colored dark red.
$R_D$ is reflection in the diagonal line $D=D_s$ of symmetry.
\item $$B_s=\nu \cup R_V(\nu) \cup R_D \circ \psi_s 
\circ(A_t) \cup \Theta, \hskip 30 pt
\nu=\bigcup_{k=0}^{K-1} \Psi_s^k.$$
The set $B_s$ is colored red/white/blue in Figure 4.9.
Here $V$ is the vertical line of symmetry,
and $\Theta$ is a finite union of square tiles.
\item $P_s=\nu \cup R_D(\nu)$. In
Figure 4.9, 
$P_s$ is colored red/white/yellow.
\item $$Q_s=R_D(\nu'), \hskip 30 pt
\nu'=T_s^{-1}\bigg(\bigcup_{k=1}^{K-1} \Psi_s^k\bigg).$$
Here $T_s$ is as in Equation \ref{tube} and $Q_s$ is
colored blue in Figure 4.9.
\end{itemize}
We can interpret all these equations as the desired
patch coverings.
\endproof

\begin{lemma}
Lemma \ref{cov0} is true when $s>1/2$.
\end{lemma}

\startproof
It only remains to treat the case when $K(s)=1$.
Figure 4.10 shows a representative picture
for $s$ in this range.

\begin{center}
\resizebox{!}{2.2in}{\includegraphics{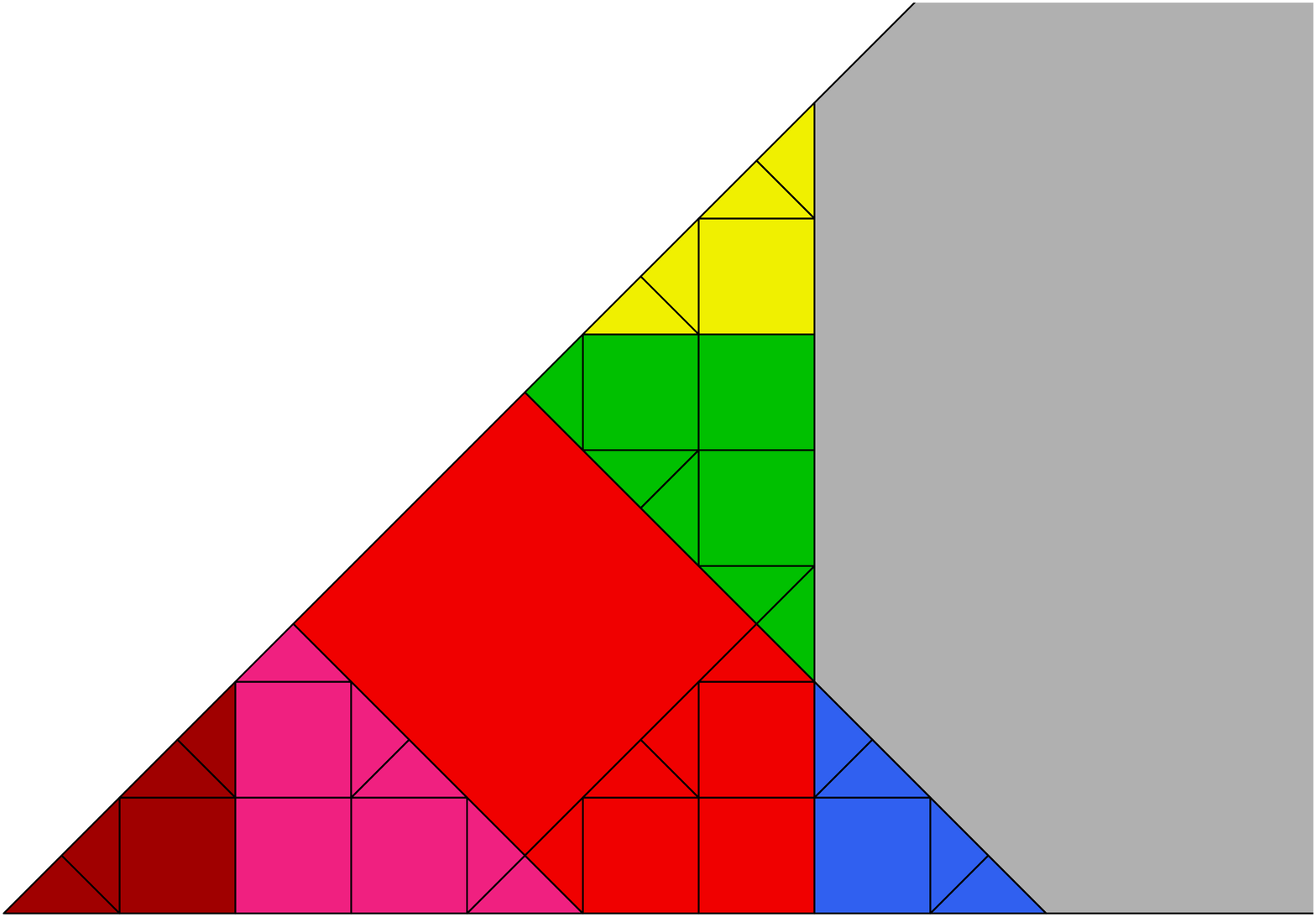}}
\newline
{\bf Figure 4.10:\/} The tiling $\Delta_s^0$ for $s=9/14$.
Here $K(s)=1$.
\end{center}

\begin{itemize}
\item $A_s=R_D \circ \phi_s(B_t)$. (This is as above.)
In Figure 4.10, $A_s$ is colored yellow/green, and
$\phi_s(B_t)$ is colored dark red/pink.
\item $$B_s=R_d \circ \phi_s(A_t) \cup \phi_s(Q_t) \cup
R_V \circ \phi_s(Q_t).$$
$B_s$ is the union of all the pieces not colored \footnote{The
grey tile, which is a central tile, is always excluded.}
yellow or green.
The first set on the right is colored red/pink,
the second set is colored dark red, and
the third set is colored blue.
\item $$P_s=\phi_s(P_t) \cup \phi_s(Q_t) \cup R_D \circ \phi_s(Q_t).$$ 
$P_s$ is the union of all the pieces not colored blue.
The first set on the right is colored pink/red/green, and the
second set is colored dark red, and the third set is colored yellow.
\item $Q_s=\phi_s(Q_t)$.  This set is colored dark red.
\end{itemize}
We can interpret all these equations as the desired
patch coverings.
\endproof

\begin{lemma}
Lemma \ref{cov0} is true when $s<1/2$ and $K(s)>1$.
\end{lemma}

\startproof
Figure 4.11 shows a representative picture
for $s$ in this range.

\begin{center}
\resizebox{!}{1.8in}{\includegraphics{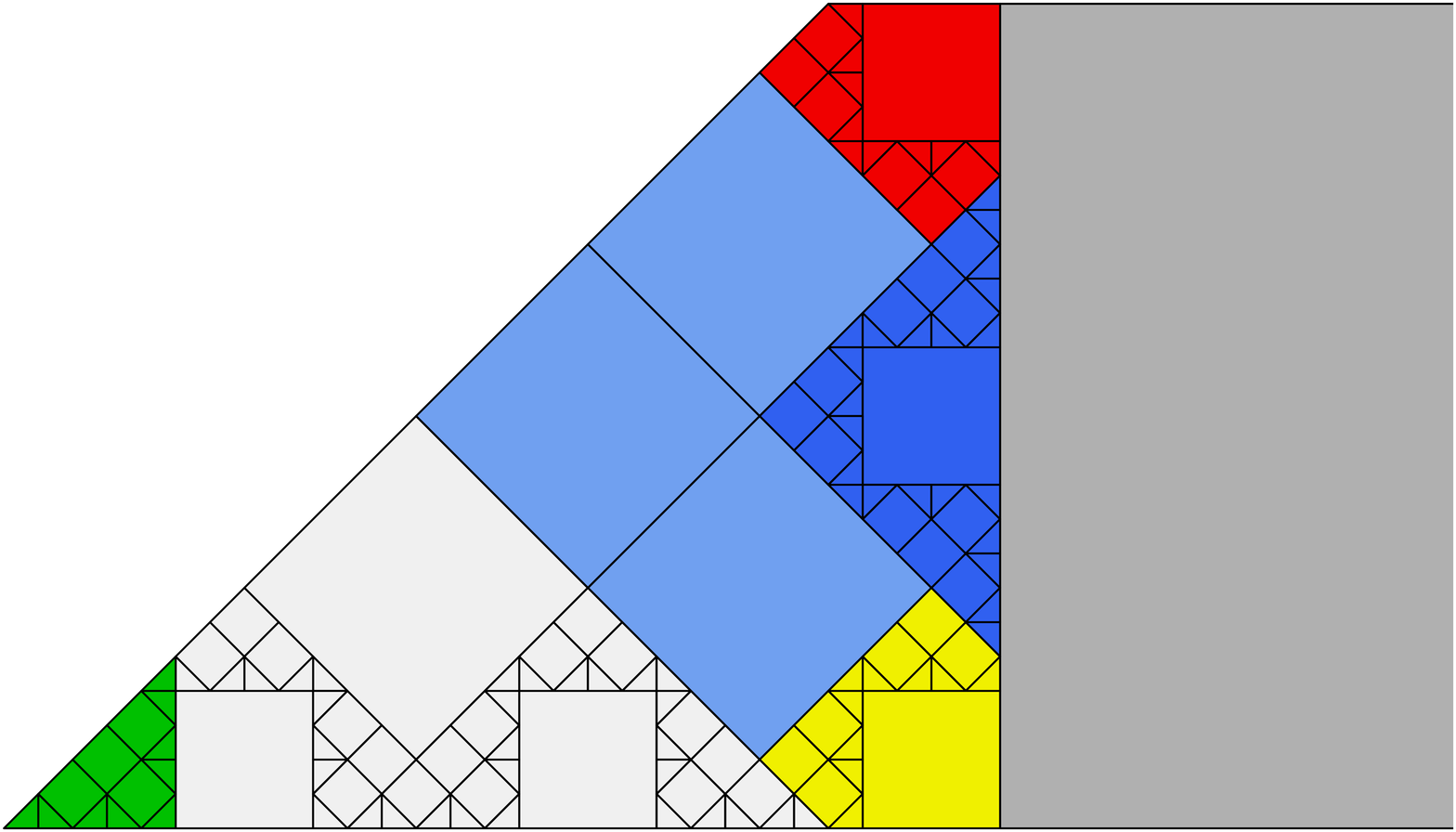}}
\newline
{\bf Figure 4.11:\/} The tiling $\Delta_s^0$ for $s=12/29$.
Here $K(s)=1$.
\end{center}

\begin{itemize}
\item $A_s=R_D \circ \phi_s(A_t)\ \cup
\ R_H \circ R_D \circ \phi_s(A_t)\ \cup
\ R_D(\nu')\ \cup\ R_H \circ R_D(\nu') \cup \Theta,$
$$ \hskip 30 pt
\nu'=\bigcup_{k=1}^{K-1} \Psi_s^k.$$
Here $\Theta$ is a finite union of square tiles,
colored light blue.
The set $A_s$ is the union of tiles in $X_s^0$
not colored white or green (or grey). 
The first set on the right is colored red.
The second set is colored yellow. The union of
the third and fourth sets is colored dark blue.
\item $$B_s=\nu \cup R_V(\nu)-\Theta, \hskip 30 pt
\nu=\bigcup_{k=0}^{K-1}\Psi_s^k.$$
Here $\Theta$ is a finite union of squares
(colored light blue) which we delete from
our union.  The point is that the top square
of each $\Psi_s^k$ lies above $B_s$.
$B_s$ is colored white and green.
\item $P_s=A_s \cup R_D(A_s)$.
\item $Q_s=\phi_s(B_t)$.
\end{itemize}
We can interpret all these equations as the desired
patch coverings.
\endproof

\begin{lemma}
Lemma \ref{cov0} is true when $s<1/2$.
\end{lemma}

\startproof
In view of what we have already shown, it suffices
to consider the case when $K(s)=1$. Figure 4.12 shows
a representative picture in this case. 

\begin{center}
\resizebox{!}{2.3in}{\includegraphics{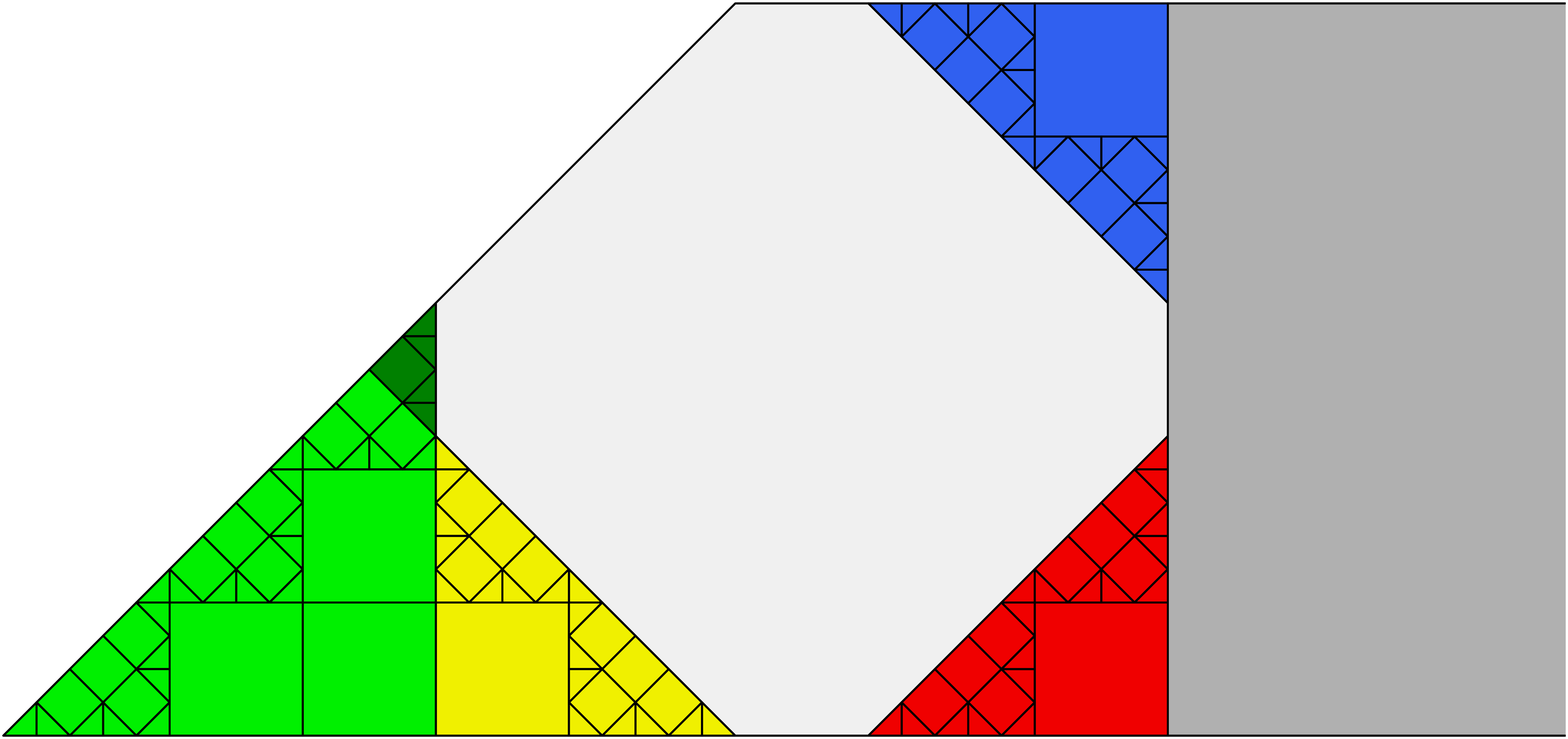}}
\newline
{\bf Figure 4.12:\/} The tiling $\Delta_s^0$ for $s=11/35$.
Here $K(s)=1$.
\end{center}

\begin{itemize}
\item $$A_s=R_D \circ \phi_s(A_t) \cup
R_H \circ R_D \circ \phi_s(A_t) \cup \phi_s(Q_t) \cup \Theta.$$
$A_s$ is colored red/white/blue/green.
$\phi_s(A_t)$ is colored yellow. The sets
on the right are respectively colored blue, red, green,
and white.
\item $B_s=\phi_s(P_t)$. This set is colored yellow and light green.
\item $$P_s=\phi_s(A_t) \cup R_D \circ \phi_s(A_t) \cup
R_H \circ R_D \circ \phi_s(A_t) \cup \Theta.$$
Here $P_s$ is the union of all tiles not colored green (or grey).
The sets on the right are, respectively, colored yellow, blue,
red, and white.
\item $Q_s=\phi_s(B_t)$.  This set is colored light and dark green.
\end{itemize}
We can interpret all these equations as the desired
patch coverings.
\endproof

We have exhausted all the cases. This completes
the proof of Lemma \ref{cov0}.

\newpage

\section{Proof of Statement 1}

\subsection{The Easy Direction}
\label{easy}

Our goal is to prove that $\widehat \Lambda_s$ is a
disjoint union of two arcs 
if and only if $s$ is oddly even.
Let $S_s$ denote the left half of
$\widehat \Lambda_s$.  By symmetry, the result we want is
equivalent to the statement that $S_s$ is an arc if and only
if $s$ is oddly even.  

First we prove that $S_s$ is an arc only if $s$ is oddly even.
This is the easier of the two directions.

\begin{lemma}
\label{wedge}
For each integer $k$ such
that $R^k(s)>1/2$, there is an octagon $O_k$ which
having one edge in the left side of $X_s$ and
one edge in the bottom side of $X_s$.
If there are two istinct indices $k$ and $\ell$ with
this property, then the octagons $O_k$ and
$O_{\ell}$ are distinct.
\end{lemma}

\startproof
Say that an octagon is {\it wedged\/} in a
parallelogram (of the kind we are considering)
if one edge of the octagon lies
in the left edge of the parallelogram and
another edge lies in the bottom edge.

Let $s_0=s$ and
$s_k=R^k(s)$.  For ease of notation, we
set $\phi_k=\phi_{s_k}$ and $\Delta_k=\Delta_{s_k}$, etc.
When $s_k>1/2$, the central tile $C_k$
of $\Delta_k$ is an octagon wedged into $X_k$.
By Theorem \ref{renorm}, the octagon
$\phi_{k-1}(C_k)$ is a tile of $\Delta_{k-1}^0$, and
is wedged into $X_{k-1}$.  

Iterating Theorem \ref{renorm}, we see that
\begin{equation}
O_k=\phi_0 \circ ... \circ \phi_{k-1}(C_k)
\end{equation}
is wedged into $X_0$. 

Suppose that $\ell>k$ is another index such that
$s_{\ell}>1/2$.  Then the two octagons
$$C_k, \hskip 30 pt
\phi_k \circ ... \circ \phi_{\ell-1}(C_{\ell})$$
are distinct because one octagon is the central
tile of $\Delta_k$ and the other one is not.
But then $O_k$ and $O_{\ell}$ are the images of
the above octagons under the same similarity.
Hence, they are distinct.
\endproof

\begin{corollary}
\label{wedge2}
If $R^n(s)>1/2$ for at least $K$ different positive
indices, then $S_s$ has at least $K+1$ connected
components.
\end{corollary}

\startproof 
The $K$ 
octagons guaranteed by Lemma \ref{wedge}
are all distinct.  Call these octagons
$O_1,...,O_K$.  Each of these octagons is
wedged into $X_s$, and so the union of
these octagons separates $X_s^0$ into
$K+1$ connected components.   We just need
to see that $S_s$ intersects each component.

Each octagon $O_j$ has two vertices in the bottom edge
of $X_s$.  At each of these vertices, the adjacent
edge of $O_j$ makes an acute angle with the bottom edge
of $X_s$.  (The angle is $\pi/4$.)  Since the 
$\Delta_s$ consists of an open dense (in fact full
measure) set of squares and semi-regular octagons,
every neighborhood of the two vertices in question
must intersect infinitely many tiles of
$\Delta_s$.  Hence, the two bottom vertices of
$O_j$ lie in $S_s$.  This proves what we want.
\endproof

What we have shown is that $S_s$ is not an arc if
$R^n(s)>1/2$ for some $n>0$.  

\begin{lemma}
If $s>1/2$, then $S_s$ is not an arc.
\end{lemma}

\startproof
Let $C$ be the trivial tile of $\Delta_s$.
The same argument as in the previous lemma shows that
$S_s$ contains the following three points.
\begin{enumerate}
\item The vertex where the vertical edge of $C$
meets the left edge of $X_s$.  This is the
top vertex of both $A_s$ and $P_s$.
\item The vertex where the slope $-1$ edge of $C$
meets the bottom edge of $X_s$.  This is the
right vertex of $B_s$ and of $Q_s$.
\item The bottom left vertex of $X_s$.  Thie is
the left vertex of $B_s$ and of $P_s$.
\end{enumerate}
Any arc connecting these $3$ points, in some
order, must cross at least twice one of the
symmetry lines from the Pinching Lemma.  Hence
$S_s$ cannot be an arc.
\endproof

\noindent
{\bf Remark:\/}
In fact $S_s$ is homeomorphic to a ''Y'' when $s>0$
and $R^n(s)<1/2$ for all $n>0$. However, we do not
prove this.
\newline

We now know that $S_s$ is an arc only if $s$ is
oddly even.  The rest of the chapter is devoted
to proving that $S_s$ is an arc when $s$ is oddly even.

\subsection{A Criterion for Arcs}

Say that a {\it marked piece\/} is an compact, embedded,
convex set with two distinguished vertices.
Say that a {\it chain\/} is a finite union
$D_1,...,D_n$ of marked pieces such that
$D_i \cap D_{i+1}$ is one point, and that this point
is one of the marked points on each of $D_i$ and $D_{i+1}$.
We also require that $D_i \cap D_j=\emptyset$ for
all other indices $i \not = j$.
We define the {\it mesh\/} of the chain to be the
maximum diameter of one of the marked pieces.

We say that a compact set
$S$ {\it fills\/} a chain $D_1,...,D_n$ if
$S \subset \bigcup D_i$ and $S$ contains 
every marked point of the chain.   The purpose
of this section is to establish the following
(certainly well known) criterion.

\begin{lemma}[Arc Criterion]
Let $S$ be a compact set.
Suppose, for every $\epsilon>0$, that $S$ fills a
chain having mesh less than $\epsilon$.
Then $S$ is an embedded arc.
\end{lemma}

We will assume that $S$ satisfies the hypotheses
of the lemma, and then show that $S$ is an arc.
First of all, $S$ is clearly connected. 

\begin{lemma}
Suppose that $S$ fills a chain $C_1,...,C_m$.
Then there is some $\epsilon>0$ with the
following property. If $S$ also fills a
chain $D_1,...,D_n$ having mesh size less
than $\epsilon$, then $S$ fills a chain
$E_1,...,E_p$ where each $E_i$ has the form
$C_j \cap D_k$.
\end{lemma}

\startproof
We can choose $\epsilon$ so small that each $D_j$ has
following properties.
\begin{itemize}
\item The diameter of $D_j$ is smaller than
the length of any edge of any $C_i$.
\item $D_j$ intersects any $C_i$ in at most $2$ edges.
\item $D_j$ cannot intersect $C_i$ and $C_k$ if
$i$ and $k$ are not consecutive indices.
\end{itemize}
For each $i$, there are unique and distinct
pieces $D_{j_1}$ and $D_{j_2}$ which
contain the two marked points of $C_i$.
We claim that
the pieces between $D_{j_1}$ and $D_{j_2}$
must have both marked points inside $C_i$.
Assuming that this is true, we form the
portion of the $E$-chain insid $C_i$ by taking the
intersections $C_i \cap D_j$ and using the
marked points of $D_j$ for $j_1<j<j_2$.
The marked
points of $C_i \cap D_{j_1}$ are the
marked point of $C_i$ inside of $D_{j_1}$
and the marked point of $D_{j_1}$ inside
$C_i$.  Similarly for $C_i \cap D_{j_2}$.
We do the same thing for each $i$ and this
gives us the conclusion of the Lemma.

Now we establish our claim.
If our claim was false, then some
$D_j$, with $j_1<j<j_2$, would have
one marked point in $C_i$.  Note that 
$D_j$ must have another marked point in
either $C_{i-1}$ or $C_{i+1}$, because
this marked point is a vertex of $S$.
Suppose that $D_j$ has its other
marked point in $C_{i+1}$.
Then $D_j \cap (C_i \cup C_{i+1})$ is
disconnected because $D_j$ does
not contain the vertex $C_i \cap C_{i+1}$.
Hence
$$D_1 \cup ... \cup D_{j-1} \cup
\big(D_j \cap (C_i \cup C_{i+1})\big) \cup
D_{j+1} \cup ... \cup D_n$$
consists of two disconnected components,
each of which intersects $S$ nontrivially.
This contradicts the connectivity of $S$.
\endproof

If the chain $C_1,...,C_m$ and
the chain $E_1,...,E_p$ are related
as in the previous lemma, we say that
$E_1,...,E_p$ {\it refines\/}
$C_1,...,C_m$.
In view of the previous result,
we can assume that $S$ fills an infinite
sequence $\{\Omega_i\}$ of chains
such that each one refines the previous one
and the mesh size tends to $0$.

For each $i$, we inductively
create a partition $P_i$ of $[0,1]$ into
intervals, such that the number of intervals coincides
with the number of marked pieces in $\Omega_i$,
in the following manner.
Once $P_i$ is created, we distribute the intervals of
$P_{i+1}$ according to how $\Omega_i$ contains $\Omega_{i+1}$.
If the $k$th piece of $\Omega_i$ contains $n_k$ pieces
of $\Omega_{i+1}$, then $P_{i+1}$ is created from $P_i$
by subdividing the $k$th interval of $P_i$ into $n_k$
intervals of equal size.  Note that the mesh
size of $P_i$ tends to $0$ as $i$ tends to $\infty$.

There is a bijective correspondence between marked
pieces in the chains and intervals in the partition.
The correspondence respects the containment and
intersection properties.  For instance, two marked
pieces intersect if and only if the corresponding
intervals share an endpoint.  Each point of $S$
is contained in an infinite nested intersection
of marked pieces, and we map this point to the
corresponding nested intersection of intervals.
This map is clearly a homeomorphism.  The inverse
map gives a parameterization of $S$ as an arc in
the plane.

\subsection{Elementary Properties of the Limit Set}

Let $A_s$ be the set from \S \ref{symm}.

\begin{lemma}
\label{fill00}
$S_s$ contains the two left vertices of $X_s$ and
the two obtuse vertices of $A_s$.
\end{lemma}

\startproof
Let $v$ be one of left vertices of $X_s$.
Since the angle of $X_s$ at $v$ is not a right
angle, there must be infinitely many squares
contained in every neighborhood of $v$.

Note that the top left vertex of $X_s^0$ is
also the top obtuse vertex of $A_s$.  So,
$S_s$ contains the top obtuse vertex of
$A_s$.  By symmetry, $S_s$ contains the
bottom obtuse vertex of $A_s$.
\endproof

By Lemma \ref{box}, $\Delta_s$ consists
entirely of squares.  As we remarked after
proving Lemma \ref{box}, these squares are
either boxes or diamonds.

\begin{lemma}
\label{cool}
Suppose $\gamma \subset X_s^0$ is a compact 
arc which connects a point
in a box to a point in a diamond.
Then $\gamma$ contains a point of $S_s$.
\end{lemma}

\startproof
Compare [{\bf S0\/}, \S 8].
By compactness, it suffices to show
that arbitrarily small perturbations of
$\gamma$ contain points of $S_s$.
Hence, we may perturb so that
$\gamma$ does not contain any vertices of
any tiles in $\Delta_s$. 
Suppose $\gamma$ does not
intersect $\widehat \Lambda$.  Then
$\gamma$ only intersects finitely many tiles,
$\tau_1,...,\tau_n$.  Moreover,
$\tau_i$ and $\tau_{i+1}$ must share
an edge.  Hence, by induction,
$\tau_1$ is a box if and only if
$\tau_n$ is a box.  But $\tau_1$ is
a box and $\tau_n$ is a diamond.
This is a contradiction.
\endproof

\begin{lemma}
\label{cool2}
Each fundamental line of symmetry
contains a point of $S_s$.
\end{lemma}

\startproof
To make the argument cleaner, we attach a
large diamond $\delta$ to the picture along
the left edge of $X_s$, and we attach a
large box $\beta$ to the picture along the
bottom edge of $X_s$.  These extra squares
are disjoint from $X_s$ except along the
relevant edges.
Once we add these
two squares, we see that each of the lines 
in question connects a diamond to a box.
$H$ connects $\delta$ to a the
leftmost central tile of $\Delta_s$ and
both $V, D, E$ all connect $\beta$ to $\delta$.
By Lemma \ref{cool}, each of these lines
contains a point of $S_s$.
\endproof

\subsection{The End of the Proof}

Suppose that $S_s$ fills some chain 
$D_1,...,D_n$.  We call this chain {\it good\/}
\begin{itemize}
\item $D_j$ is disjoint from the interiors of 
the edges of
$\partial A_s$, for $j=1,...,n$.
\item The first marked point of $D_1$ is the bottom
left vertex of $X_s$.
\item The last marked point of $D_n$ is the top
left vertex of $X_s$.
\end{itemize}

\begin{lemma}
\label{good3}
$S_s$ fills a good chain.
\end{lemma}

\startproof
Our chain has two pieces.  We set $D_1=B_s$,
the triangle from \S \ref{symm},
and we let $D_2$ be the set from
Corollary \ref{good2}.
\endproof

\begin{lemma}
\label{good4}
Let $t=R(s)$.  Suppose
$S_t$ fills a good chain having
mesh $m$.  Then $S_s$ fills a good
chain having mesh at most $m/\sqrt 2$.
\end{lemma}

\startproof
Let $\Omega_t$ be the good chain filled
by $S_t$. We use the notation from
the Filling Lemma, Equation \ref{tube},
and Theorem \ref{renorm}.  We make our
construction in $5$ steps.
\newline
\newline
{\bf Step 1:\/}
For $j=0,...,K$, we define
\begin{equation}
\Omega_j=T_s^j \circ \phi_s(\Omega_t).
\end{equation}
Figure 5.1 illustrates our construction.
The individual chains
$\Omega_0,...,\Omega_{K-1}$ piece together to
make one long chain because the second disk
of $\Omega_j$ touches the common edge between
$\Psi_s^j$ and $\Psi_s^{j+1}$ only at the
bottom vertex of this edge.

\begin{center}
\resizebox{!}{2.5in}{\includegraphics{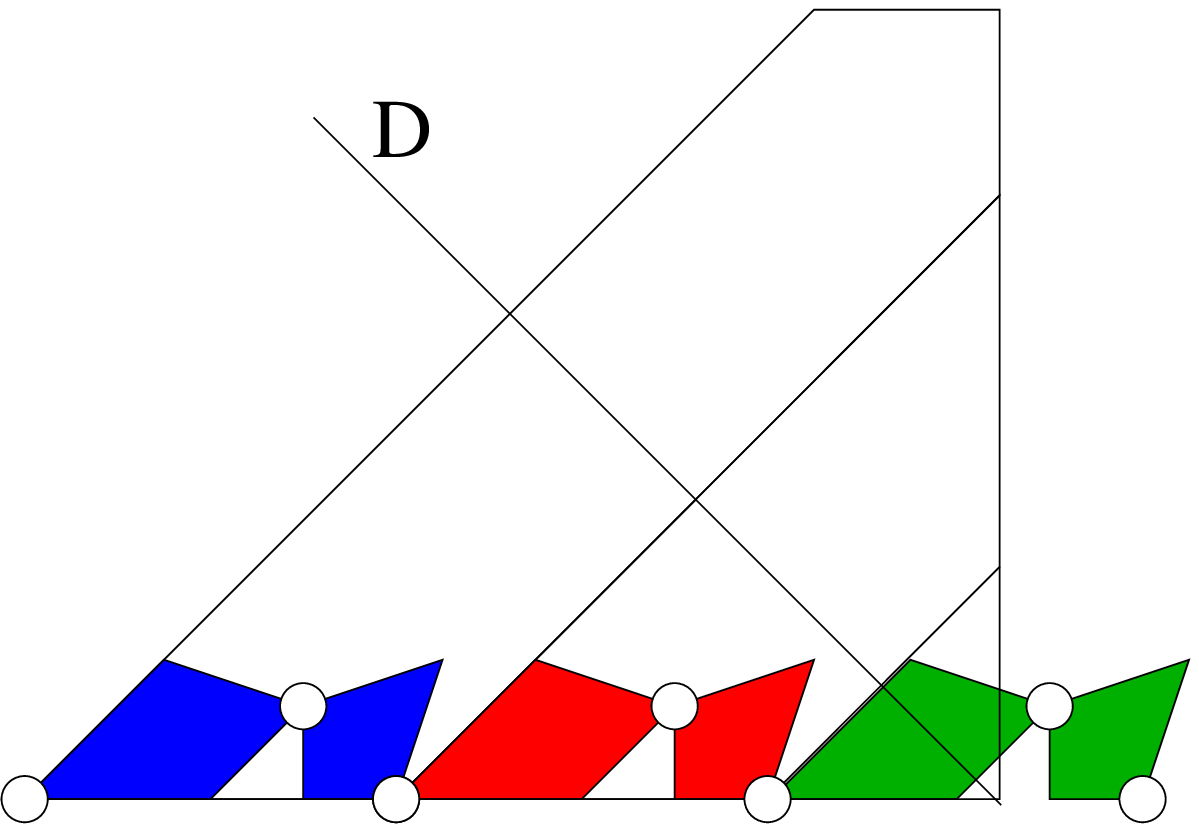}}
\newline
{\bf Figure 5.1:\/} Step 1: The chains
$\Omega_j$ for $j=0,...,K$. Here $K=2$.
\end{center}

\noindent
{\bf Step 2:\/}
The problem with $\Omega_K$ is that some of it sticks
over the edge of $X_s^0$.  This is the green set in
Figure 5.1.   However, we know from the Pinching Lemma
that $S_s$ intersects the line $D_s$ in a single point.
All other points of $D_s$ must have neighborhoods
contained in finitely many squares.  For this reason,
we can make the essentially the same construction
as in Corollary \ref{good2} to produce a convex disk
$U \subset \Psi_s^K$ such that 
$S_s \cap \Psi_s^K \subset U$ and
$U \cap D_s$ is the single point which
belongs to $S_s$.  

We now improve $\Omega_K$ as follows.  We intersect
each piece of $\Omega_K$ with the set $U$ and
throw out all those after the first one which is
disjoint from $U$.

\begin{center}
\resizebox{!}{2.6in}{\includegraphics{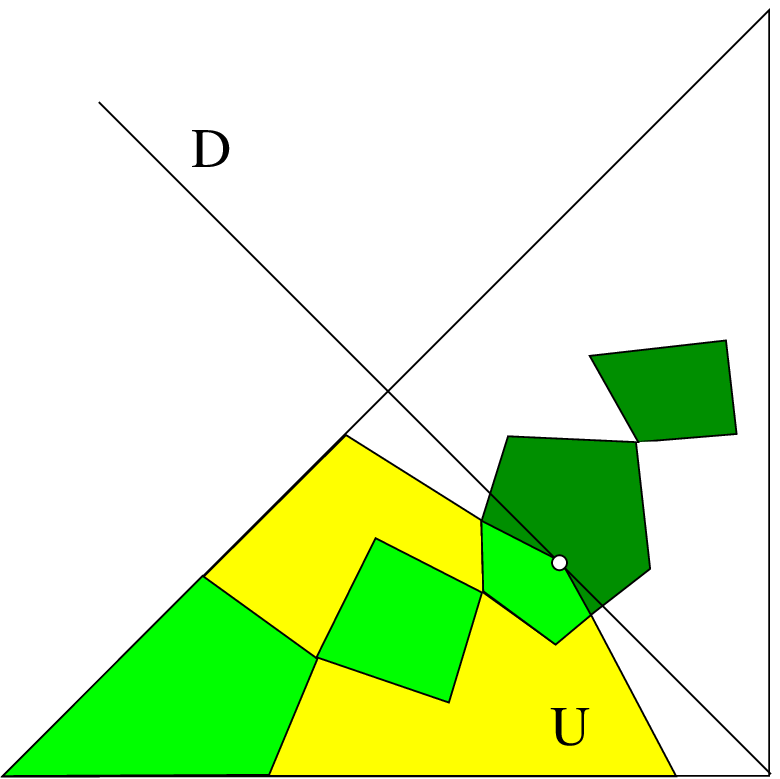}}
\newline
{\bf Figure 5.2:\/} Step 2: Improving $\Omega_K$.
\end{center}
  The result is a chain which joins
the bottom vertex of the edge 
$\Psi_s^{K-1} \cap \Psi_s^k$
to the point $S_s \cap D_s$.
Figure 5.2 shows the construction.
We call this improved chain $\Omega_K'$.

Define
\begin{equation}
\Upsilon_0 = \Omega_0,...,\Omega_{K-1},\Omega_K'.
\end{equation}
By construction, this chain is filled by
the portion of $S_s$ beneath the line $D_s$.
\newline
\newline
{\bf Step 3:\/}
Define
\begin{equation}
\Upsilon_2=\Upsilon_0,\Upsilon_1; \hskip 30 pt
\Upsilon_1=R_D(\Upsilon_0)
\end{equation}

That is, we continue our chain by reflecting
it across the line $D$.  The resulting chain
contains $S_s$, by symmetry, but we are not
quite done.

\begin{center}
\resizebox{!}{2.6in}{\includegraphics{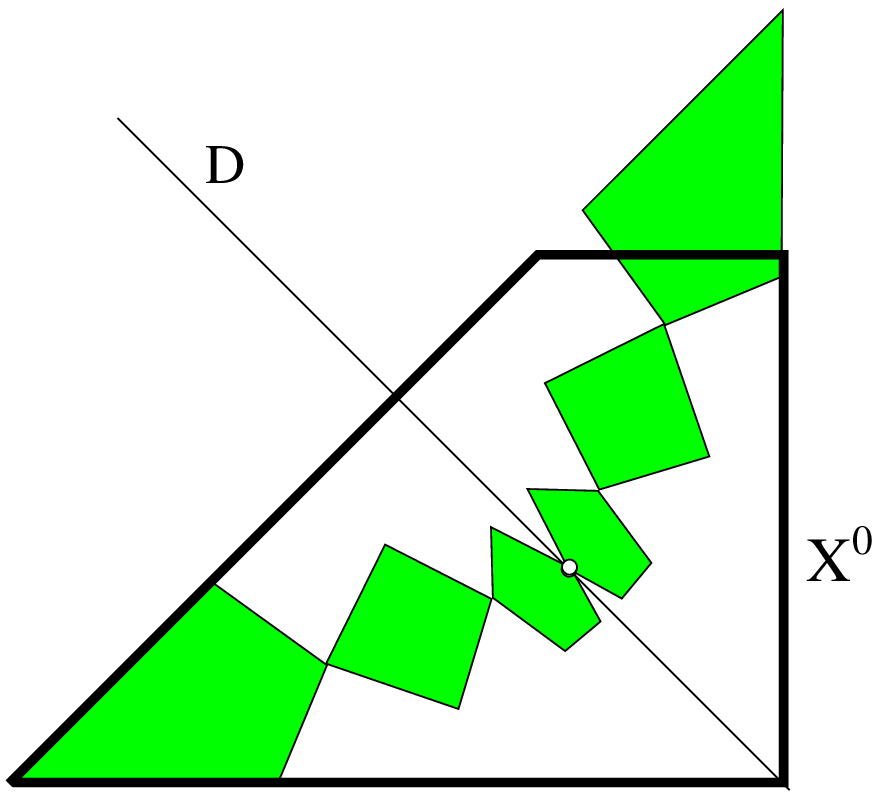}}
\newline
{\bf Figure 5.3:\/} Step 3: Extending by reflection
\end{center}

Some of the final pieces of
$\Upsilon_1$ might not lie in $X_s$.  The problem
is that $X_s$ is not symmetric with respect to
$R_D$. The portion below $D$ is larger than the
portion above $D$.
\newline
\newline
{\bf Step 4:\/}
We finish the construction by a method very similar
to what we did in Step 2.  We simply intersect
all the pieces of $\Upsilon_1$ with $X_s^0$,
and let $\Upsilon_1'$ and omit all those pieces
which come after the first one which has trivial
intersection with $X_s^0$.  We set
$\Upsilon_3=\Upsilon_0,\Upsilon_1'$.
By construction, $\Upsilon_3$ is a chain filled by
$S_s$.   Moreover, since $\phi_s$ contracts
distances by some $\lambda_s<1/\sqrt 2$. we see
that the mesh of $\Upsilon_3$ is less than
$m/\sqrt 2$.
\newline
\newline
{\bf Step 5:\/}
The chain $\Upsilon_3$ might not be clean.  To
remedy this, we shrink the pieces slightly (away
from the marked points) so that they are all
disjoint from the interiors of the edges of
$A_s$.  What allows us to do this is the
Shield Lemma combined with compactness.
The final chain has all the desired properties.
\endproof

Note that if $s$ is oddly even, then so is
$R(s)$.  The chains in Lemma \ref{good3} all
have mesh size less than $2$.  It now follows
from iterating Lemma \ref{good4} that
$S_s$ fills a good chain having mesh size
less than $\epsilon$, for any given
$\epsilon_0$.  Our Arc Criterion how shows
that $S_s$ is an arc.  This completes the proof
of Statement 1 of the Main Theorem.

\section{Proof of Statement 2}

\subsection{The Main Argument}

In this chapter we suppose throughout that
$s \in (0,1)$ is irrational and
$R^n(s)>1/2$ only finitely often.
This means that there is some
$n$ such that $R^n(s)$ is oddly
even.  Let $f(s)$ denote the smallest
$n$ with this property.  Let $S_s$
denote the left half of
$\widehat \Lambda_s$, as in the
previous chapter.  Our goal
is to show that $S_s$ is a finite forest.

\begin{lemma}
Both $S_s \sqcap A_s$
and $S_s \sqcap B_s$ are
finite unions of arcs.
\end{lemma}

\startproof
When $f(s)=0$, the result is true by Statement 1
of the Main Theorem.
Suppose by induction that this lemma is true
for all $s$ such that $f(s)<N$.  If $s$ is
chosen so that $f(s)=N$, then let $t=R(s)$.
Then $f(t)=N-1$.  By induction, both
$S_t \sqcap A_t$ and
$S_t \sqcap B_t$ are finite unions of arcs.
By the Filling Lemma,
both $\widehat \Lambda_s \sqcap A_s$
and $\widehat \Lambda_s \sqcap B_s$ are
contained in finite unions of
similar copies of $\widehat \Lambda_t \sqcap A_t$
and $\widehat \Lambda_t \sqcap B_t$.   
The argument is very similar to what we did
to show that $S_s$ is connected in the 
oddly even case.   Hence, both
both $S_s \sqcap A_s$
and $S_s \sqcap B_s$
are finite unions of arcs.

\begin{corollary}
$S_s$ is a finite union of arcs.
\end{corollary}

\startproof
We have $S_s=(S_s \sqcap A_s) \cup (S_s \sqcap B_s)$.
\endproof

Following this section, the entire chapter is
devoted to proving the
following result.
\begin{theorem}[No Loops]
Let $s \in (0,1)$ be any irrational number.
Then $S_s$ contains no embedded loops.
\end{theorem}
But a finite union of arcs which contains
no embedded loops must be a finite forest.
Modulo the No Loops Theorem, this completes
of Statement 2 of the Main Theorem.

We prove the No Loops Theorem using the same kind of
extremality argument we gave for the Pinching Lemma.
One case of this argument is much harder than the others,
and so we split it off and tackle it first.  Once
we build the machinery for this one case, the
rest of the proof is easy.

\subsection{Interaction with the Octagrid}

Suppose that $s<1/2$ and $t=R(s)<1/2$.  We consider
the octagrid components defined in \S \ref{octagrid}.
Here is the main result of this section.
\begin{lemma}
An embedded loop in $S_s$ lies in a
single octagrid component.
\end{lemma}

We will prove this through a series of smaller results.
Say that an {\it octagrid edge\/} is an edge of one
of the octagrid components.
Let $u=R(t)$.

\begin{lemma}
\label{o2}
$S_s$ intersects each octagrid segment in at most one point
provided that $u<1/2$.
\end{lemma}

\startproof
Let $G$ be an octagrid component and
let $\sigma$ be an edge of $G$.
Just as in the proof of Lemma \ref{octagrid},
we can use bilateral symmetry and the map $T_s$
from Equation \ref{tube} to reduce to the case
where $G \subset Z_s^0$.  If one endpoint of
$\sigma$ is contained in a square of the pyramid
associated to $\Delta_s$, then $G$ is entirely
contained in that square.  This case is trivial.
So, $\sigma$ must be one of the $8$ segments
emanating from the center of the 
bottom left tile of the extended pyramid
associated to $\Delta_s$.

By Theorem \ref{renorm}, it suffices to prove
that $\sigma'=\phi_s^{-1}(\sigma)$ intersects at most
one point of $S_t$, where $t=R(s)$.  One
endpoint of $\sigma'$ is the center $c'$ of the central
tile of $Z_t$.  We think of $\sigma'$ as
pointing away from $c'$. 
Since $u<1/2$, the central
tile of $Z_t$ is a square.  
The $4$ nontrivial cases are as follows.
\begin{enumerate}
\item $\sigma'$ is horizontal and points left.
\item $\sigma'$ has slope $1$ and points down.
\item $\sigma'$ is vertical and points down.
\item $\sigma'$ has slope $1$ and points up.
\end{enumerate}
In the other $4$ cases, $\sigma'$ is entirely
contained in the central tile of $Z_t$.
In the $4$ nontrivial cases, $\sigma'$ is longer than
an octagrid edge associated to the parameter $t$,
but this does not bother us.

The line $D_t$ of symmetry contains the point $c'$.
Reflection in $D_t$ reduces Case 3 to Case 1 and
Case 4 to Case 2.  So, it suffices to consider 
Cases 1 and 2.  
In Case 1,  $\sigma'$ is a segment of $H^0$.  In this case,
we apply the Pinching Lemma.  In Case 2,
we apply Theorem \ref{renorm} again:
$\sigma''=\phi_t^{-1}(\sigma')$ is the left half of
$H$ and intersects $S_u$ at most once by the
Pinching Lemma.  
\endproof

\begin{lemma}
$S_s$ intersects each octagrid segment in at most one point
provided that $u>1/2$.
\end{lemma}

\startproof
The proof here is the same, except that
the central tile of $Z_t$ is an octagon.
In this case, we immediately renormalize
so that $\sigma''=\phi_t^{-1}(\sigma')$ is
one of the $8$ segments emanating from the
center $(0,0)$ of the trivial tile of
$\Delta_u$.  After reflecting in the origin,
we arrive at $2$ nontrivial cases.
\begin{enumerate}
\item $\sigma''$ is horizontal and points left.
\item $\sigma''$ has slope $1$ and points downward.
\end{enumerate}

In the first case $\sigma''$ is simply the portion of $H$
lying to the left of the origin.  This case follows
from the Pinching Lemma applied to the parameter $u$.
In the second case,
$\sigma''$ agrees with the line $E_u$ of symmetry
outside the trivial tile of $\Delta_u$.  Again, this
follows from the Pinching Lemma applied to the
parameter $u$.
\endproof

We record the obvious corollary.

\begin{corollary}
\label{o3}
$S_s$ intersects each octagrid segment in at most one point.
\end{corollary}

Suppose now that $\gamma$ is an embedded loop in $S_s$.

\begin{lemma}
$\gamma$ cannot link any vertex of an octagrid component.
\end{lemma}

\startproof
$\gamma$ cannot surround any square in the pyramid
associated to $s$, because this pyramid is a topological
disk sharing an arc with $\partial X_s$. 

Say that a {\it peripheral tile\/} is a tile of
the extended pyramid which is not a tile of the
pyramid.  The peripheral tiles are colored blue
in Figures 5.8 and 5.9.

Say that a {\it blocker\/} is an octagrid edge which
lies entirely in a tile and has a vertex in
$\partial X_s^0$.  These segments are
disjoint from $S_s$ and hence from
$\gamma$.
Each peripheral tile has at least
one blocker.  But then $\gamma$ cannot surround
a peripheral tile because it would have to
intersect a blocker.

The remaining vertices of the octagrid (meaning,
the vertices of the octagrid components) have the
following structure.  At least one octagrid
edge emanating from the vertex lies entirely inside
a tile.  For this reason, $\gamma$ cannot link
any vertex of an octagrid component.   
\endproof

\begin{center}
\resizebox{!}{2.6in}{\includegraphics{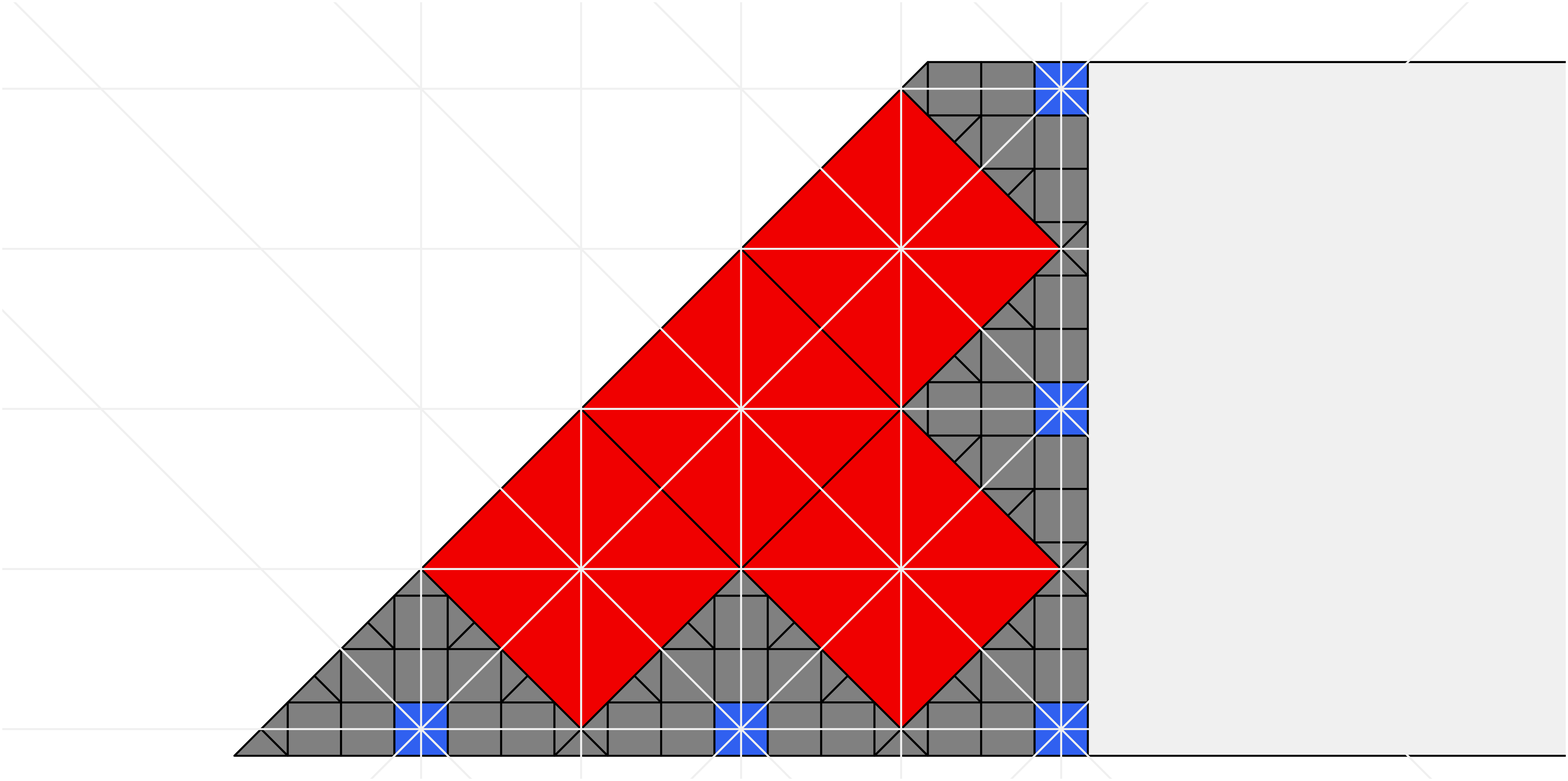}}
\newline
{\bf Figure 6.1\/}: The octagrid for $s=13/32$.
\end{center}  

\begin{center}
\resizebox{!}{2.8in}{\includegraphics{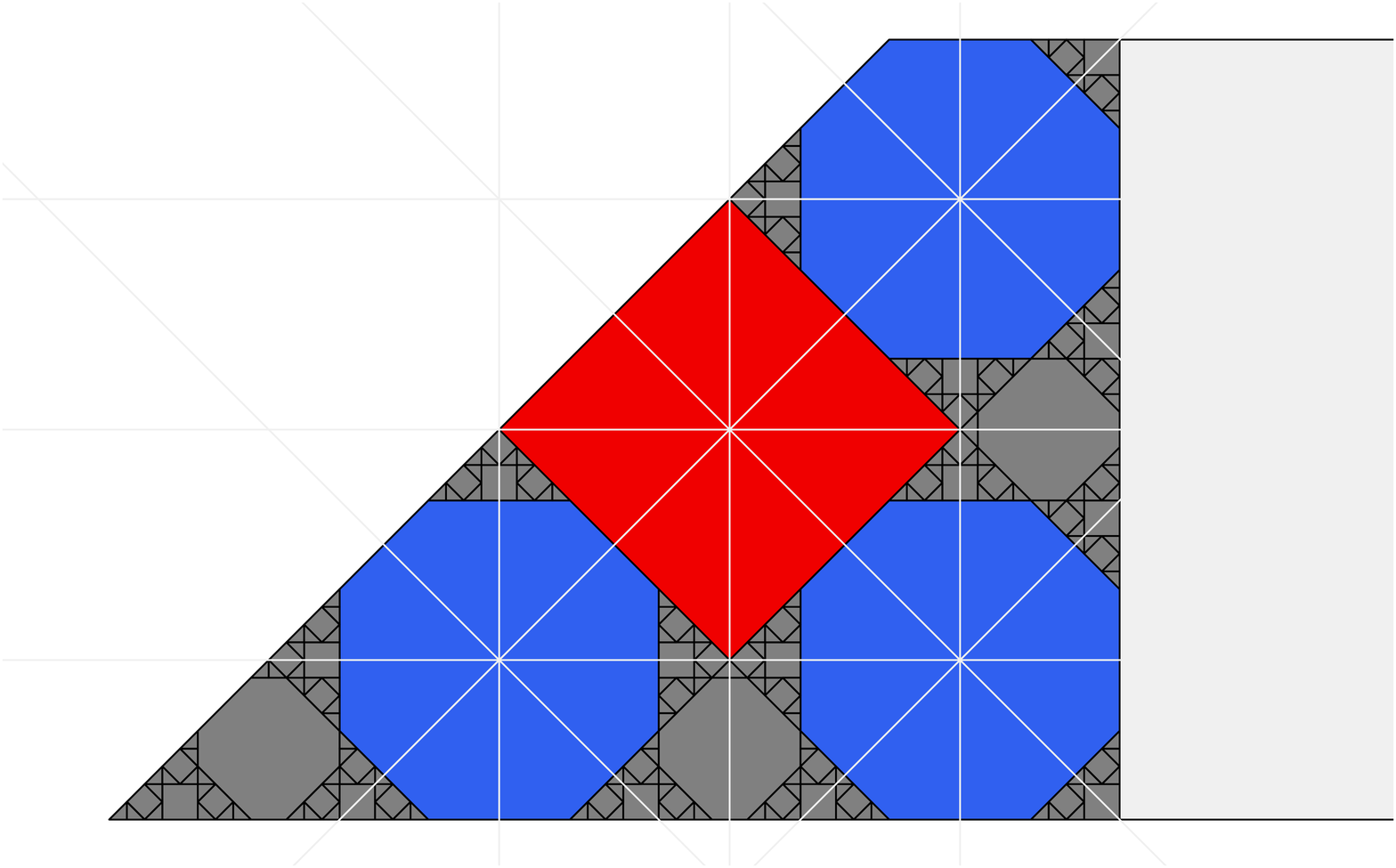}}
\newline
{\bf Figure 6.2\/}: The octagrid for $s=22/57$.
\end{center}  

If $\gamma$ is not contained in a single octagrid
component, and $\gamma$ does not link any of the
vertices of octagrid components, then $\gamma$
must cross some octagrid edge twice.  We have
already ruled this out.  Hence $\gamma$ is
contained in a single octagrid component.
This completes the proof of Lemma \ref{o2}.

\subsection{No Embedded Loops}
\label{noloop}

Now we prove the No Loops Theorem.
Our method of proof is very much like what we did for
the Pinching Lemma. 
Say that a {\it counterexample\/}
is a pair $\Omega=(\gamma,s)$ where $\gamma$ is an embedded
loop in $S_s$.  We define
\begin{equation}
\lambda(\Omega)=\frac{{\rm diam\/}(\gamma)}{g(s)}, \hskip 30 pt
g(s)={\rm diam\/}(X_s^0).
\end{equation}

\begin{lemma}
\label{lazy22}
For any $\epsilon>0$ there exists a
counterexample $\Omega$, whose parameter
is less than $1/2$, such that $\lambda(\Omega)>\lambda-\epsilon$.
\end{lemma}

\startproof
As in the proof of Lemma \ref{lazy2}, it suffices
to consider the case of a counterexample
$\Omega$ having $s \in (1/2,\sqrt 2/2]$.
In this case, the reflection
$R_V$ maps the region to the
right of $V$ into $Z_s^0$ and
the reflection $R_D$ maps the region above
$D_s$ into $Z_s^0$.  Figure 6.3 shows a
fairly typical example.

\begin{center}
\resizebox{!}{2in}{\includegraphics{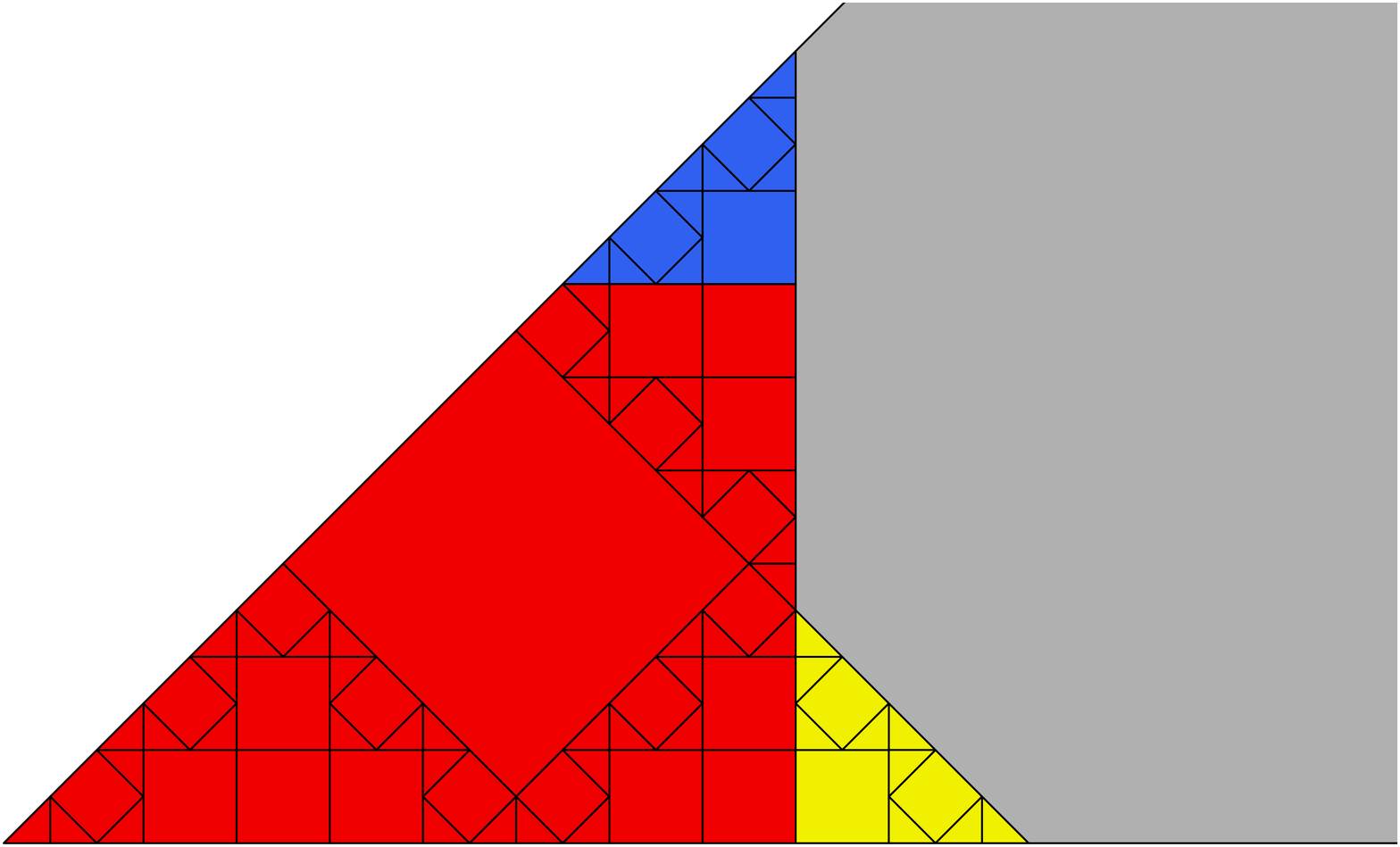}}
\newline
{\bf Figure 6.3:\/} 
$Z_s^0$ (red) and the components of $X_s^0-Z_s$ for $s=11/17$.
\end{center}

By the Pinching Lemma, $\gamma$ can intersect
each of $V$ and $D_s$ at most once.
Hence $\gamma$ lies to one side or the other
of each of these lines.  So, by symmetry,
we can assume that $\gamma \subset Z_s^0$.
The rest of the proof is as in Lemma \ref{lazy2}.
\endproof

As in the proof of the Pinching Lemma, we
now analyze potential counterexamples
having parameter $s \in (1/4,1/2)$.
First suppose $s \in (1/4,1/3)$. 
For $s$ in this range, we have $t=R(s)>1/2$.
Let $O_s$ be the image, under the $\phi_s^0$
of the central tile of $\Delta_s$.  This
octagon separates $X_s^0$ into three regions
having disjoint closures, as shown in Figure 6.3.
One of the regions, colored red in Figure 6.3,
is precisely $Z_s^0$.

\begin{center}
\resizebox{!}{1.6in}{\includegraphics{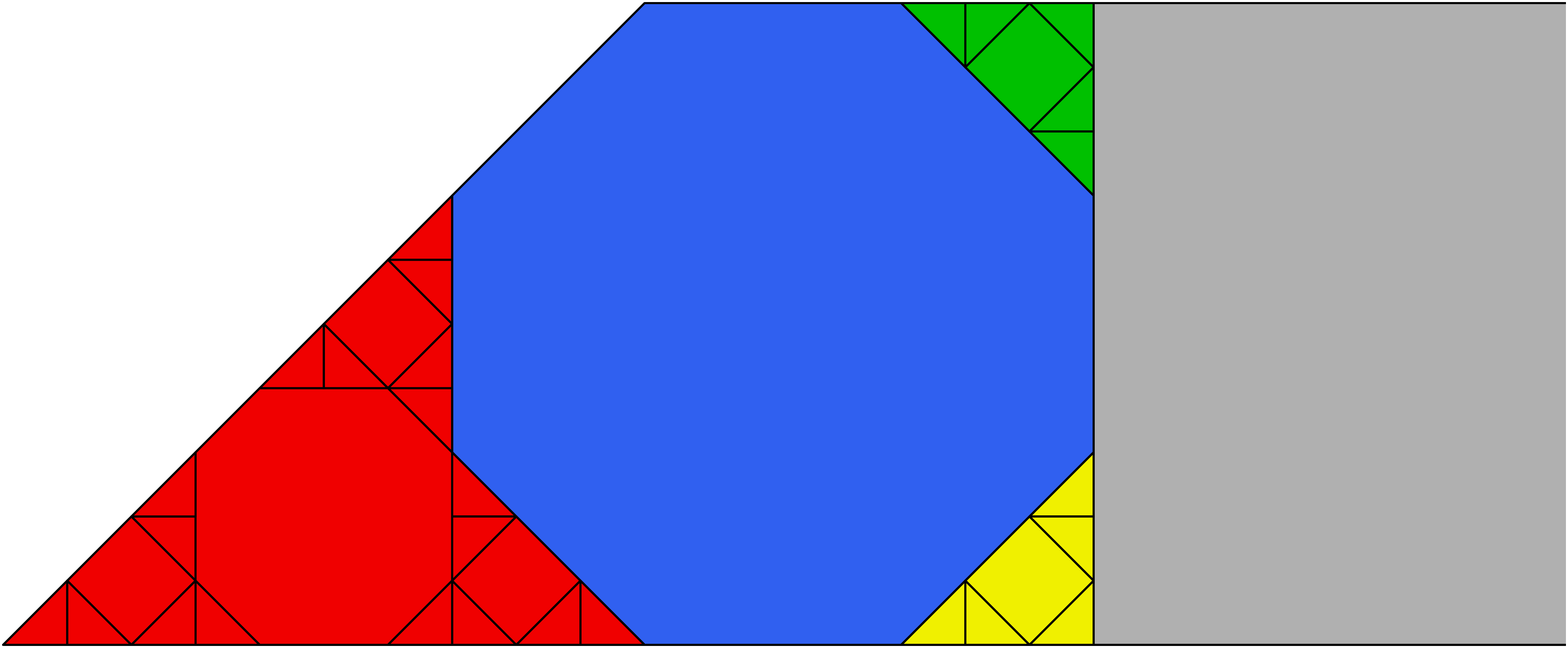}}
\newline
{\bf Figure 6.4:\/} 
$Z_s^0$ (red) and other regions for $s=5/17$.
\end{center}

The reflections $R_D$ and $R_V$
respectively carry the other regions into the
one contined in $Z_s$.  Hence, we may assume
by symmetry that $\gamma \subset Z_s$.
We now get the same contradiction as we
got in the proof of the Pinching Lemma.

For $s \in (1/3,1/2)$ we have
$R(s)<1/2$.  Lemma \ref{o2} applies,
and so we know that our counterexample
is contained in a single octagrid
component.  But then, by Lemma \ref{octagrid},
we can assume by symmetry that our counterexample
lies in $Z_s$.  This gives us the same
contradiction as in the Pinching Lemma.

This completes the proof that $\widehat \Lambda_s$
has no embedded loops.

\newpage

\section{Proof of Statement 3}

\subsection{The Limit Set is Perfect}

Recall that a closed set $C$ is {\it perfect\/} if
every point $p \in C$ is an accumulation point
of $C-\{p\}$.  

\begin{lemma}
\label{perfect}
When $s$ is irrational,
$\widehat \Lambda_s$ is a perfect set.
\end{lemma}

\startproof
If $\widehat \Lambda_s$ is not perfect,
then there is some $p \in \widehat \Lambda_s$
and come open disk $U$ containing $p$ such that
$U \cap \widehat \Lambda_s=p$.
The open set $U$ must contain infinitely
many tiles of $\Delta_s$, 
because $\widehat \Lambda_s \cap U$ is nonempty.
Therefore, if we write $\Delta_s$ as in the Covering Lemma,
the image of some patch must have $p$ as an accumulation
point.  Choosing $\epsilon$ small enough, we can guarantee
that there exists an $\epsilon$-patch
$(K,\psi,\epsilon)$ such that $\psi(K) \subset U$.

If $K$ is a 
triangle, then two of the vertices $v_1$ and $v_2$ of $K$
have acute angles. (The angle is $\pi/4$.)
These vertices must be
accumulation points of infinitely many tiles,
because all the tiles are squares and
semi-regular octagons.
But then $\psi(v_1)$ and $\psi(v_2)$ are accumulation
points of infinitely many squares in
$\Delta_s$.  Hence $\widehat \Lambda_s \cap U$
contains at least $2$ points. This
is a contradiction.

If $K$ is a pentagon, then $K$ has $2$ vertices
$v_1$ and $v_2$ with obtuse angles.
(The angle is $3 \pi/4$.)  If $v_1$ is
not an accumulation point of infinitely
many tiles of $\Delta_u \cap K$, then $v_1$
is the vertex of some octagon of $\Delta_u \cap K$,
 But then there are two
new acute vertices $w_1$ and $w_2$ which must
be accumulation points of infintely many tiles
of $\Delta_u \cap K$.  This gives us the
same contradiction as above.  The only way out
of the contradiction is for both $v_1$ and
$v_2$ to be accumulation points of infinitely
many tiles of $\Delta_u \cap K$, but this
is again a contradiction. There is no way out.
\endproof

\subsection{Overview for the Rest of the Proof}

We call $s \in (0,1)$ an
{\it octagonal parameter\/} if 
$R^n(s)>1/2$ infinitely often.
The reason for the name is that,
thanks to Theorem \ref{tiling},
$\Delta_s$ contains infinitely many
octagons if and only if $s$ is
octagonal.
Our remaining goal is to prove
that $\widehat \Lambda_s$ is a Cantor
set when $s$ is octagonal.  Let $S_s$
denote the left half of $\widehat \Lambda_s$,
as usual.

We already know that $S_s$ is closed and perfect.
It remains to show that $S_s$ is totally
disconnected 
for any octagonal $s$.  The proof 
is a bootstrap argument.
The basic idea is that renormalization
tends to make connected components of
$S_s$ larger.  So, we will
start with the assumption that $S_s$
has a nontrivial connected component when
$s$ is octagonal, and ultimately we will produce
a new octagonal parameter $u$ for which
$S_u$ has a very large 
and egregious kind of connected component,
which we can rule out.  Throughout the proof,
$K_s$ will stand for a {\it symmetric piece\/},
one of the sets $A_s,B_s,P_s,Q_s$ from
\S \ref{symm}.

\subsection{Unlikely Sets}

\begin{lemma}
\label{barrier}
Let $s \in (0,1)$ be irrational. 
Let $K_s$ be a symmetric piece.
Let $e$ be any edge of $K_s$. Then
some tile of $\Delta_s \cap K_s$
has an edge in $e$.
\end{lemma}

\startproof
We do this in a case-by-case way.
\begin{enumerate}
\item Suppose $s<1/2$ and $K_s \in \{B_s,Q_s\}$.
Applying Theorem \ref{renorm} repeatedly, we
see that there are infinitely many tiles which
have edges in the bottom edge of $X_s$ and
infinitely many tiles which have edges in
the left edge of $X_s$. Eventually these tiles
lie in both $B_s$ and $Q_s$. This takes care
of two out of three edges of $B_s$ and $Q_s$.
The third edge, in each case, follows from
bilateral symmetry.  

\item Suppose $s<1/2$ and $R(s)>1/2$ and
$K_s \in \{A_s,P_s\}$.  In this case,
there is an octagon having the
same width as the central tiles, and this octagon
has edges in all the sides of $A_s$.  Likewise,
the same octagon has edges in all the sides of
$P_s$.  See Figure 4.12.

\item Suppose $s<1/2$ and $R(s)<1/2$ and
$K_s \in \{A_s,P_s\}$.
In this case, each edge of $A_s$ has a segment
which is an edge of one of the tiles in the
extended pyramid.  The same goes for
$P_s$.  See Figures 3.8 and 3.9.

\item Suppose $s>1/2$.
The argument given in Case 1 takes care of
$P_s$ and $B_s$.
The case of $A_s$ and $Q_s$ follows
from inversion symmetry.  See Remark (i)
in \S \ref{symm}.
\end{enumerate}
Thus exhausts the possibilities.
\endproof

Let $C$ be a nontrivial connected subset of
$\widehat \Lambda_s$.  Let
$\cal K$ be a patch cover, as in
the Covering Lemma. What we mean
is that $\cal K$ is a finite union of
patches and tiles, as in
Equation \ref{patcov}.  We call $C$
{\it bad\/} with respect to 
$\cal K$ if $C$ does not intersect
the interiors of the images of
any of the patches.  That is,
$C$ is disjoint from the interiors 
of all the sets $\psi_j(K_j)$.

Call $C$ {\it unlikely\/} if,
for every $\epsilon>0$, there is a
patch covering $\cal K$ of scale less
than $\epsilon$, so that $C$ is bad
with respect to $\cal K$. So, $C$ is
bad with respect to an infinite sequence
of patch covers, having scale tending
to $0$.

\begin{lemma}
\label{unlikely}
Let $s \in (0,1)$ be irrational.
There do not exist any unlikely subsets of
$\widehat \Lambda_s$.
\end{lemma}

\startproof
We will suppose some unlikely set $C$ exists and get a contradiction.
By taking a suitable
subset of $C$, we can assume that $C$ is a line
segment contained in
the boundary of some particular patch of one
of the sequence of bad covers.  
By Lemma \ref{barrier}, some segment of $C$ lies
in a tile boundary.  Further shrinking $C$, we
can assume that $C$ is one edge of some
tile $\tau_1$.

The midpoint $m$ of $C$ is the accumulation point of
infinitely many tiles of $\Delta_s$. These
tiles are all disjoint from $\tau_s$. Hence,
there is some patch $(\psi,K,u)$ so that
$m \in \psi(K)$.  But $C$ cannot intersect
the interior of $\psi(K)$.  Hence, one
edge of $\psi(K)$ lies in the line containing $C$.
Shrinking the scale as needed, we can assume that
one edge of $\psi(K)$ is contained in $C$.

Some tile $\tau_2$ of $\psi(K) \cap \Delta_s$
has an edge in $C$, by Lemma \ref{barrier}.
But then some point of $C$ lie on the
common boundary of $\tau_1$ and $\tau_2$ and
cannot belong to the limit set.
This is a contradiction.
\endproof

\subsection{Tails}

For each symmetric set $K_u \in \{A_u,B_u,P_u,Q_u\}$, and each
subset $C \subset S_u$, let
$C \sqcap K_u$ denote
the set of points $p \in C$ such that
every neighborhood of $p$ contains infinitely
many tiles of $K_u \cap \Delta_u$.  Note that
$C \sqcap K_u$ might be
a proper subset of $C \cap K_u$, but
the two sets agree on the interior of $K_u$.
The former set is easier to pull back using
the patches from the Covering Lemma.

We say that the symmetric piece $K_u$ has a
{\it tail\/} if some
connected component of $S_u \sqcap K_u$
contains both an interior point of
$K_u$ and a boundary point of $K_u$.
Conceptually, we think of a tail as a little arc which
joins a boundary point of $K_u$ to
an interior point, but we don't actually
know that our connected set is path connected.

\begin{lemma}
\label{cross1}
Suppose $S_s$ is not
totally disconnected. Then for all
sufficiently large $n$ the parameter
$u=R^n(s)$ has the following property.
At least one of the $4$ symmetric
pieces $K_u$ has a tail.
\end{lemma}

\startproof
If $S_s$ is not totally disconnected,
then it has a closed connected subset $C$ having
more than one point.  By Lemma \ref{unlikely},
once $n$ is large enough and $u=R^n(s)$,
we can find a patch $(K,\psi,u)$ such that
some point $p \in C$ lies in the interior of
$K'=\psi(K)$ and some point of $C$ lies outside
of $K'$.

Consider the set $C \sqcap K'$.  For the sake
of exposition, assume first that $C$ is
path connected.  Then we can find some path
$\alpha \in C$ joining $p$ to some point
of $C$ lying outside $K'$.   Let
$\beta$ be the maximal initial portion of
$\alpha$ which lies in $K'$.
By construction
$\beta \subset C \sqcap K'$ and $\beta$ joins
$p$ to a point $q \in \partial K'$.
By definition
of a patch,
$\psi^{-1}(\beta)$ is a path in 
$S_u \sqcap K_u$ joining 
$\psi^{-1}(p)$ to $\psi^{-1}(q)$.
The former point lies in the interior
of $K$ and the latter point lies on the
boundary.  This gives $K_u$ a tail.

We follow the same outline when $C$ is merely connected.
Define an $\epsilon$-{\it chain\/} to be
a sequence of points $p_0,...,p_N \in C$
such that $\|p_k-p_{k+1}\|<\epsilon$
for all $k$.  We say that this chain
joins $p_0$ to $p_N$.
Any open neighborhood of a connected
set is path connected.  Hence,
for any $m$, there is a
$(1/m)$-chain $\alpha_m$ joining $p$ to some
point of $C$ outside of $K'$.
Let $\beta_m$ be the maximal initial
portion of $\alpha_m$ which lies entirely
in $K'$.  Passing to a subsequence, we
can assume that the endpoints of
$\beta_n$ converge.  One of the endpoints
is always $p$.  Let $q$ be the limit of
the other endpoints. By construction
$q \in \partial K'$.

Let $\beta$ be the connected component
of $C \sqcap K'$ containing $p$.
By construction $q \in C \sqcap K'$.
Moreover, for every $\epsilon>0$ there is
an $\epsilon$-chain connecting $p$ to $q$.
Hence $q \in \beta$. Now
we pull back by $\psi$ to give $K_u$ a tail,
as in the path connected case.
\endproof

\noindent
{\bf The Chain Trick:\/}
We will have many occasions below to make
the same kind of argument as we just gave,
and the geometry behind the argument is always
clearer in the path connected case.  For
this reason, we will give the arguments
in the path connected case and then remark
that the same trick as used above -- i.e.,
using an  infinite sequence
of chains in place of a path -- handles the general
case.  For reference, we call this the
{\it Chain trick:\/}

\subsection{Crosscuts}

Let $K$ be a (solid) triangle and
let $\Lambda \subset K$ be a closed set.
We say that a {\it crosscut\/} for
$(K,\Lambda)$ is a connected subset of
$\Lambda$ which contains two points
$p,q \in \partial K$.  We require that
$p$ and $q$ are not both contained in
the interior of the same edge of $K$.

In the following result, it is not
really essential that $s$ is octagonal
and that $s,u>1/2$.  However, in this case,
the symmetric 
pieces associated to the two parameters
are all right-angled isosceles triangles.
This makes the geometry easier and
cuts down on the number of cases to
consider.

\begin{lemma}
Suppose that $s>1/2$ is octagonal and
$S_s$ is not
totally disconnected. Then for all
sufficiently large $n$ the parameter
$u=R^n(s)$ has the following property
provided that $u>1/2$:
The pair
$(K_u,K_u \sqcap S_u)$
has a crosscut for some symmetric piece $K_u$.
\end{lemma}

\startproof
By Lemma \ref{cross1}, we can assume
without loss of generality that some
$K_s$ has a tail $C$.
As remarked above, we will assume that
$C$ is a path; the Chain Trick then
finishes the proof.

Let $\epsilon=\|p-q\|$ and choose $n$
so large that the patch covering
in the Covering Lemma corresponding to
$u=R^n(s)$ has scale much smaller than
$\epsilon$.  There is a patch
$(K_j,\psi_j,u)$ such that
$K'=\psi_j(K_j)$ contains $p$.
There are several cases to consider,
as shown in Figure 7.1.  (The third part
of Figure 7.1 just shows one of the several
possibilities for that case.)\begin{center}

\resizebox{!}{.8in}{\includegraphics{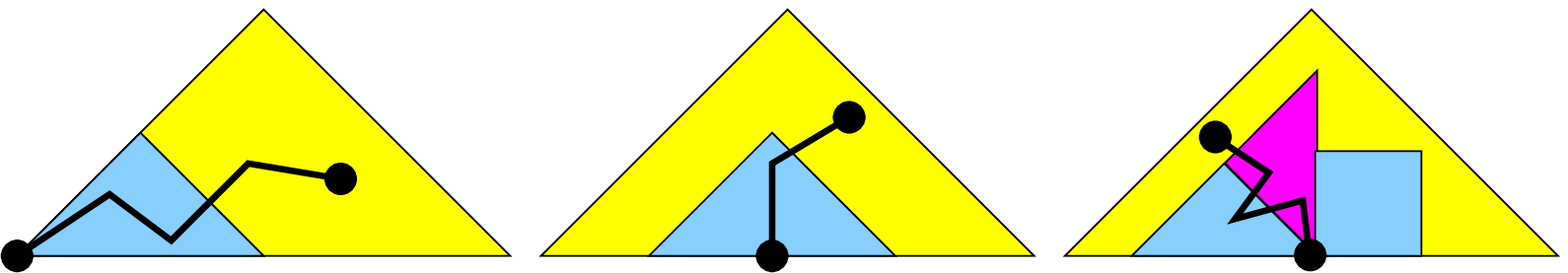}}
\newline
{\bf Figure 7.1:\/} Possibilities for the patch covering.
\end{center}

\begin{itemize}
\item Suppose that $p$ is a vertex of $K_s$.
Then $p$ is also a vertex of $K'=\psi_j(K_j)$.
The two sides of $K'$ incident to $p$ lie
in sides of $K$.  But then the path $C$
must exit through the third side of $K'$
in order to reach $q$.  The pullback
$\psi_j^{-1}(C)$ is a crosscut of $K_j$.
\item Suppose $p$ lies in the interior of an edge of
$K$ and also in the interior of an edge of $K'$.
Then, again, $C$ must
exit $K'$ through one of the other edges of $K'$.
Now we pull back as before.
\item Suppose $p$ lies in the interior of an edge
of $K$ and is a vertex of $K'$.  In this
case, we can assume that $p$ is not contained
in the interior of an edge of the image
any other patch.  So, a neighborhood of $p$ in
$K$ is covered by finitely many tiles of
$\Delta_s$ and either $1$, $2$, or $4$
patch images, each of which has $p$ as a vertex.
In any case, $\gamma$ must exit this neighborhood,
and some initial portion of $\gamma$ makes a
crosscut in one of the patch images.  Now we
pull back as before.
\end{itemize}
This exhausts the possibilities.
\endproof

The rest of the proof involves ruling out
the existence of crosscuts associated to octagonal
parameters.

\subsection{The Proof Modulo one Case}

\begin{lemma}
Suppose that $s>1/2$ is an octagonal parameter
and $K_s$ is a symmetric piece with a crosscut.
Then the crosscut cannot contain an acute vertex
of $K_s$.
\end{lemma}

\startproof
Suppose the crosscut contains a vertex $p$.
Using the various kinds of bilateral symmetry discussed
in \S \ref{symm} we reduce to the case when 
$p$ is the bottom left vertex of $X_s$ and
$K_s$ is one of $B_s$ or $P_s$.  (These are the
two pieces having this point as a vertex.)

By Corollary \ref{wedge2}, there are infinitely many
octagons wedged into $X_s$. These octagons shrink
down to $p$, and isolate $p$ from the rest of
$S_s$.  So, $p$ is its own
connected component of $S_s$.
\endproof

We call a crosscut of $K_s$ {\it acute\/} if it
contains points $p,q$, where $p$ lies in the
interior of a short side of $K_s$ and
$q$ lies in the interior of a long side of $K_s$.
We use this name, because we think of a path
joining $p$ to $q$ and subtending one of the
acute angles.

\begin{lemma}[Acute]
\label{acute}
Suppose that $s>1/2$ is an octagonal parameter
and $K_s$ is a symmetric piece with a crosscut.
Then the crosscut cannot be acute.
\end{lemma}

\startproof
We prove this in the next section.
\endproof

Say that a
{\it right crosscut\/} is a crosscut of
$K_s$ which contains the right-angled
vertex of $K_s$.

\begin{lemma}
Suppose that $s>1/2$ is an octagonal parameter
and $K_s$ is a symmetric piece with a crosscut.
The crosscut cannot contain the right-angled
vertex of $K_s$.
\end{lemma}

\startproof
Let $C$ be a crosscut which supposedly
contains the right-angled vertex $p$ of $K_s$.
We will give the proof when $C$ is a path,
and then the Chain Trick handles the general case.
Since $s$ is octagonal, there are
arbitrarily large choices of $n$ for
which $u=R^n(s)>1/2$,  If we choose
$n$ sufficiently large, then by
the Covering Lemma we can find a patch
$(\psi_j,K_j,u)$ so that
$K'=\psi_j(K_j)$ has $p$ as a vertex,
shares two edges with $K$, and has
diameter (say) $\epsilon/10$.
This is shown in Figure 7.2.

\begin{center}
\resizebox{!}{1in}{\includegraphics{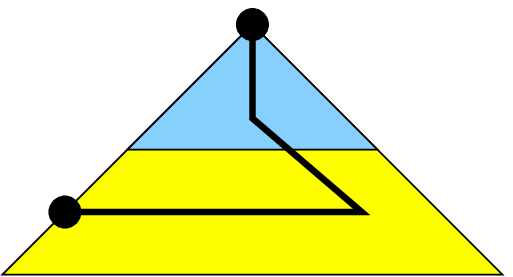}}
\newline
{\bf Figure 7.2:\/} The pieces $K$ and $K'$.
\end{center}

But then $C$ must exit $K'$ from the hypotenuse.
The maximal initial portion of $C$ contained in
$K'$ is an acute crosscut (with respect to either
acute angle). But, Lemma \ref{acute} says that
these cannot exist.
\endproof

Before we give the final argument, we
single out some important points.
\newline
\newline
{\bf The Proper Nesting Property:\/}
Suppose $K'=\psi_j(K_j) \subset K_s$ is the image
of some patch which arises in the conclusion
of the Covering Lemma for the parameters $s$
and $u$.  If $K'$ and $K_s$ share a vertex,
then this vertex has the same type (acute
or right) with respect to both triangles.
One sees this just by inspecting the equations
used in the proof of Lemma \ref{cov0}.
Call this the {\it proper nesting property\/}.
\newline
\newline
{\bf Images of Crosscuts:\/}
If $(K_j,\psi_j,u)$ is a patch
for the parameter $s$, and $K_j$ has
a right crosscut, then some connected subset
$S_s \sqcap K'$ contains
points on the interiors of both the short edges
of $K'=\psi_j(K_j)$. We will abuse our
terminology and say that $K'$ has a
crosscut, even though technically $K'$ is a
similar copy of a symmetric piece associated
to the different parameter $u$.  With this
terminology, a crosscut for $K'$ is a subset 
of $S_s$.
\newline

We call a crosscut of $K_s$ {\it right\/} if it
contains points $p,q$, where $p$ and $q$ 
respectively lie in the interiors of the
two short sides of $K_s$.   Modulo proving
Lemma \ref{acute}, we have ruled out all
types of crosscut except right crosscuts.
So, if Statement 3 of the Main Theorem is
false, then there is an infinite sequence
of octagonal parameters $s_1,s_2,...$ with
the following properties.
\begin{itemize}
\item $s_{k+1}=R^{n_k}(s_k)$ for some $n_k>0$.
\item Some symmetric piece $K_{s_k}$ has
a right crosscut for $k=1,2,3...$
\item The crosscut associated to $K_{s_k}$
is disjoint from a neighborhood of the
right-angled vertex of $K_{s_k}$.  (Otherwise,
we reduce to the previous case.)
\end{itemize}

Once one symmetric piece $K_s$ has a crosscut,
all similar copies of $K_s$ have crosscuts.
Moreover, by Lemma \ref{cov0},
each similar copy of a
symmetric pieces for the parameter
$s_j$ contains a similar copy of a
symmetric piece for the parameter $s_k$
as long as $k>j$.

Using the proper nesting property together
with the Pidgeonhole Principle, we can fnd two nested
(similar copies of)
symmetric pieces which both have right crosscuts.
We can arrange that the smaller piece is so small
that its crosscut is completely disjoint from
the crosscut of the larger piece, as shown in Figure 7.3.

\begin{center}
\resizebox{!}{1in}{\includegraphics{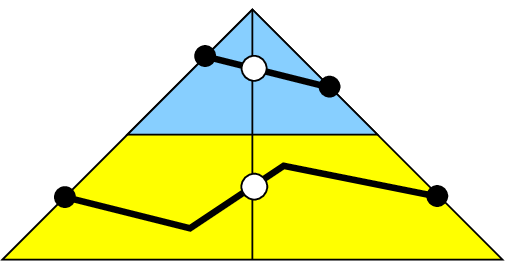}}
\newline
{\bf Figure 7.3:\/} Nested symmetric pieces
\end{center}

Call $s$ the parameter associated to the larger of
the two symmetric pieces.  By construction, both
the crosscuts we have found must cross the line
of symmetry of $K_s$.  Moreover, both crosscuts
are subsets of $S_s$.  Hence,
$S_s$ crosses the line of
symmetry for $K_s$ at least twice.  This
contradicts the Pinching Lemma.

The only way out of the contradiction is that
Statement 3 of the Main Theorem is true.
This completes the proof of the Main Theorem,
modulo the proof of Lemma \ref{acute}.

\subsection{No Acute Crosscuts}

Here we complete the proof of the Main
Theorem by establishing Lemma \ref{acute}.
If some symmetric piece $K_s$ has an
acute crosscut, then we can use bilateral
symmetry to reduce to the case when
$S_s$ has a connected
subset $C$ containing a point on the
bottom edge of $X_s$ and a point on the
left edge of $X_s$.  We abbreviate this
by saying that the parameter 
$s$ has an acute crosscut.
At this point we no longer care about
the condition that $s>1/2$.  As above,
we will treat the case when $C$ is a path;
the Chain Trick handles the general case.

\begin{lemma}
If an octagonal parameter $s$ has an acute
crosscut, then there is another octagonal
parameter $u$ which has an acute crosscut
not entirely contained in $Z_s$.
\end{lemma}

\startproof
If $s$ has an acute crosscut contained
in $Z_s$ and $t=R(s)$ also has an acute
crosscut.  If the acute crosscut for
$t$ lies in $Z_t$ we can repeat the
procedure.  Every one or two steps of
the procedure, the distance between the
endpoints of the crosscut increases by a
factor of at least $\sqrt 2$.
So, eventually we reach a stage where
the crosscut cannot lie in $Z_s$.
\endproof

Say that an {\it egregious crosscut\/} is 
a connected subset of $S_s$
which contains a point in the top edge of
$X_s$ and a point on the bottom edge of
$X_s$. 

\begin{lemma}
\label{egg1}
If $s<1/2$ has an acute crosscut which does
not lie in $Z_s$ then $R(s)$ has an
egregious crosscut.
\end{lemma}

\startproof
Define
\begin{equation}
\label{diag}
Z_s^* = (Z_s^0 \cup R_D(Z_s^0)) \cap X_s
\end{equation}
Figure 7.4 shows this set.

\begin{center}
\resizebox{!}{1.6in}{\includegraphics{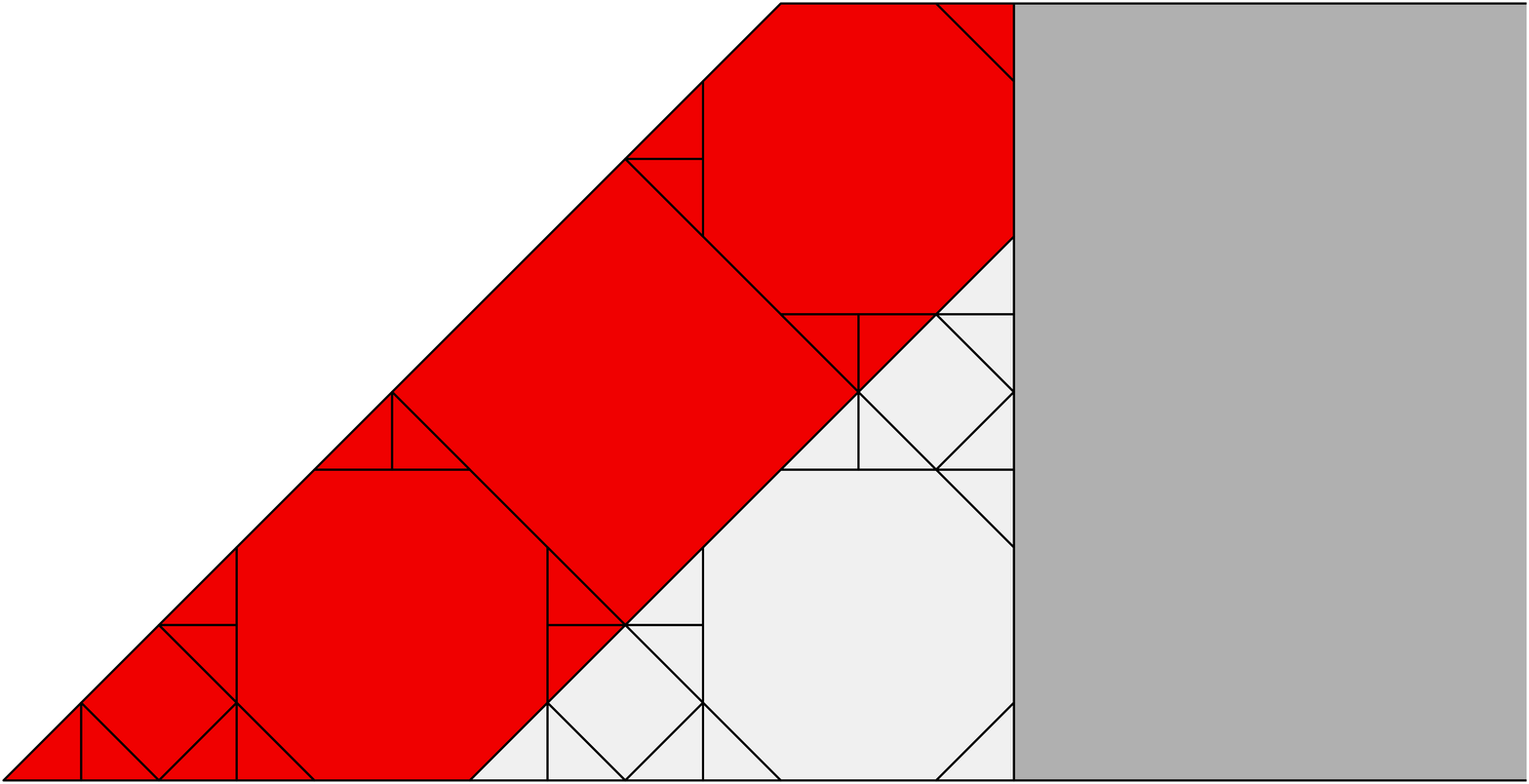}}
\newline
{\bf Figure 7.4:\/} The set $Z^*_s$ for $s=5/13$.
\end{center}

$C$ starts out on the left edge of $Z_s^*$ and
must eventially reach the bottom edge of $X_s$.
But then $C$ must cross the right edge of 
$Z_s^*$.  By symmetry, the left branch $\phi_s^0$ of the
map $\phi_s$ from Theorem \ref{renorm} extends
to all of of $Z_s^*$.   When we pull back the
maximial initial portion of $C$ lying in $Z_s^*$,
we get an egregious crosscut for $t$.
\endproof

\begin{lemma}
\label{egg2}
If $s>1/2$ has an acute crosscut which does
not lie in $Z_s$, then one of
$R(s)$ or $R^2(s)$ has an egregious crosscut.
\end{lemma}

\startproof 
Let $K$ be the layering constant for $s$.
When $K>1$ we define 

\begin{equation}
Z_s^* = (Z_s^0 \cup R_V(Z_s^0)) \cap X_s
\end{equation}

\begin{center}
\resizebox{!}{1.5in}{\includegraphics{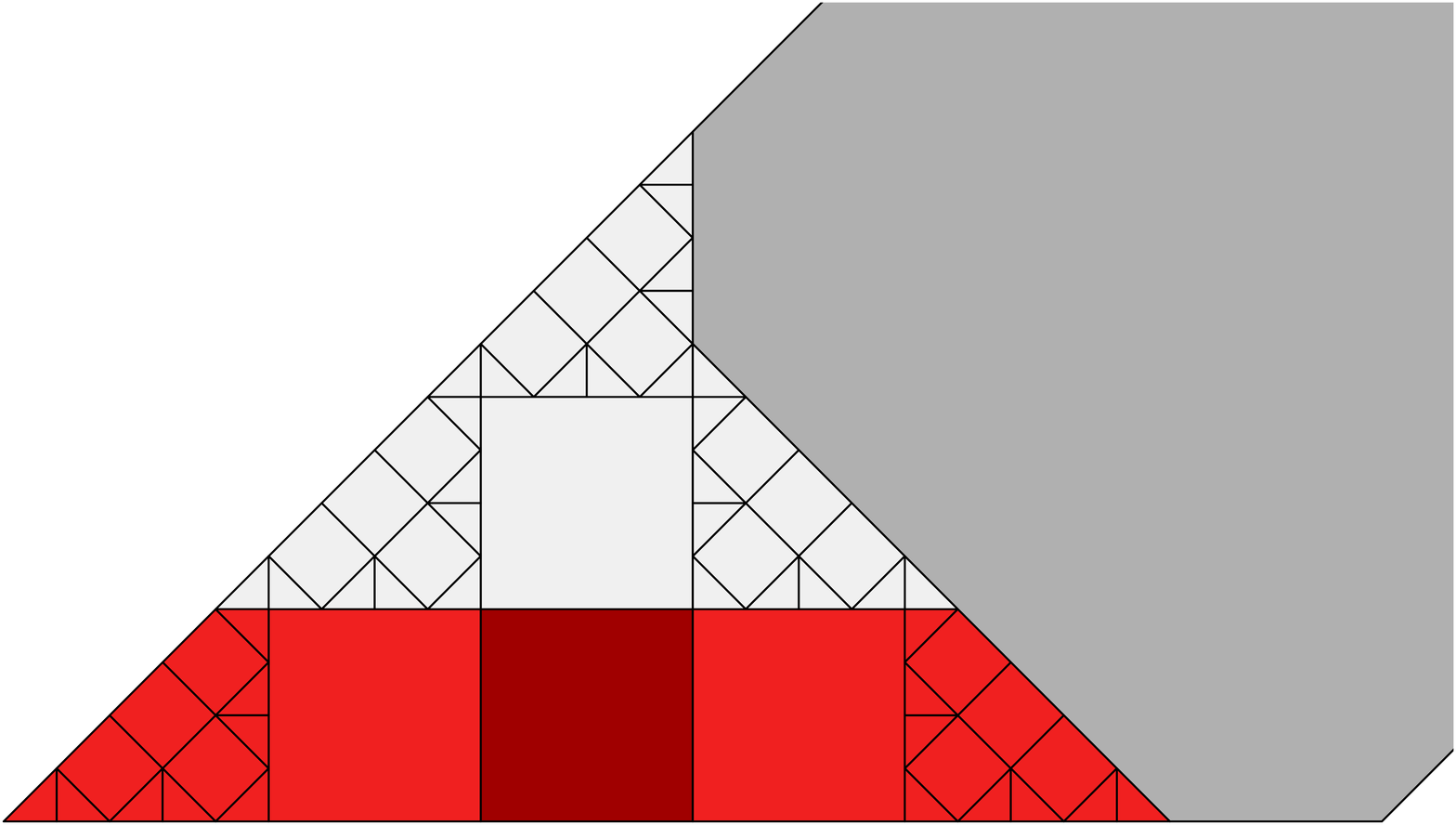}}
\newline
{\bf Figure 7.5:\/} The set $Z^*_s$ for $s=11/13$.
\end{center}

If $C$ starts out above $Z_s^*$ then $C$ must
cross both the top and bottom edges of $Z_s^*$ to
reach the bottom edge of $X_s$.  If $C$ starts
out in $Z_s^*$ but does not lie entirely in
$Z_s^0$, then $C$ must exit the top of $Z_s^*$.
The problem is that $C$ cannot penetrate through
the *dark red) central tile of $Z_s^0$.  So, in this case
as well, $C$ crosses both the top and the bottom
of $Z_s^*$.  Moreover, $C$ must cross both
sides on the same side (left or right) of the
central tile of $Z_s^0$. Reflecting in $V$,
we can assume that $C$ crosses both sides
of $Z_s^*$ to the left of the central tile
of $Z_s^0$.  But then some connected subset
of $C$ lies in $Z_s^0$ and contains points
both on the top edge and the bottom edge.
Now we pull back by $\phi_s$, as above.

\begin{center}
\resizebox{!}{1.5in}{\includegraphics{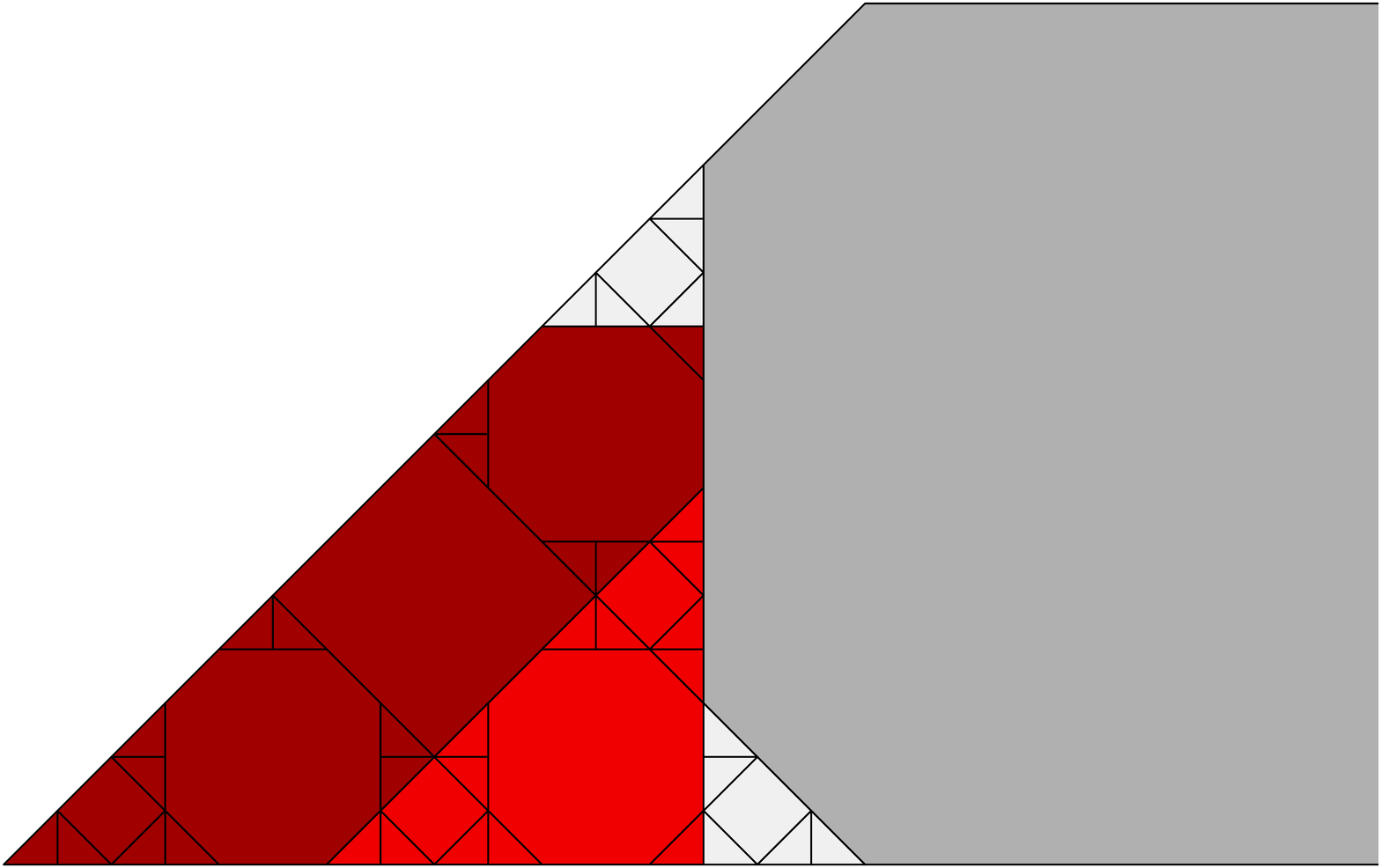}}
\newline
{\bf Figure 7.6:\/} The set $Z_s$ for $s=11/13$.
\end{center}

Now suppose that $K=1$.  Figure 7.6 shows a
representative example. The red set is
$Z_s^0$. The dark red subset is a similar
copy of $Z_t^*$, defined in Equation \ref{diag}.
Call this similar copy $\Omega$.
If $C$ exits $Z_s^0$ then one of two things
must happen.
\begin{enumerate}
\item $C$ contains points on the left and right
diagonal edges of $\Omega$.
\item $C$ contains points on the top edge and
right diagonal edge of $\Omega$.
\end{enumerate}

In the first case, we proceed
as in Lemma \ref{egg1}, and get an egregious
crosscut for $u=R(t)=R^2(s)$.
In the second case, we also
proceed as in Lemma \ref{egg1}, but we get
something different:  After reflecting through
the origin, we get a connected
set $C' \subset S_u$ which has points on the bottom
edge of $X_u^0$ and some point lying above
the line $H$ of symmetry.  (The point is that
$C$ crosses the diagonal midline of $\Omega$.)
But then $C''=C' \cup R_H(C')$ is an
egregious crosscut for $u$.
\endproof

Combining what we have proved so far, we get the following.
\begin{corollary}
It some octagonal parameter $s$ has an acute crosscut,
then some octagonal parameter $u<1/2$ has an
egregrious crosscut.
\end{corollary}

\begin{lemma}
\label{egg3}
If $s<1/2$ has an egregious crosscut then
$t=R(s)<1/2$ and $t$ has an egregious 
crosscut.
\end{lemma}

\startproof
When $s<1/2$ and $R(s)>1/2$, it
follows from Theorem \ref{renorm} (or from
a direct calculation) that
there is a large octagon
which separates the top edge of $X_s^0$ from the
bottom edge of $X_s^0$. See Figure 7.7.

\begin{center}
\resizebox{!}{1.5in}{\includegraphics{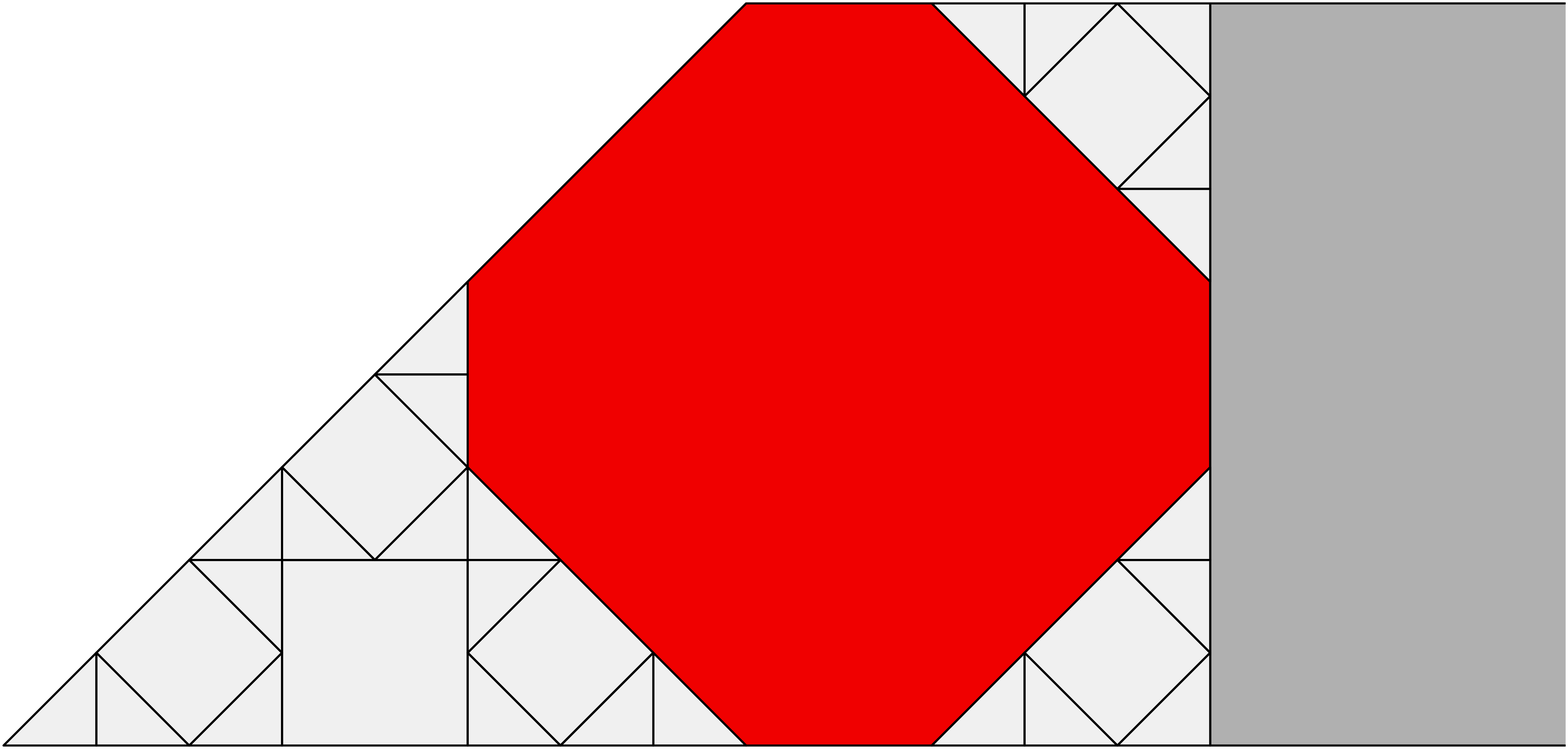}}
\newline
{\bf Figure 7.7:\/} The large octagon.
\end{center}

Now we know that $R(s)<1/2$.
In this case, $C$ has to cross the set
\begin{equation}
Z_s'=R_D(Z_s^0) \cap  X_s.
\end{equation}
$C$ starts out in the top edge of
$Z_s'$ and exits through the bottom edge.

\begin{center}
\resizebox{!}{1.8in}{\includegraphics{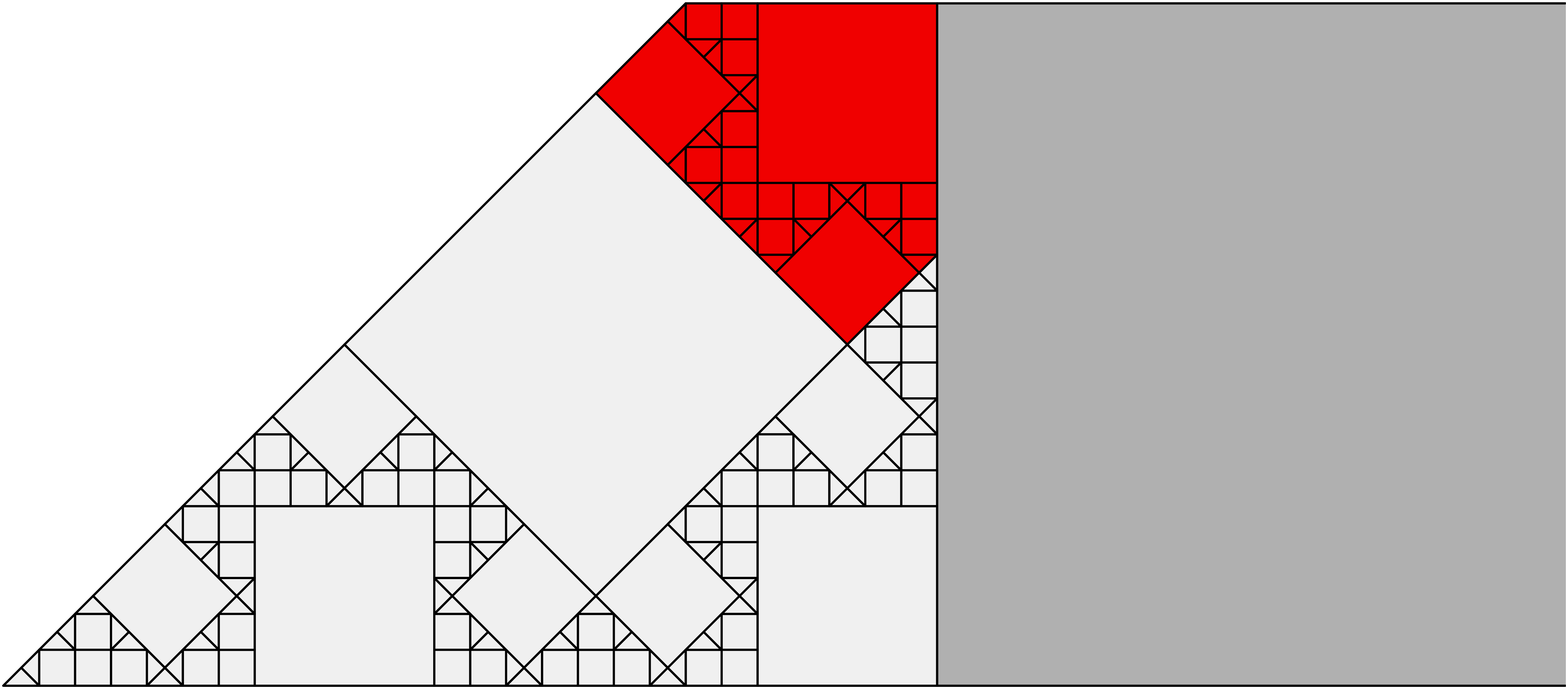}}
\newline
{\bf Figure 7.8:\/} $Z_s'$ (red) for $s=19/52$.
\end{center}

The set
\begin{equation}
C'=\phi_s^{-1} \circ R_D(C)
\end{equation}
connects a point on the bottom edge of $X_t$ with a
point that lies above the line $H$ of symmetry.
In this case (as in the proof of
Lemma \ref{egg2}) $C' \cup R_H(C')$ is an
egregious crosscut for $t$. 
\endproof

Iterating Lemma \ref{egg3}, we see that
$R^n(s)<1/2$ for all $n$ when $s$ has
an egregious crosscut.  So, $s$ cannot
be octagonal in this case!
If we assume that Statement 3 of the Main 
Theorem is false, then we have reached a contradiction.

\newpage

\section{References}

[{\bf AG\/}] A. Goetz and G. Poggiaspalla, {\it Rotations by $\pi/7$\/}, Nonlinearity {\bf 17\/}
(2004) no. 5 1787-1802
\newline
\newline
[{\bf AKT\/}] R. Adler, B. Kitchens, and C. Tresser,
{\it Dynamics of non-ergodic piecewise affine maps of the torus\/},
Ergodic Theory Dyn. Syst {\bf 21\/} (2001) no. 4 959-999
\newline
\newline
[{\bf BKS\/}] T. Bedford, M. Keane, and C. Series, eds.,
{\it Ergodic Theory, Symbolic Dynamics, and Hyperbolic Spaces\/}, Oxford University Press, Oxford (1991).
\newline
\newline
[{\bf H\/}] H. Haller, {\it Rectangle Exchange Transformations\/}, Monatsh Math. {\bf 91\/}
(1985) 215-232
\newline
\newline
[{\bf Hoo\/}] W. Patrick Hooper, {\it Renormalization of Polygon Exchage Maps arising from Corner Percolation\/} Invent. Math. 2012.
\newline
\newline
[{\bf K\/}], M. Keane, {\it Interval Exchange Transformations\/}, Math Z. {\bf 141\/}, 25-31 (1975).
\newline
\newline
[{\bf LKV\/}] J. H. Lowenstein, K. L. Koupsov, F. Vivaldi, {\it Recursive Tiling and Geometry of piecewise rotations by $\pi/7$\/}, nonlinearity {\bf 17\/} (2004) no. 2.
[{\bf Low\/}] J. H. Lowenstein, {\it Aperiodic orbits of piecewise rational
rotations of convex polygons with recursive tiling\/}, Dyn. Syst. {\bf 22\/}
(2007) no. 1 25-63
\newline
\newline
[{\bf R\/}] G. Rauzy, {it Exchanges d'intervalles et transformations induites\/},
Acta. Arith. {\bf 34\/} 315-328 (1979)
\newline
\newline
[{\bf S0\/}] R.E. Schwartz {\it The Octagonal Pet I: Hyperbolic Symmetry and Renormalization\/},
preprint (2012)
\newline
\newline
[{\bf S1\/}] R.E. Schwartz {\it Outer Billiards, Quarter Turn Compositions, and Polytope Exchange 
Transformations\/}, preprint (2011)
\newline
\newline
[{\bf S2\/}] R. E. Schwartz, {\it Outer Billiards on Kites\/},
Annals of Math Studies {\bf 171\/}, Princeton University Press (2009)
\newline
\newline
[{\bf S3\/}] R. E. Schwartz, {\it Outer Billiards on the Penrose Kite:
Compactification and Renormalization\/}, Journal of Modern Dynamics, 2012.
\newline
\newline
[{\bf T\/}] S. Tabachnikov, {\it Billiards\/}, Soci\'{e}t\'{e} Math\'{e}matique de France, 
``Panoramas et Syntheses'' 1, 1995
\newline
\newline
[{\bf V1\/}] W. Veech, {\it The metric theory of interval exchange transformations I: Generic spectral properties\/}, Amer. Journal of Math. {\bf 106\/}, 1331-1359 (1984)
\newline
\newline
[{\bf V2\/}] W. Veech, {\it The metric theory of interval exchange transformations II: Approximation by Primitive Interval Exchanges\/} Amer. Journal of Math {\bf 106\/}, 1361-1387 (1984)
\newline
\newline
[{\bf VL\/}] F. Vivaldi and J. H. Lowenstein, {|it Arithmetical properties of a family
of irrational piecewise rotations\/}, {\it Nonlinearity\/} {\bf 19\/}:1069--1097 (2007).
\newline
\newline
[{\bf Y\/}] J.-C. Yoccoz, {\small {\it Continued Fraction Algorithms for Interval Exchange Maps: An Introduction\/}\/}, Frontiers in Number Theory, Physics, and Geometry Vol 1, P. Cartier, B. Julia, P. Moussa, P. Vanhove (editors) Springer-Verlag 4030437 (2006)
\newline
\newline
[{\bf Z\/}] A. Zorich, {\it Flat Surfaces\/}, Frontiers in Number Theory, Physics, and Geometry Vol 1, P. Cartier, B. Julia, P. Moussa, P. Vanhove (editors) Springer-Verlag 4030437 (2006)

\end{document}